\documentclass[11pt]{article}
\usepackage{a4}
\usepackage{latexsym}
\usepackage{amsfonts}
\usepackage[noxcolor]{pstricks}
\usepackage{amssymb}
\usepackage{amscd}
\usepackage{amsmath}
\usepackage{epsfig}
\usepackage{xypic}
\newtheorem{theorem}{Theorem}%[section]
\newtheorem{lemma}[theorem]{Lemma}
\newtheorem{proposition}[theorem]{Proposition}
\newtheorem{corollary}[theorem]{Corollary}
\newtheorem{definition}[theorem]{Definition\rm}
\newtheorem{remark}{Remark\/}
\newtheorem{example}{Example\/}
\newtheorem{conjecture}[theorem]{Conjecture\/}
\newtheorem{notation}{Notation\/}
\begin{document}
\title{\textbf{Zeta functions and monodromy for surfaces that are general for a toric idealistic cluster}}
\author{Ann Lemahieu and Willem Veys\footnote{Ann Lemahieu, Willem Veys, K.U.Leuven, Departement Wiskunde,
Celestij\-nenlaan 200B, B-3001 Leuven, Belgium, email:
lemahieu@mathematik.uni-kl.de, wim.veys@wis.kuleuven.ac.be. The
research was partially supported by the Fund of Scientific Research
- Flanders (G.0318.06) and MEC PN I+D+I MTM2007-64704.}
\date{}} \maketitle {\footnotesize \emph{\textbf{Abstract.---} In
this article we consider surfaces that are general with respect to a
$3$-dimensional toric idealistic cluster. In particular, this means
that blowing up a toric constellation provides an embedded
resolution of singularities for these surfaces. First we give a
formula for the topological zeta function directly in terms of the
cluster. Then we study the eigenvalues of monodromy. In particular,
we derive a useful criterion to be an eigenvalue. In a third part we
prove the monodromy and the holomorphy conjecture for
these surfaces.}} %is to
%prove the monodromy and the holomorphy conjecture for % that we study
%predicts a relation between the poles of the topological zeta
%function associated to a hypersurface and the eigenvalues of
%monodromy on the hypersurface. In this article In order to do that, we compute the topological
%zeta function and we analyse the eigenvalues of monodromy.  We then
%provide a proof for the
\\ \\
 ${}$ \begin{center}
\textsc{Contents}
\end{center} ${}$\\
$\begin{array}{lr}
\mbox{1. Introduction} & \qquad 1 \\
\mbox{2. Toric clusters} & \qquad 4 \\
\mbox{3. Conjectures} & \qquad 7\\
\mbox{4. Computation of the topological zeta function} & \qquad 10 \\
\mbox{5. Analysis of $\chi(E_i^{\circ})$} & \qquad 14\\
\mbox{6. Determination of the sign of $\chi(E_i^{\circ})$} & \qquad 22\\
\mbox{7. The monodromy conjecture for candidate poles of order $1$} & \qquad 30\\
\mbox{8. The monodromy conjecture for candidate poles of order $2$
or $3$}
& \qquad 34\\
\mbox{9. The holomorphy conjecture} & \qquad 40
%\\ \mbox{10. Remarks} & \qquad
\end{array}$
%\\ Let $E_j, 1 \leq
%j \leq r$, be the irreducible exceptional components created by this
%resolution. For such surfaces we give an algorithm to compute the
%topological zeta function directly out of the toric cluster. In
%particular, we obtain an explicit formula for the topological Euler
%characteristic of the spaces $E_j^{\circ}$, where $E_j^{\circ}:= E_j
%\setminus (\cup_{i \in \{1,\cdots,r\} \setminus \{j\}} E_i)$. By
%combinatorial arguments we determine the sign of
%$\chi(E_j^{\circ})$. This classification permits us to prove the
%monodromy conjecture in this context.}}
\\ \\
 ${}$ \begin{center}
\textsc{1. Introduction}
\end{center} ${}$\\
In \cite{weil}, Weil introduced some zeta functions
$\mathcal{Z}(K,f)$ that are integrals over a $p$-adic field $K$ and
that are associated to a polynomial $f(\underline{x}) \in
K[\underline{x}]$. Using embedded resolution of singularities, Igusa
showed that these zeta functions are rational and he studied their
poles (see \cite{Igusa75} and \cite{Igusa78}). One can define the
analogous integrals over $K=\mathbb{R}$ or $\mathbb{C}$. Also these
zeta functions are rational (see for example \cite{atiyah} and
\cite{bernstein69}) and it is known that their poles are contained
in the set of roots - and roots shifted by a negative integer - of
the Bernstein polynomial $b_f$. By Malgrange (\cite{malgrange}), if
$\alpha$ is a root of $b_f$, then $e^{2 \pi i \alpha}$ is an
eigenvalue of the local monodromy of $f$ at some point of
$f^{-1}(0)$. So when $K= \mathbb{R}$ or $\mathbb{C}$, then the poles
of the zeta function induce eigenvalues of the local monodromy. This
result was a motivation to study this relation at the $p$-adic side.
The study of concrete examples made it natural to propose the
following conjecture.
\\ \\
\textbf{Monodromy conjecture.} (\cite{Igusa88})  Let $F \subset
\mathbb{C}$ be a number field and $f \in F[\underline{x}]$. For
almost all $p$-adic completions $K$ of $F$, if $s_0$ is a pole of
$\mathcal{Z}(K,f)$, then $e^{2 \pi i Re(s_0)}$ is an eigenvalue of
the local monodromy of $f$ at some point of the hypersurface $f=0$.
\\ \\
Loeser verified this conjecture for plane curves (see
\cite{loeser1}). He also gave a proof for a class of polynomials in
higher dimensions; the polynomial should be nondegenerate with
respect to its Newton polyhedron and should satify some numerical
conditions (\cite{loeser2} and Section $3$).
\\ \indent When Denef and Loeser introduced the topological zeta function in
$1992$ in \cite{DenefLoeser1}, an analogous version of the monodromy
conjecture arose. This monodromy conjecture relates the poles of the
topological zeta function $Z_{top,f}$ associated to a polynomial
function or a germ of a holomorphic function $f$ with the
eigenvalues of monodromy of the hypersurface $f=0$.
\\ \\ \textbf{Monodromy conjecture.} If $s_0$ is a pole
of $Z_{top,f}$, then $e^{2\pi i s_0}$ is an eigenvalue of the local
mo\-no\-dromy of $f$ at some point of the hypersurface $f=0$.
\\ \\
By the original definition of the topological zeta function, it
follows that the monodromy conjecture for the Igusa zeta function
implies the monodromy conjecture for the topological zeta function.
Artal Bartolo, Cassou-Nogu\`{e}s, Luengo and Melle Hern\'{a}ndez
proved the monodromy conjecture for some surface singularities, such
as the superisolated ones (see \cite{luengo}), and for
quasi-ordinary polynomials in \cite{luengo2}. The second author
provided results in \cite{Veysconfigurations}, \cite{Veys2} and
\cite{Veys3}, and together with Rodrigues in \cite{rove03}. In
\cite{CRAS}, the authors consider the same context as in this paper
but they had to impose a restricting condition on the surfaces. Via
geometrical arguments they showed that the monodromy conjecture
holds for candidate poles of the topological zeta function of order
$1$ that are poles.
\\ \\
There are more conjectures relating the poles of the topological
zeta function (Igusa zeta function) and the eigenvalues of
monodromy. There exist the rational functions $Z_{top,f}^{(r)}$ ($r
\in \mathbb{Z}_{> 0}$) that are variants of the topological zeta
function and that play a role in the holomorphy conjecture,
which was stated by Denef. \\
\\ \noindent \textbf{Holomorphy conjecture.} (\cite{denef91}) If $r \in \mathbb{Z}_{> 0}$ does not divide
the order of any eigenvalue of the local monodromy of $f$ at any
point of $f^{-1}\{0\}$, then $Z_{top,f}^{(r)}$ is holomorphic on
$\mathbb{C}$.
\\ \\
Originally the holomorphy conjecture was formulated for the Igusa
zeta function. We refer to \cite{denef91} for the inspiration. Denef
showed that the conjecture is true for the relative invariants of a
few prehomogeneous vector spaces. The second author proved the
conjecture for plane curves (see \cite{Veys93A}) and together with
Rodrigues for homogeneous polynomials (see \cite{rove01}).
\\ \\
Although the monodromy conjecture and/or holomorphy conjecture has
been proven for these kinds of singularities, one did not get a
better understanding of the deep reason why the conjectures hold for
them. Until now, the attempts are thus restricted to prove the
conjecture for classes of singularities. \\ \indent This article
deals with the class of surfaces that are general with respect to a
$3$-dimensional toric idealistic cluster. This implies that we work
with surfaces for which there exists an embedded resolution of
singularities by blowing up in points that are orbits for the action
of the torus, i.e. in a toric constellation. We refer to Section $2$
for a recap about clusters and in Section $3$ we explain the objects
that play the main role in the conjecture. In Section $4$ we show
how the topological zeta function can be computed directly in terms
of the toric cluster for the surfaces that we consider. We use the
embedded resolution provided by the blowing up of the constellation.
Let $\pi: Z \rightarrow \mathbb{C}^3$ be that resolution of such a
surface $f=0$ and let $E_j, j \in S$, be the irreducible components
obtained by this resolution of which $E_1,\cdots,E_r$ are the
exceptional ones. We will denote $E_j^{\circ}:= E_j \setminus
(\cup_{i \in S \setminus \{j\}} E_i)$, for $j \in S$. We write $N_j$
and $\nu_j - 1$ for the multiplicities of $E_j$ in the divisor on
$Z$ of $f \circ \pi$ and $\pi^\ast (dx \wedge dy \wedge dz)$,
respectively. The numbers $-\nu_j/N_j, j \in S$, form a complete
list of candidate poles of $Z_{top,f}$.
\\ \indent We compute in particular the Euler characteristic
of the spaces $E_j^{\circ}, 1 \leq j \leq r$, in terms of the
cluster. They show up in A'Campo's formula for the eigenvalues of
mo\-no\-dromy and they are very relevant for the monodromy
conjecture. In a fifth section we analyse these Euler
characteristics. Our goal is to determine when these numbers are
less than or equal to $0$. A geometric argument will show that we
can reduce this job to the investigation of a finite number of
families of constellations. We complete Section $5$ with
combinatorial preparations. These make it possible to determine the
sign of the Euler characteristics that we are looking for. We carry
this out in Section $6$. We then
prove the following result. \\ \\
\textbf{Theorem.} If $\chi(E_j^{\circ}) > 0$, then $e^{-2\pi i
\frac{\nu_j}{N_j}}$ is an eigenvalue of monodromy of $f$.
\\ \\
Using this result, we prove in Section $7$ the monodromy conjecture
for candidate poles of order $1$ that are poles and in Section $8$
the monodromy conjecture for candidate poles of order $2$ or
$3$ that are poles. Hence, we obtain: \\ \\
\textbf{Theorem.} Let $f$ be a germ of a polynomial map that is
general with respect to a $3$-dimensional toric idealistic cluster.
If $s_0$ is a pole of $Z_{top,f}$, then $e^{2\pi i s_0}$ is an
eigenvalue of mo\-no\-dromy of $f$ at some point of the hypersurface
$f=0$.
\\ \\
In Section $9$ we prove the holomorphy conjecture for these
surfaces.
\\ \\
\textbf{Theorem.} Let $f$ be a germ of a polynomial map that is
general with respect to a $3$-dimensional toric idealistic cluster.
If $r \in \mathbb{Z}_{> 0}$ does not divide the order of any
eigenvalue of the local monodromy of $f$ at any point of $f=0$, then
$Z_{top,f}^{(r)}$ is holomorphic on $\mathbb{C}$.
%\\ \\
%\noindent Apart from proving these conjectures for new cases, we get
%some very explicit results about the eigenvalues of monodromy. In
%the last section we emphasise some results that are particular for
%this context and we point to some details that might be interesting
%to investigate in other contexts too.
\\ \\  ${}$
\begin{center} \textsc{2. Toric clusters}
\end{center} ${}$\\
In this section we introduce the terminology of infinitely near
points, (toric) clusters etc. according to \cite{clusters}. %We also
%comment some known results about clusters that are of interest for
%the mathematics that will be developed further on in this work.
We would like to refer to \cite{clusters} for some histo\-ri\-cal
notes on clusters. See also \cite{casas1}, \cite{casas2}, \cite{EC},
\cite{lipman1}, \cite{lipman2}, \cite{lipman3} and \cite{Zariski}
for more details on the theory of clusters.
\\ \\
\textbf{\emph{$2.1.$ Clusters.}}--- Let $X$ be a nonsingular variety
of dimension $d \geq 2$ and let $Z$ be a variety obtained from $X$
by a finite succession of point blowing-ups. A point $Q \in Z$ is
said to be \emph{infinitely near} to a point $P \in X$ if $P$ is in
the image of $Q$; we write $Q \geq P$. A \emph{constellation} is a
finite sequence $\mathcal{C}:=\{Q_1,Q_2,\cdots,Q_{r}\}$ of
infinitely near points of $X$ with $Q_1 \in X=:X_0$ and each
$Q_{j+1}$ is a point on the variety $X_j$ obtained by blowing up
$Q_{j}$ in
$X_{j-1}$, $j \in \{1,\cdots,r-1\}$. The variety $X(\mathcal{C}):=X_r$
obtained by blowing up $Q_r$ in $X_{r-1}$ is called the \emph{sky}.\\
\indent The relation `$\geq$' gives rise to a partial ordering on
the points of a constellation. In the case that they are totally
ordered, so $Q_{r} \geq \cdots \geq Q_1$, the constellation
$\mathcal{C}$ is called a \emph{chain}. For every $Q_j$ in
$\mathcal{C}$, the subsequence $\mathcal{C}^j:=\{Q_i \mid Q_j \geq
Q_i\}$ of $\mathcal{C}$ is a chain. The integer $l(Q_j):=\#
\mathcal{C}^j - 1$ is called the \emph{level} of $Q_j$. In
particular $Q_1$ has level $0$. If no other point of $\mathcal{C}$
has level $0$ then $Q_1$ is called the \emph{origin} of
$\mathcal{C}$. We will always work with constellations that have an
origin and we will also denote the origin of the constellation by
$o$. If $Q_j \geq Q_i$ and $l(Q_j)=l(Q_i)+1$, we will write $Q_j
\succ Q_i$ or $j \succ i$.
\\ \indent For each $Q_i \in \mathcal{C}$, denote the exceptional divisor of the blowing-up in $Q_i$
by $E_i$, as well as its strict transform at some intermediate stage
(including the final stage) $X_j$, $i \leq j \leq r$. The total
transform at some intermediate stage (including the final stage)
will be denoted by $E_i^*$. If $Q_j \in E_i$, then one says that
$Q_j$ is \emph{proximate} to $Q_i$. This will be denoted as $Q_j
\rightarrow Q_i$ or $j \rightarrow i$. As $E_i = E_i^* - \sum_{j
\rightarrow i} E_j^*$, it follows that also
$\{E_1^*,\cdots,E_{r}^*\}$ is a basis of the group of divisors
with exceptional support $\oplus_{j=1}^{r} \mathbb{Z}E_j$. \\
\indent A pair $\mathcal{A}:=(\mathcal{C},\underline{m})$ consisting
of a constellation $\mathcal{C}:=\{Q_1,\cdots,Q_{r}\}$ and a
sequence $\underline{m}:=(m_1,\cdots,m_{r})$ of nonnegative integers
is called a \emph{cluster}. One calls $m_j$ the \emph{weight} or
\emph{multiplicity} of $Q_j$ in the cluster and we write
$D(\mathcal{A}):=\sum_{j=1}^{r}m_j E_j^*$. Introducing the numbers
$v_j$, $1 \leq j \leq r$, by setting $m_j:=v_j - \sum_{j \rightarrow
i}v_i$, allows us to write also $D(\mathcal{A})=\sum_{j=1}^{r}v_j
E_j$. The idea of clusters is to express that a system of
hypersurfaces is passing through the points of the constellation
with (at least) the given multiplicities. This explains why we are
interested in the ideals
\[I(v_1,\cdots,v_r)=\{g \in \mathcal{O}_{X,o} \mbox{ $|$ } \nu_j(g)
\geq v_j, 1 \leq j \leq r\}\] for a point $o \in X$, discrete
valuations $\nu_1,\cdots,\nu_r$ on $k(X)$ and $(v_1,\cdots,v_r) \in
\mathbb{Z}^r$. If we want that these ideals principalise by blowing
up the points of the constellation, we require the ideals to be
finitely supported. Formally, an ideal $I$ in $\mathcal{O}_{X,o}$ is
called \emph{finitely supported} if $I$ is primary for the maximal
ideal $\mathsf{m}$ of $\mathcal{O}_{X,o}$ - so supported at the
closed point - and if there exists a constellation $\mathcal{C}$ of
infinitely near points of $X$ such that
$I\mathcal{O}_{X(\mathcal{C})}$ is an invertible sheaf.
\\ \indent On the other hand, given a finitely supported
ideal $I$, one can associate a cluster to it. Let
$\mathcal{C}_I=:\{Q_1,\cdots,Q_{r}\}$ be the constellation of base
points of $I$, i.e. the minimal constellation $\mathcal{C}$ such
that
$I\mathcal{O}_{X(\mathcal{C})}$ is an invertible sheaf. %Let $I$ be
%an ideal in a local ring $(\mathcal{O}_{X,Q},\mathsf{m})$, then the
%\emph{order} ord$_Q I$ \emph{of $Q$ at $I$} is defined as max$\{n
%\mbox{ $|$ } I \subset \mathsf{m}^n\}$. One associates a cluster
%$\mathcal{A}_I:=(\mathcal{C}_I,\underline{m})$ to $I$ as follows:
%for $j \in \{0,\cdots,r-1\}$, let $Q_{j}^-$ be the point in the
%constellation for which $Q_j \succ Q_{j}^-$ and denote its weight by
%$m_{j}^-$. The weight $m_j$ of $Q_j$ is defined by induction by
%setting:
%\begin{enumerate}
%\item $m_0:=\mbox{ord}_{Q_0}I_{Q_0}$ with $I_{Q_0}:=I$;
%\item $m_j := \mbox{ord}_{Q_j}I_{Q_j}$ with
%$I_{Q_j}:=(x)^{-m_{j}^-}I_{Q_{j}^-}\mathcal{O}_{X,Q_j}$,
%\end{enumerate}
%where $x$ is a generator of the principal ideal
%$\mathsf{m}_{Q_{j}^-}\mathcal{O}_{X,Q_j}$.
Let $m_j$ be the order of the point $Q_j, 1 \leq j \leq r$ in the
strict transform of the ideal $I$ in $\mathcal{O}_{X_j,Q_j}$. Then
the ideal sheaf $I\mathcal{O}_{X(\mathcal{C}_I)}$ is associated to
$-D(\mathcal{A}_I):=\sum_{j=1}^r m_j E_j^*$.
\\ \indent
If $\mathcal{C}$ is a constellation with origin at $Q_1$, the
cluster $\mathcal{A}:=(\mathcal{C},\underline{m})$ is called
\emph{idealistic} if there exists a finitely supported ideal $I$ in
$\mathcal{O}_{X,Q_1}$ such that $I\mathcal{O}_{X(\mathcal{C})}$ is
the ideal sheaf associated to $-D(\mathcal{A})$. %We call
%\emph{galaxy} of $\mathcal{C}$ the set of idealistic clusters on
%$\mathcal{C}$.
For an idealistic cluster $\mathcal{A}$, Lipman proved that there
exists a unique finitely supported complete ideal $I_{\mathcal{A}}$
such that $I_{\mathcal{A}}
\mathcal{O}_{X(\mathcal{C})}=\mathcal{O}_{X(\mathcal{C})}(-D(\mathcal{A}))$,
namely that given by the direct image of
$\mathcal{O}_{X(\mathcal{C})}(-D(\mathcal{A}))$ in $X$, see
\cite{lipman2}. %Note that the galaxy of a constellation is a
%semigroup. Indeed, if $\mathcal{A}_1:=(\mathcal{C},\underline{m_1})$
%and $\mathcal{A}_2:=(\mathcal{C},\underline{m_2})$ are idealistic
%clusters, then also
%$\mathcal{A}:=(\mathcal{C},\underline{m_1}+\underline{m_2})$ is
%idealistic. One has $I_{\mathcal{A}}=I_{\mathcal{A}_1} *
%I_{\mathcal{A}_2}$, where $*$ denotes the completion of the product
%of ideals.
\\ \\
\textbf{\emph{$2.2$. Toric clusters in $\mathbb{C}^3$.}}--- From now
on suppose that $X$ is the affine toric variety $\mathbb{C}^3$.
%Let $M$ be a $d$-dimensional lattice $(d \geq 2)$, $N$ its dual
%lattice and $\sigma$ a regular cone in $N_{\mathbb{R}}$. Consider
%the smooth affine toric variety $X:=$ Spec
%$\mathbb{C}[\check{\sigma}\cap M]$.
A $3$-dimensional \emph{toric constellation} of infinitely near
points with origin $Q_1$ is a constellation
$\mathcal{C}:=\{Q_1,Q_2,\cdots,Q_{r}\}$ such that each $Q_j$ is a
$0$-dimensional orbit in the toric variety $X_j$ obtained by blowing
up $Q_{j-1}$ in $X_{j-1}$, $2 \leq j \leq r$. Blowing up in orbits
of smooth varieties corresponds to making star subdivisions of the
fan corresponding to the variety (see for example \cite{oda}). In
this way each blowing-up in a $0$-dimensional orbit induces the
creation of three cones of dimension $3$ and thus of three new
$0$-dimensional orbits. Hence, the choice of a point $Q_i$ in a
toric chain is equivalent to the choice of an integer $a_i \in
\{1,2,3\}$, which determines a $3$-dimensional cone in the fan. A
tree with a root such that each vertex has at most three following
adjacent vertices is called a \emph{$3$-nary tree}. The above
observation shows that there is a natural bijection between the set
of $3$-dimensional toric constellations with origin and the set of
finite $3$-nary trees with a root, with the edges labeled with
positive integers not greater than $3$, such that two edges with the
same source have different labels.
\\ \indent A cluster $\mathcal{A}:=(\mathcal{C},\underline{m})$ is called
\emph{toric} if the constellation $\mathcal{C}$ is toric. As we want
the finitely supported ideals to be supported in the $0$-dimensional
orbit, they should be invariant under the action of the torus and
thus be monomial. \\ \\ \textbf{\emph{$2.3.$ Properties.}}--- In
this subsection we recall some properties about clusters, in
particular about toric clusters.
\\ We will prove the monodromy and holomorphy conjectures for the class of surfaces
for which the following theorem holds (see \cite{clusters}).
\begin{theorem} The canonical map from the sky of the
constellation of base points of a finitely supported ideal $I$ to
$X$ is an embedded resolution of the subvariety of $(X,o)$ defined
by a general enough element in $I$.
\end{theorem}
We will call these `general enough' elements \emph{general for $I$}
or \emph{general for \emph{$\mathcal{C}_I$}}. %Notice that not all
%surfaces for which the blowing up of a toric constellation gives an
%embedded resolution are general for some finitely supported ideal.
%Take for example $x^5+y^5+z^5+xyz=0$.
\\ \\In the case of toric clusters, there exists a combinatorial characterisation
for the idealistic clusters. Fix a point $Q_i$ in a toric
$3$-dimensional constellation $\mathcal{C}$ and some integers $a, b$
such that $a,b \in \{1,2,3\}$ and $a \neq b$. For $s,t \in
\mathbb{Z}_{\geq 0}$, let $Q_i(a^s,b^t)$ be the terminal point of
the chain with origin $Q_i$ coded by $(a,\cdots,a,b,\cdots,b)$ where
$a$ appears $s$ times and $b$ appears $t$ times. If $t=0$, it is
denoted by $Q_i(a^s)$. The point $Q_i(a^s,b^t)$ may not belong to
$\mathcal{C}$. A point $Q_j \in \mathcal{C}$ that is infinitely near
to $Q_i$ is said to be \emph{linearly proximate} to $Q_i$, if
$Q_j=Q_i(a,b^t)$, with $a,b$ and $t$ as above. We denote this
relation by $Q_j \twoheadrightarrow Q_i$ or $j \twoheadrightarrow
i$. Then we have that $Q_j$ is linearly proximate to $Q_i$ if and
only if there exists a $1$-dimensional orbit $\mathsf{l}$ in $B_i$
such that $Q_j$ belongs to the strict transform of the closure of
$\mathsf{l}$ in $E_i$. This explains the terminology. Denote
$M_{Q_i}(a,b):=\sum_{t \geq 0}m_{{Q_i}(a,b^t)}$. Campillo,
Gonzalez-Sprinberg and Lejeune-Jalabert show the following.
\begin{enumerate}
\item A toric cluster $\mathcal{A}=(\mathcal{C},\underline{m})$ is idealistic if and only if
for each point $Q_i$ of the constellation $\mathcal{C}$ and for each
pair of integers $a$ and $b$ such that $a, b \in \{1,2,3\}$ and $a
\neq b$, the following inequality is satisfied:
\[M_{Q_i}(a,b) + M_{Q_i}(b,a) \leq m_{Q_i}.\]
These inequalities are called the \emph{linear proximity
inequalities}. %They generalise in dimension bigger than $2$ for
%toric clusters the inequalities that Enriques and Chisini obtained
%in dimension $2$.
\item Let $\mathcal{A}=(\mathcal{C},\underline{m})$ be a
$3$-dimensional toric idealistic cluster with associated divisor
$D(\mathcal{A})=\sum_{j=1}^r m_j E_j^*=\sum_{j=1}^r v_j E_j$ and let
$\nu_1,\cdots,\nu_r$ be the induced discrete valuations. Such a
valuation is called Rees for the ideal
$I(\underline{v}):=I(v_1,\cdots,v_r)$ if it is a valuation induced
by an irreducible component of the exceptional divisor of the
normalised blowing-up $\overline{Bl_{I(\underline{v})} X}$ of
$I(\underline{v})$. Then \vspace*{-0.1cm}
\begin{eqnarray} \forall
Q_i \in \mathcal{C}: m_i^{2} \geq \sum_{j \rightarrow i} m_j^2
\qquad \qquad \mbox{ and}
\end{eqnarray}
\vspace*{-0.5cm}
\begin{eqnarray} \nu_i \mbox{ is Rees for $I(\underline{v})$ if and only if } m_i^{2}
> \sum_{j \rightarrow i} m_j^{2}.
\end{eqnarray}
\end{enumerate}
To a monomial ideal $I$ one can associate a Newton polyhedron
$\mathcal{N}_I$. It is the union of the compact faces of the convex
hull of $m + \mathbb{Z}^3_{\geq 0}$ as $m$ runs through the set of
exponents of monomials in $I$. We refer to \cite{kempf} for the
proofs of the following properties.
\begin{enumerate}
\item The facets of $\mathcal{N}_I$ correspond with the Rees
valuations of $I$. \item A monomial ideal is complete if and only if
it contains every monomial whose exponent is a point of
$\mathcal{N}_I + \mathbb{Z}^3$.
\end{enumerate}
%${}$\\  \textbf{\emph{$2.4.$ Example.}}---
\newpage
\begin{example}
\emph{${}$ \\${}$
\begin{tabular}{p{4.5cm}p{8cm}}
\begin{picture}(35,10)(-8,-2)
\put(20,-25){\circle*{1}}
\put(30,-15){\circle*{1}}\put(10,-15){\circle*{1}}
\put(20,-25){\line(1,1){10}} \put(20,-25){\line(-1,1){10}}
\put(18,-29){\footnotesize{$Q_1$}}
\put(32,-15){\footnotesize{$Q_3$}} \put(4,-15){\footnotesize{$Q_2$}}
\put(28,-20){\scriptsize{\emph{3}}}
\put(11,-20){\scriptsize{\emph{1}}} \put(10,-15){\line(1,1){10}}
\put(10,-15){\line(-1,1){10}}\put(20,-5){\circle*{1}}\put(0,-5){\circle*{1}}
\put(22,-5){\footnotesize{$Q_5$}}  \put(-6,-5){\footnotesize{$Q_4$}}
\put(18,-10){\scriptsize{\emph{2}}}
\put(1,-10){\scriptsize{\emph{1}}}
\end{picture}
& Suppose $d=3$ and $\mathcal{C}$ is the constellation pictured at
the left. It represents the following resolution process: by blowing
up in the origin $Q_1$ we get an exceptional variety $E_1 \cong
\mathbb{P}^2$. In $E_1$ there are two points in which we blow up,
namely $Q_2$ and $Q_3$. The labels indicate in which affine chart
the points of the constellation are created.
\end{tabular}
\\ For example the point
$Q_2$ is the origin of the affine chart induced by the edge going
out of $Q_1$ with label $1$. After blowing up in $Q_2$ we get an
exceptional variety $E_2 \cong \mathbb{P}^2$, where again we blow up
in two points.
%\end{example}\hfill $\square$
\\ \\
%\vspace*{-1cm} \[\]
The induced valuations are represented by the following vectors in
the lattice $\mathbb{N}^3$:
\begin{eqnarray*}
\nu_1 \leftrightarrow (1,1,1) \quad \nu_2 \leftrightarrow (1,2,2)
\quad \nu_3 \leftrightarrow (2,2,1) \quad \nu_4 \leftrightarrow
(1,3,3) \quad \nu_5 \leftrightarrow (2,3,4).
\end{eqnarray*}
Consider the following multiplicities for the points of this
constellation: \[(m_1,m_2,m_3,m_4,m_5)=(3,2,1,1,1) \mbox{ or }
(v_1,v_2,v_3,v_4,v_5)=(3,5,4,6,9).\] Saying that a monomial $x^a y^b
z^c$ passes through $Q_j$ is saying that $\nu_j(x^a y^b z^c) \geq
v_j$, for $1 \leq j \leq 5$. The
induced hyperplanes define a Newton polyhedron. %We picture the Newton
%polyhedron  $\mathcal{N} + \check{\sigma}$.
%\begin{center}
%\epsfig{figure=Newtonpolytoop.eps, width= 7cm, height= 6cm}
%\end{center}
%\vspace*{-2cm}
Now let $I_{\mathcal{A}}$ be the ideal generated by
the monomials whose exponents are in this Newton polyhedron. We find
\[I_{\mathcal{A}}=(x^6,y^3,z^4,x^3y,x^2y^2,yz^2,y^2z,x^3z,xz^2,xyz).\]
The blowing-up of the constellation gives an embedded resolution for
a general element of $I_{\mathcal{A}}$, such as for example
$h(x,y,z):=x^6+y^3+z^4+x^3y+x^2y^2+yz^2+y^2z+x^3z+xz^2-xyz$.} \hfill
$\square$ ${}$\end{example} \vspace*{-0.5cm} ${}$\\ \\ ${}$
%\newpage
\begin{center}
\textsc{3. Conjectures}
\end{center} ${}$\\
Let $f$ be a complex polynomial in $d$ variables and let $\pi : Z
\rightarrow \mathbb{C}^d$ be an embedded resolution of singularities
of $f^{-1} \{ 0 \}$. We write $E_j, j \in S$, for the irreducible
components of $\pi^{-1}(f^{-1}\{ 0 \})$ and we denote by $N_j$ and
by $\nu_j - 1$ the multiplicities of $E_j$ in the divisor on $Z$ of
$f \circ \pi$ and $\pi^\ast (dx_1 \wedge \ldots \wedge dx_d)$,
respectively. The couples $(\nu_j, N_j), j \in S$, are called the
numerical data of the embedded resolution $(Z,\pi)$. We denote also
$E_j^{\circ}:= E_j \setminus (\cup_{i \in S \setminus \{j\}} E_i)$,
for $j \in S$. Let the $E_j$, $j \in J:=\{1,\cdots,r\} \subset S$,
be the exceptional irreducible components of $\pi^{-1}(\{ 0 \})$.
\\ \\
\textbf{\emph{$3.1.$ Monodromy.}}--- We assume that $f(b)=0$. Take
$\epsilon
> 0$ small enough such that the open ball $B_{\epsilon}$ with radius $\epsilon$ around
$b$ in $\mathbb{C}^d$ intersects the fibre $f^{-1}(0)$
transversally. Then choose $\epsilon \gg \eta > 0$ such that for $t$
in the disc $D_{\eta} \subset \mathbb{C}$ around the origin, the
fibre $f^{-1}(t)$ intersects $B_{\epsilon}$ transversally. Write
$X:=f^{-1}(D_{\eta}) \cap B_{\epsilon}$, $X_t:=f^{-1}(t) \cap
B_{\epsilon}$ for $t \in D_{\eta}$ and $D_{\eta}^*:=D_{\eta}
\setminus \{0\}$ for the pointed disc. Milnor showed that $f_{|_{X
\setminus X_0}}: X \setminus X_0 \rightarrow
D_{\eta}^*$ is a locally trivial fibration, see \cite{milnor}. A fibre %$f_{|_{X
%\setminus X_0}}^{-1}(t)$
$X_t$ of this bundle is called \emph{Milnor fibre of $f$ at $b$}. We
will denote it by $F_b$.
%\begin{example}
%\hfill $\square$
%\end{example}
Consider the loop $\gamma$ %: [0,1] \rightarrow D_{\eta}: t \mapsto
%\eta e^{2\pi i t}$ in $D_{\eta}^*$
encircling the origin once counterclockwise. Since $f_{|_{X
\setminus X_0}}$ is a locally trivial fibration, the loop $\gamma$
lifts to a diffeomorphism $h$ of the Milnor fibre $F_b$, which is
well determined up to homotopy. In this way $\gamma$ induces an
automorphism $h^*: H^i(F_b,\mathbb{C}) \rightarrow
H^i(F_b,\mathbb{C})$, $i \geq 0$, that is called the
\emph{monodromy transformation}. %$H^0(X_{\delta},\mathbb{C})=?$
%\begin{definition} The \emph{zeta function of monodromy} at the origin $\zeta_f$ associated
%to the polynomial $f$ is
%\[\zeta_f(t):=\prod_{n \geq 0}\left(\mbox{det}(id^* - t h^*;
%H^n(X_{\delta}, \mathbb{C})) \right)^{{(-1)}^{(n+1)}}.\]
%\end{definition}
\\
\\ The surfaces for which we will prove the monodromy
conjecture have exactly one isolated singularity in the origin. A
result of Milnor (see \cite{milnor}) then says that
$H^i(F_0,\mathbb{C})=0$, for $i \neq 0$ and $i \neq d-1$, and
$H^0(F_0,\mathbb{C})=\mathbb{C}$ with trivial monodromy action. The
formula of A'Campo (\cite{A'Campo}) describes the characteristic
polynomial of the monodromy action on $H^{d-1}(F_0,\mathbb{C})$ in
terms of an embedded resolution of the hypersurface $f^{-1}(0)$.
\\We may suppose that $\pi$ is an isomorphism outside
the inverse image of the origin. %Say that the $E_j$, $j \in
%J:=\{1,\cdots,r\}$, are the irreducible exceptional components of
%$\pi^{-1}(\{ 0 \})$.
\begin{theorem}(A'Campo) The characteristic polynomial of the monodromy action on
$H^{d-1}(F_0,\mathbb{C})$ is equal to
\[\left[\frac{\prod_{j=1}^r(1-t^{N_j})^{\chi(E_j^{\circ})}}{1-t}\right]^{(-1)^{d-1}}.\]
\end{theorem}
${}$\\
\textbf{\emph{$3.2.$ Topological zeta function.}}--- In $1992$ Denef
and Loeser created a new zeta function which they called the
topological zeta function because of the topological
Euler--Poincar\'e characteristic $\chi(\cdot)$ turning up in it. %Roughly said, the topological zeta
%function $Z_{top,f}$ associated to a polynomial function $f:
%\mathbb{C}^d \rightarrow \mathbb{C}$ (or to the germ $f:
%(\mathbb{C}^d,0) \rightarrow (\mathbb{C},0)$ of a holomorphic
%function) is a function containing information that we can pick out
%of every chosen embedded resolution of $f^{-1} \{ 0 \} \subset
%\mathbb{C}^d$.
It is associated to a complex polynomial $f$ with $f(0)=0$. If $E_I
:= \cap_{i \in I} E_i$ and $E_I^{\circ}:= E_I \setminus (\cup_{j
\notin I} E_j)$, then they introduced it in \cite{DenefLoeser1} in
the following way.
\begin{definition} The local \emph{topological zeta function
associated to $f$} is the rational function in one complex
variable\[Z_{top,f} (s)  := \sum_{I \subset S} \chi (E_I^{\circ}
\cap \pi^{-1}\{ 0 \}) \prod_{i \in I} \frac{1}{N_i s+ \nu_i}.\]
\end{definition}
%For a polynomial function $f$ there is a global version replacing
%${E_I}^{\circ} \cap \pi^{-1}\{ 0 \}$ by ${E_I}^{\circ}$. When we do
%not specify, we mean the local one.
%\\ \\
Denef and Loeser proved that every embedded resolution gives rise to
the same function, so the topological zeta function is a
well-defined singularity invariant (see \cite{DenefLoeser1}). Once
the motivic Igusa zeta function was introduced, they proved this
result alternatively in \cite{DenefLoeser2} by showing that this
more general zeta function specialises to the topological one. There
exists a global version, replacing $E_I^{\circ} \cap \pi^{-1}\{ 0
\}$ by $E_I^{\circ}$.
\\ \\
\newpage \noindent \textbf{\emph{$3.3.$ Monodromy conjecture.}}--- One calls
$\alpha$ an \emph{eigenvalue of monodromy of $f$ at $b \in
f^{-1}\{0\}$} if $\alpha$ is an eigenvalue for some $h^*: H^i(F_b,
\mathbb{C}) \rightarrow
H^i(F_b, \mathbb{C})$. %We can now
%formulate the monodromy conjecture:
%\begin{conjecture}(\emph{Monodromy conjecture}) If $s$ is a pole of
%$Z_{top,f}$, then $e^{2\pi i s}$ is an eigenvalue of monodromy of
%$f$ for some point of the hypersurface $f=0$.
%acting on the cohomology (in some dimension) of the Milnor fibre of
%$f$ associated to some point of the hypersurface $f=0$.
%\end{conjecture}
%The zeta functions $\zeta_f$ determine all the eigenvalues of the
%monodromy transformations. For isolated singularities this follows
%from the fact that the cohomology groups
%$H^n(X_{\delta},\mathbb{C})$ are all $\{0\}$ except for $n=0$ and
%$n=d-1$. The zeta function then becomes
%\begin{eqnarray*}
%\zeta_f(t)=\frac{\mbox{det}(id^* - t h^*; H^{d-1}(X_{\delta},
%\mathbb{C}))^{{(-1)}^{d}}}{1-t}. \end{eqnarray*}
%the fact that the Milnor fibre is
%then connected an that there is only cohomology in degrees $0$ and
%$d-1$.
%Also for non-isolated singularities the zeroes and poles of
%$\zeta_f$ are the eigenvalues of monodromy of $f$ (see e.g.
%\cite{denef}).
\begin{conjecture} (Monodromy Conjecture) If $s_0$ is a pole of $Z_{top,f}$, then $e^{2\pi i s_0}$ is an eigenvalue of monodromy of
$f$ at some point of the germ at $0$ of the hypersurface $f=0$.
\end{conjecture}
Let $f$ be a polynomial that is general with respect to a
$3$-dimensional toric idealistic cluster. Consider the embedded
resolution $\pi : Z \rightarrow \mathbb{C}^3$ of $f^{-1} \{ 0 \}$
that corresponds to the blowing up of the constellation. We fix a
candidate pole $s_0=-\nu_j/N_j$ of $Z_{top,f}$. If $E_j$ is not an
exceptional component, then $\nu_1=1$ and $N_1=1$. As $1$ is always
an eigenvalue of the local monodromy of $f$, this candidate pole
does not pose any difficulty. If $s_0=-\nu_j/N_j$ is a candidate
pole of $Z_{top,f}$ induced by an exceptional component $E_j$, then
we write $\nu_j/N_j$ as $a/b$ such that $a$ and $b$ are coprime. We
define the set $J_b:=\{j \in J \mbox{ $|$
    } \mbox{ }  b \mbox{ divides } N_j \}$. It follows from A'Campo's formula that
\begin{center} $e^{2\pi i s_0}$ is an eigenvalue of monodromy of $f$ at the origin $0$ \\
$\Updownarrow$
\\ $\sum_{j \in J_b} \chi(E_j^{\circ}) \neq 0.$
\end{center} \vspace{0.3cm}
In general there can be a lot of cancelations which make that
$\sum_{j \in J_b} \chi(E_j^{\circ}) = 0$. To control this, we will
determine when $\chi(E_j^{\circ})$ is positive, negative or zero. We
will see that the cases where $\chi(E_j^{\circ}) \leq 0$ are very
rare in this context.
\\ \\
\textbf{\emph{$3.4.$ Holomorphy conjecture.}}--- For every $r \in
\mathbb{Z}_{> 0}$, one can define a variant $Z_{top,f}^{(r)}$ of the
topogical zeta function that is also a rational function in one
complex variable.
\begin{definition}
\[Z_{top,f}^{(r)} := \sum_{\substack{I \subset S \\ \forall i \in I: r \mid N_i}} \chi
(E_I^{\circ} \cap \pi^{-1}\{ 0 \}) \prod_{i \in I} \frac{1}{N_i s+
\nu_i}.
\]
\end{definition}
The functions $Z_{top,f}^{(r)}$ are limits of more general Igusa
zeta function associated to a polynomial and a character (see
\cite{denef91}). In particular $Z_{top,f}^{(1)}=Z_{top,f}$. Clearly
they are holomorphic on $\mathbb{C}$ if and only if they do not have
a pole. The holomorphy conjecture stated by Denef predicts the
following relation.
\begin{conjecture}(Holomorphy conjecture) If $r \in \mathbb{Z}_{>
0}$ does not divide the order of any eigenvalue of monodromy of $f$,
then $Z_{top,f}^{(r)}$ is holomorphic on $\mathbb{C}$.
\end{conjecture}
In Section $9$ we provide a proof of the holomorphy conjecture for
the surfaces we are studying. Again, the classification of
$\chi(E_j^{\circ})$ according to the sign will be the key to solve
the conjecture.
\\ \\ \newpage ${}$\begin{center} \textsc{4. Computation of the topological zeta
function}
\end{center} ${}$\\
Given a germ of a polynomial function $f$ in $d$ variables over
$\mathbb{C}$, its topological zeta function $Z_{top,f}$ can be
calculated by computing an embedded resolution. If $f$ is
nondegenerate with respect to its Newton polyhedron, then there
exists also the formula for $Z_{top,f}$ in terms of its Newton
polyhedron, see \cite{DenefLoeser1}. In our context, we show that,
directly from the tree that represents the toric constellation, one
can read all
information needed to write down the topological zeta function. %The
%well-developed formulas for genera and for the topological
%Euler-Poincar\'{e} characteristic for curves on surfaces make that
%surfaces are very nice varieties to develop a formula for.
\\ \\
Concretely, we consider a toric idealistic cluster in $\mathbb{C}^3$
and a complex polynomial $f$ in three variables in a finitely
supported ideal such that the cluster gives an embedded resolution
for the surface $S:=V(f) \subset \mathbb{C}^3$. To determine the
topological zeta function of $f$, we determine the numbers $\chi
(E_I^{\circ})$. We will denote the strict transform of $S$ by
$\hat{S}$, whatever the stage is, and we will denote the curves
$\hat{S} \cap E_i$ by $C_i$. We will write $p_a$ for the geometric genus.\\
\indent First of all, notice that when blowing up in a point of
multiplicity $m$ on $S$, and $E$ being the created exceptional
divisor, the curve $\hat{S} \cap E$ has degree $m$. Another
important observation is that if $Q \in E$, then the multiplicity of
$Q$ on $\hat{S} \cap E$
is equal to the multiplicity of $Q$ on $\hat{S}$. %This follows from
%the fact that $E$ is transversal to $\hat{S}$ in $Q$. If not, the
%singularities could never be resolved by blowing up in points.
\\ \indent We give a formula for the topological zeta function but first
we illustrate the computation by following the embedded resolution
process of the following toric constellation. We think that such
concrete pictures are very useful to understand the computation of
the $\chi (E_I^{\circ})$ in general.
\begin{example}\label{exampleresolution} \emph{Consider the toric
constellation represented by the following tree. \vspace*{-0.7cm}}
\begin{center}
\begin{pspicture}(-2.5,-0.5)(3.5,4.5)
\psline{-}(0,0)(0,1.5) \psline{-}(0,1.5)(-1.5,2.5)
\psline{-}(0,1.5)(1.5,2.5)\psline{-}(1.5,2.5)(0,3.5)
\psline{-}(1.5,2.5)(1.5,3.5)\psline{-}(1.5,2.5)(3,3.5)\psdot(0,0)\psdot(0,1.5)\psdot(-1.5,2.5)
\psdot(1.5,2.5)\psdot(3,3.5)\psdot(0,3.5)\psdot(1.5,3.5)\rput(0.4,0){\footnotesize{$Q_1$}}\rput(0.4,1.3){\footnotesize{$Q_2$}}
\rput(-1.5,2.8){\footnotesize{$Q_3$}}\rput(1.7,2.1){\footnotesize{$Q_4$}}
\rput(-0.1,3.8){\footnotesize{$Q_5$}}\rput(1.5,3.8){\footnotesize{$Q_6$}}\rput(3.1,3.8){\footnotesize{$Q_7$}}\rput(-0.3,0.7){\scriptsize{\emph{1}}}
\rput(-1,1.8){\scriptsize{\emph{1}}}\rput(0.5,2.1){\scriptsize{\emph{2}}}
\rput(0.3,3){\scriptsize{\emph{1}}}\rput(1.7,3){\scriptsize{\emph{2}}}\rput(2.7,3){\scriptsize{\emph{3}}}
\end{pspicture} \end{center}
\vspace*{-0.3cm} \emph{Let $S$ be a surface in $\mathbb{C}^3$ that
is general for the above toric constellation. We follow the
resolution process and we picture the intersections that are
relevant in the calculation of the numbers $\chi (E_I^{\circ})$. The
gray curve (that can be reducible) pictured in the ambient $E_j$
represents the curve $C_j$.}
%\begin{center}
\begin{center}
\includegraphics{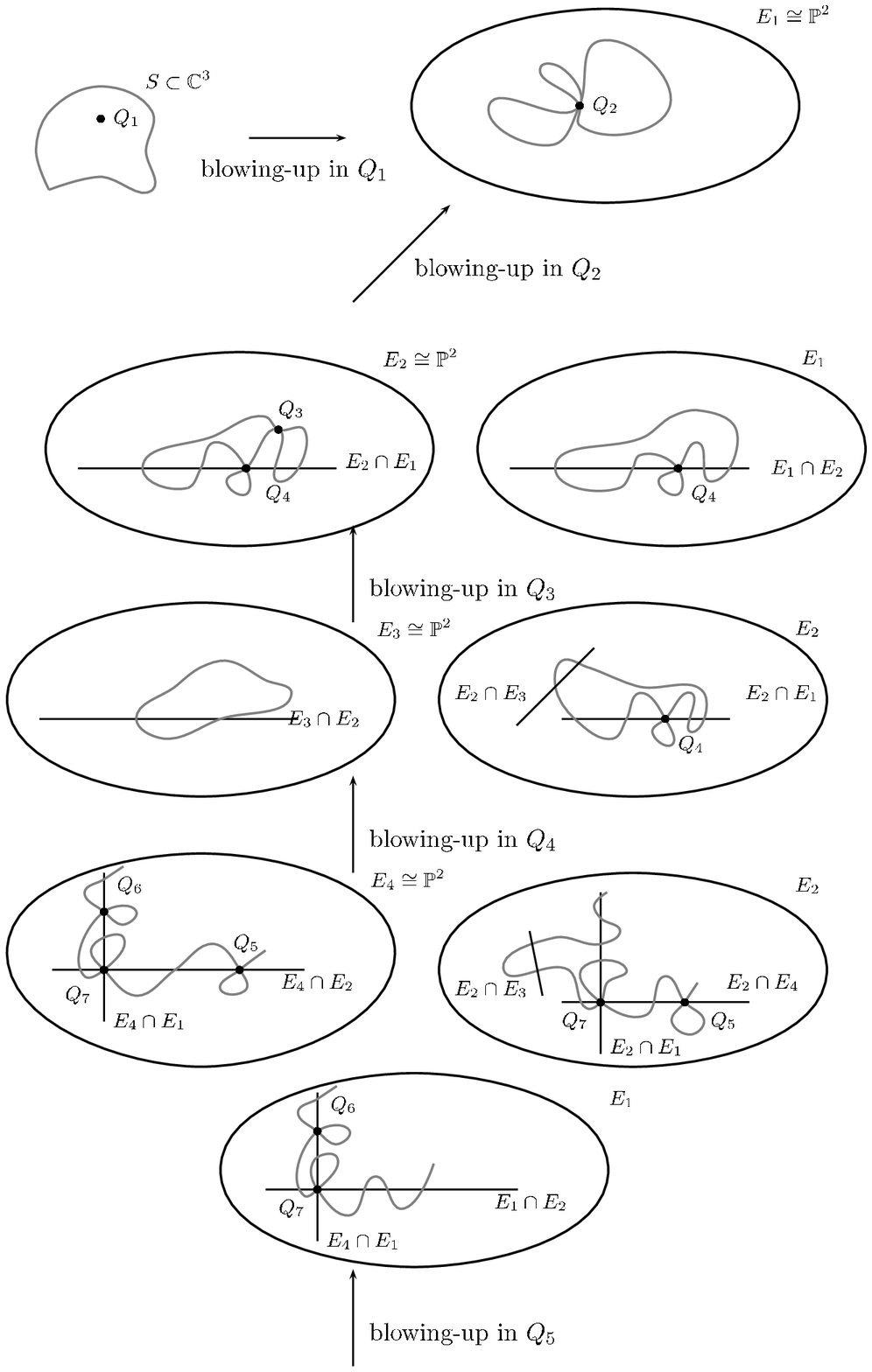}
\end{center}
\begin{center}
\includegraphics{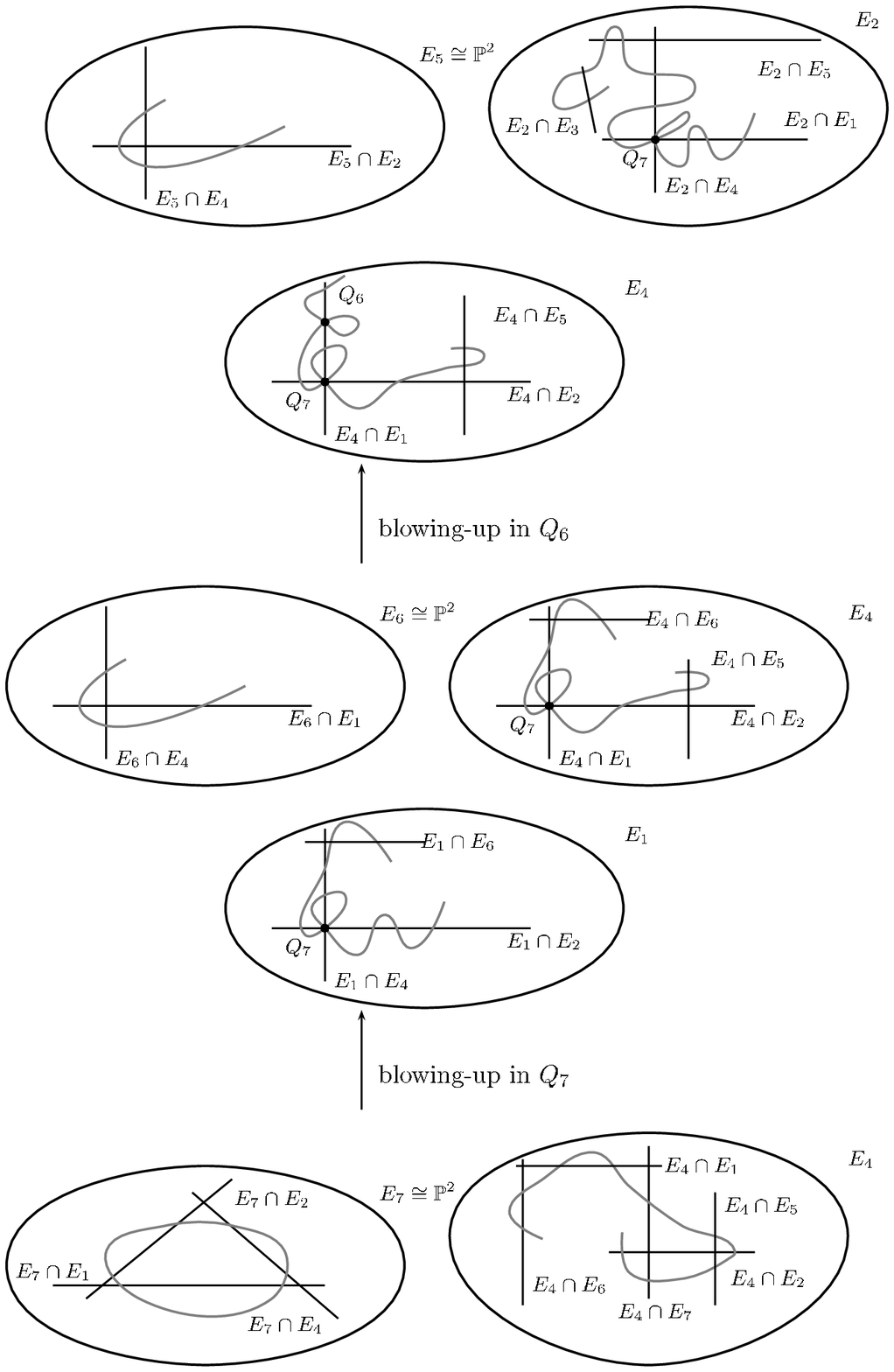}
\end{center}
\begin{center}
\newpage \vspace*{-2cm} \includegraphics{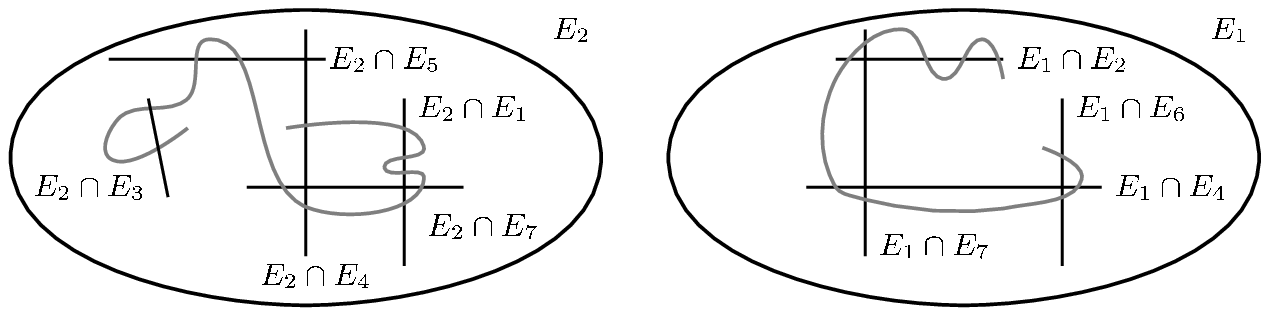}
\end{center}
%\end{center}
\end{example}\vspace*{-1.5cm}
\[\]\hfill $\square$
${}$\\ \\We now proceed to the computation of the $\chi
(E_I^{\circ})$. We will write $m_j$ for the multiplicity of the
point $Q_j$ on $\hat{S}$ and $E_0$ for the strict transform
$\hat{S}$.
\begin{enumerate}
\item \textbf{$I:=\{0,i,j\}$ with $0 < i < j$ and $j \rightarrow i$. }\\ \\
From the number of intersection points of $C_j$ and $E_i$ in $E_j
\cong \mathbb{P}^2$, we subtract the number of points in which we
will blow up. Then we get $\chi (E_I^{\circ})=m_j -
\sum_{\substack{k \succ j
\\ k\rightarrow i}}(C_j\cdot(E_i\cap E_j))_{Q_k}$. We can conclude
\[\chi (E_I^{\circ} )=m_j-\sum_{\substack{k \rightarrow i \\ k \twoheadrightarrow j}} m_k.\]
%Note that this number is positive because the cluster is satisfying
%the proximity inequalities as it is idealistic.
\item \vspace*{-0.3cm} \textbf{$I:=\{i,j,k\}$ with $0 \neq i < j < k, k \rightarrow i$ and $k \rightarrow
j$.}\\ \\The contribution to $\chi (E_I^{\circ})$ comes from the
intersection point of $E_i \cap E_j \cap E_k$ unless it is a point
in which we will blow up. We can express this as follows:
\[\chi (E_I^{\circ})=1- \# \{l \mbox{ $|$ } l \rightarrow i, l \rightarrow j \mbox{ and } l \rightarrow k\}.\]
\item \textbf{$I:=\{0,i\}$ with $0 \neq i$.} \\ \\
We look at $E_i$ in the final stage. There we have to subtract from
$E_0 \cap E_i$ the intersection points with the other exceptional
components.
\[\chi (E_I^{\circ} )=\chi(C_i)-\sum_{j\rightarrow
i}\chi(\overset{\circ}{\widehat{E_0 \cap E_i \cap
E_j}})-\sum_{i\rightarrow j}\chi(\overset{\circ}{\widehat{E_0 \cap
E_i \cap E_j}}).\] We have $\chi({C_i})=2-2p_a({C_i})$ for the
nonsingular ${C_i}$ that can be irreducible or reducible. This leads
to the formula
\begin{eqnarray*}
\chi (E_I^{\circ} ) & = & m_i(3-m_i) +
\sum_{j\rightarrow i}m_j(m_j-1) \\
& &- \sum_{j\rightarrow i}(m_j-\sum_{\substack{k \rightarrow i\\ k
\twoheadrightarrow j}} m_k)-\sum_{i\rightarrow
j}(m_i-\sum_{\substack{k \rightarrow j\\ k \twoheadrightarrow i}}
m_k).
\end{eqnarray*}
\item \vspace*{-0.3cm}\textbf{$I:=\{i,j\}$ with $0 \neq i < j, j \rightarrow i$.}\\
\\
We compute the contribution from the configuration in $E_j \cong
\mathbb{P}^2$.
\begin{eqnarray*}
\chi (E_I^{\circ} ) & = & 2-\left(\chi(\overset{\circ}{\widehat{E_0
\cap E_i \cap E_j}})+\#A_{ij}+\#B_{ij}-\#C_{ij} \right)
\\
&= & 2 - (m_j - \sum_{\substack{k \rightarrow i \\ k
\twoheadrightarrow j}} m_k) -\#A_{ij}-\#B_{ij}+\#C_{ij},
\end{eqnarray*}
 with
\vspace*{-0.7cm}\begin{eqnarray*}
A_{ij} & := & \{k\mbox{ $|$ } k \succ j, k \rightarrow i \}\\
B_{ij} & := & \{k\mbox{ $|$ } k \neq i, j \rightarrow k \}\\
C_{ij} & := & \{k\mbox{ $|$ } k \succ j, k \rightarrow i \mbox{ and } \exists l : l \neq i, k \rightarrow l \mbox{ and } j \rightarrow l \}.\\
\end{eqnarray*}
\item \vspace*{-0.5cm}\textbf{$I:=\{i\}$ with $i \neq 0$.}
\\ \\ We look in $E_i \cong \mathbb{P}^2$ and find
\begin{eqnarray*}
\chi (E_I^{\circ}) &= & 3-\left(\chi(\overset{\circ}{\widehat{E_0 \cap E_i}})+\#A_{i}+2\#B_{i}-\begin{pmatrix} \#B_i\\
2 \end{pmatrix} \right)
\\ & =& 3 + m_i(m_i-3) -
\sum_{j\rightarrow i}m_j(m_j-1) \\
& &+ \sum_{j\rightarrow i}(m_j-\sum_{\substack{k \rightarrow i\\ k
\twoheadrightarrow j}} m_k)+\sum_{i\rightarrow
j}(m_i-\sum_{\substack{k \rightarrow j\\ k \twoheadrightarrow i}}
m_k) - \\ & & \#A_{i}-2\#B_{i}+\begin{pmatrix} \#B_i\\
2 \end{pmatrix},
\end{eqnarray*}
 with
\vspace*{-0.7cm}\begin{eqnarray*} A_i & := & \{k \mbox{ $|$ } k
\succ i \mbox{ and } \nexists l: i
\rightarrow l \mbox{ and } k \rightarrow l\} \\
B_i & := & \{k \mbox{ $|$ } i \rightarrow k  \}.
\end{eqnarray*}
\item For $I$ not of the form of one of the sets described above, $\chi(E_I^{\circ} )=0$.
\end{enumerate}
${}$\\ Also the numerical data are completely determined by the
tree. We obtain the numbers $N_i$ via the recursive formula
$N_i=m_i+\sum_{i \rightarrow j}N_j$. For the $\nu_i$, we find $\nu_i
= \sum_{i \rightarrow j}(\nu_j - 1)+3$.
%define the following functions
%\begin{eqnarray*}
%WT_1(x,y,z) & = & (x+y+z+2,y,z),\\
%WT_2(x,y,z) & = & (x,x+y+z+2,z),\\
%WT_3(x,y,z) & = & (x,y,x+y+z+2)
%\end{eqnarray*}
%and code the points of the constellation as follows:
%$Q=Q(i_1,i_2,\cdots,i_k)$ if starting at the origin, one has to
%follow the edges coded by $i_1,i_2,\cdots,i_k$ to arrive at the
%point $Q$ in the tree. Then $\nu(Q(i_1,i_2,\cdots,i_k))$ is equal to
%the $j$-\emph{th} component of
%\[(WT_j\circ WT_{i_k}\circ \cdots \circ WT_{i_1})(0,0,0)\] with $j
%\in \{1,2,3\}$.
${}$
\\ \\  ${}$ \begin{center} \textsc{5. Analysis of $\chi(E_i^{\circ})$}
\end{center} ${}$\\
In order to investigate the conjectures, we study the expression for
$\chi(E_i^{\circ})$ that we obtained in the previous section:
\[\chi (E_i^{\circ}) =  m_i(m_i-3) -
\sum_{j\rightarrow i}m_j(m_j-1) + \sum_{j\rightarrow
i}(m_j-\sum_{\substack{k \rightarrow i\\ k \twoheadrightarrow j}}
m_k)\] \vspace*{-0.4cm}\[ \qquad \quad +\sum_{i\rightarrow j}(m_i-\sum_{\substack{k \rightarrow j\\
k \twoheadrightarrow i}}
m_k) +3-\#A_{i}-2\#B_{i}+\begin{pmatrix} \#B_i\\
2 \end{pmatrix},\]
 with $A_i =  \{k \mbox{ $|$ } k \succ i \mbox{ and } \nexists l: i
\rightarrow l \mbox{ and } k \rightarrow l\}$ and $B_i = \{k \mbox{
$|$ } i \rightarrow k  \}$.
\\
Notice that the linear proximity inequalities imply that
$m_j-\sum_{\substack{k \rightarrow i\\ k \twoheadrightarrow j}} m_k
\geq 0$, for all $j \rightarrow i$ and that $m_i-\sum_{\substack{k
\rightarrow j\\ k \twoheadrightarrow i}} m_k \geq 0$ for all $i
\rightarrow j$. Moreover for a point $Q_j$ with maximal level in the
set of the points that are proximate to $Q_i$, we have
$m_j-\sum_{\substack{k \rightarrow i\\ k \twoheadrightarrow j}} m_k
=m_j > 0$ and so $\sum_{j\rightarrow i}(m_j-\sum_{\substack{k
\rightarrow i\\ k \twoheadrightarrow j}} m_k) > 0$.
\\Let $T:=3-\#A_{i}-2\#B_{i}+\begin{pmatrix} \#B_i\\
2 \end{pmatrix}$. Then $T$ takes the following values:
\begin{center}
\begin{tabular}[t]{|c|c|c|} \hline $ \quad  \#B_i$ \quad & \quad $\#A_i$ \quad & \quad
$T$ \quad
\\ \hline \vspace*{-0.4cm}& & \\
$0$ & $0$ & $3$  \\
$0$ & $1$ & $2$\\
$0$ & $2$ & $1$ \\
$0$ & $3$ & $0$ \\
 \hline \vspace*{-0.4cm}& & \\
$1$ & $0$ & $1$  \\
$1$ & $1$ & $0$\\
 \hline \vspace*{-0.4cm}& & \\
$2$ & $0$ & $0$ \\
 \hline \vspace*{-0.4cm}& & \\
$3$ & $0$ & $0$ \\
\hline
\end{tabular}
\end{center}
\begin{center}
\emph{Table 1}
\end{center}
We want to investigate when $\chi(E_i^{\circ}) \leq 0$.
%\begin{example}
%We consider the constellation.... Maximum of the distance to the
%point $(1/2,1/2,1/2)$ in the compact polytope....
%
%\end{example}
%According the constellation we can get polytopes of very different
%forms in arbitrary dimension. Then we treat this quadratic integer
%problem in another way.
A priori there are infinitely many constellations to consider. The
first result in this section will permit us to reduce our study to a
finite number of cases. Secondly we will rewrite $\chi (E_i^{\circ}
)$ and via combinatorics we will analyse this new description.
\begin{lemma} \label{finite}
Let $\mathcal{A}=(\mathcal{C},\underline{m})$ be a $3$-dimensional
toric idealistic cluster and let $Q_i$ be a point of the
constellation $\mathcal{C}$. If $\#\{t \in \mathbb{Z}_{\geq 0}\mbox{
$|$ } Q_i(a,b^t) \in \mathcal{C}\} + \#\{t \in \mathbb{Z}_{\geq
0}\mbox{ $|$ } Q_i(b,a^t) \in \mathcal{C}\} \geq 3$ for all $a, b
\in \{1,2,3\}$ with $a \neq b$, then $\chi(E_i^{\circ})
> 0$.
\end{lemma}
${}$\\
\emph{Proof.} \quad If $\#\{t \in \mathbb{Z}_{\geq 0}\mbox{ $|$ }
Q_i(a,b^t) \in \mathcal{C}\} + \#\{t \in \mathbb{Z}_{\geq 0}\mbox{
$|$ } Q_i(b,a^t) \in \mathcal{C}\} \geq 3$ for all $a, b \in
\{1,2,3\}$ with $a \neq b$, then it follows that $m_i > 3$ except
when there are exactly $6$ points - that have multiplicity $1$ -
that are proximate to $Q_i$ and such that $\#\{t \in
\mathbb{Z}_{\geq 0}\mbox{ $|$ } Q_i(1,3^t) \in \mathcal{C}\}=\#\{t
\in \mathbb{Z}_{\geq 0}\mbox{ $|$ } Q_i(2,1^t) \in
\mathcal{C}\}=\#\{t \in \mathbb{Z}_{\geq 0}\mbox{ $|$ } Q_i(3,2^t)
\in \mathcal{C}\}=2$, up to permutation of the labels. In that case
$m_i$ can be equal to $3$ and one then finds that $\chi
(E_i^{\circ}) > 0$. \\ \indent When $m_i>3$ we construct a new
cluster. We define $m_i':=m_i-3$, $m_j':=m_j-1$ for all $j$ for
which $j \rightarrow i$ and we do not change the weights of the
other points in $\mathcal{C}$. Let $\mathcal{C}'$ be the
subconstellation of $\mathcal{C}$ that contains exactly the points
$Q_j$ of $\mathcal{C}$ for which $m_j > 1$ and let $\mathcal{A'}$ be
the cluster $(\mathcal{C}',\underline{m'})$. Then also
$\mathcal{A'}$ satisfies the linear proximity inequalities and thus
$\mathcal{A'}$ is a toric idealistic cluster. Let us now consider a
surface $\mathcal{S}$ that is general with respect to $\mathcal{A}$
and a surface $\mathcal{S'}$ that is general with respect to
$\mathcal{A'}$. Blowing up the point $Q_i$ provides two curves
$C_i=E_i \cap \hat{S}$ and $C_i'=E_i \cap \hat{S'}$ in the
exceptional variety $E_i \cong \mathbb{P}^2$ of degree $m_i$ and
$m_i'$, respectively. From Bezout's formula it follows that $m_i
m_i' \geq \sum_{\substack{j \rightarrow i \\ Q_j \in \mathcal{C}'}}
m_j m_j'$.
The latter sum is also equal to $\sum_{\substack{j \rightarrow i \\
Q_j \in \mathcal{C}}} m_j m_j'$. We can conclude that $\chi
(E_i^{\circ} )
> 0$. \hfill $\blacksquare$
\\ \\
This lemma will allow us to work with a finite number of families of
constellations. We represent these families in List $1$. We first
explain some notations. \\  \indent To save place, from now on we
draw the clusters from left to right. So if there is an edge between
$Q_i$ and $Q_j$ and if $Q_j$ is at the right from $Q_i$, then $Q_j >
Q_i$. If there exists an edge with label $x$ between points of the
chain $\mathcal{C}^i:=\{Q_j \mid Q_i \geq Q_j\}$, then we will
simply say that `label $x$ appears below $Q_i$'.
\\  \indent The constellations are listed according to the number of
points $Q_j$ for which $Q_j \succ Q_i$ (indicated by a roman
number). %We have given a roman number to the constellations from
%List $1$, according to the number of points $Q_j$ for which $Q_j
%\succ Q_i$.
We only draw the subconstellation that shows $Q_i$ and the points
$Q_j$ that are proximate to $Q_i$ and for which holds that $j \succ
i$ or for which there exists a point $Q_k$ such that $k \succ i$ and
$j \succ k$. By drawing `$--$' going out of a point $Q_j$ for which
$j \rightarrow i$, we mean that there can exist a point $Q_k$ for
which $k > j$ and $k \rightarrow i$. We also draw the symbol `$--$'
arriving in the point $Q_i$ when $Q_i$ is not necessarily the
origin. When $Q_j$ is a point of the constellation, we will denote
its multiplicity by $m_{Q_j}$ or by $m_j$.
\\ \indent List $1$
contains the constellations we should study - according to Lemma
\ref{finite} - up to permutation of the labels. In constellations
II9, II10 and II11, we mean by \^{3} that label $3$ should not occur
at that place, so $\#\{t \in \mathbb{Z}_{\geq 0}\mbox{ $|$ }
Q_i(2,3^t) \in \mathcal{C}\}=2$.
\\ \\
\begin{pspicture}(0,-1)(4.9,1)
\rput(0.2,0.8){\textbf{01}} \psline{-,linestyle=dashed}(0,0)(0.5,0)
\psdot(0.5,0) \rput(0.5,0.3){\footnotesize{$Q_i$}}
\end{pspicture}
\begin{pspicture}(0,-1)(4.9,1)
\rput(0.2,0.8){\textbf{I1}} \psline{-,linestyle=dashed}(0,0)(0.5,0)
\psdot(0.5,0) \rput(0.5,0.3){\footnotesize{$Q_i$}}
\psline{-}(0.5,0)(1.5,0) \psdot(1.5,0)
\rput(1,0.2){\scriptsize{\emph{1}}}
\end{pspicture}
\begin{pspicture}(0,-1)(4.9,1)
\rput(0.2,0.8){\textbf{I2}} \psline{-,linestyle=dashed}(0,0)(0.5,0)
\psdot(0.5,0) \rput(0.5,0.3){\footnotesize{$Q_i$}}
\psline{-}(0.5,0)(1.5,0) \psdot(1.5,0) \psline{-}(1.5,0)(2.5,0)
\psdot(2.5,0)
\rput(1,0.2){\scriptsize{\emph{1}}}\rput(2,0.2){\scriptsize{\emph{2}}}\psline{-,linestyle=dashed}(2.5,0)(3,0)
\end{pspicture}\\ \\
\begin{pspicture}(0,-1)(4.9,1)
\rput(0.2,0.8){\textbf{I3}} \psline{-,linestyle=dashed}(0,0)(0.5,0)
\psdot(0.5,0) \rput(0.5,0.3){\footnotesize{$Q_i$}}
\psline{-}(0.5,0)(1.5,0) \psdot(1.5,0) \psline{-}(1.5,0)(2.5,0.5)
\psdot(2.5,0.5) \psline{-}(1.5,0)(2.5,-0.5) \psdot(2.5,-0.5)
\rput(1,0.2){\scriptsize{\emph{1}}}\rput(2,0.45){\scriptsize{\emph{2}}}\rput(2,-0.5){\scriptsize{\emph{3}}}
\psline{-,linestyle=dashed}(2.5,0.5)(3,0.5)
\psline{-,linestyle=dashed}(2.5,-0.5)(3,-0.5)
\end{pspicture}
\begin{pspicture}(0,-1)(4.9,1)
\rput(0.2,0.8){\textbf{II1}} \psline{-,linestyle=dashed}(0,0)(0.5,0)
\psdot(0.5,0) \rput(0.45,0.3){\footnotesize{$Q_i$}}
\psline{-}(0.5,0)(1.5,0.5) \psline{-}(0.5,0)(1.5,-0.5)
\psdot(1.5,0.5)
\psdot(1.5,-0.5)\rput(1,0.45){\scriptsize{\emph{1}}}\rput(1,-0.5){\scriptsize{\emph{2}}}
\end{pspicture}
\begin{pspicture}(0,-1)(4.9,1)
\rput(0.2,0.8){\textbf{II2}} \psline{-,linestyle=dashed}(0,0)(0.5,0)
\psdot(0.5,0) \rput(0.45,0.3){\footnotesize{$Q_i$}}
\psline{-}(0.5,0)(1.5,0.5) \psline{-}(0.5,0)(1.5,-0.5)
\psdot(1.5,0.5)
\psdot(1.5,-0.5)\rput(1,0.45){\scriptsize{\emph{1}}}\rput(1,-0.5){\scriptsize{\emph{2}}}
\psline{-}(1.5,0.5)(2.5,0.5)\psdot(2.5,0.5)\psline{-,linestyle=dashed}(2.5,0.5)(3,0.5)
\rput(2,0.7){\scriptsize{\emph{2}}}
\end{pspicture}
\\ \\
\begin{pspicture}(0,-1)(4.9,1)
\rput(0.2,0.8){\textbf{II3}} \psline{-,linestyle=dashed}(0,0)(0.5,0)
\psdot(0.5,0) \rput(0.45,0.3){\footnotesize{$Q_i$}}
\psline{-}(0.5,0)(1.5,0.5) \psline{-}(0.5,0)(1.5,-0.5)
\psdot(1.5,0.5)
\psdot(1.5,-0.5)\rput(1,0.45){\scriptsize{\emph{1}}}\rput(1,-0.5){\scriptsize{\emph{2}}}
\psline{-}(1.5,0.5)(2.5,0.5)\psdot(2.5,0.5)\psline{-,linestyle=dashed}(2.5,0.5)(3,0.5)
\rput(2,0.7){\scriptsize{\emph{3}}}
\end{pspicture}
\begin{pspicture}(0,-1)(4.9,1)
\rput(0.2,0.8){\textbf{II4}} \psline{-,linestyle=dashed}(0,0)(0.5,0)
\psdot(0.5,0) \rput(0.45,0.3){\footnotesize{$Q_i$}}
\psline{-}(0.5,0)(1.5,0.5) \psline{-}(0.5,0)(1.5,-0.5)
\psdot(1.5,0.5)
\psdot(1.5,-0.5)\rput(1,0.45){\scriptsize{\emph{1}}}\rput(1,-0.5){\scriptsize{\emph{2}}}
\psline{-}(1.5,0.5)(2.5,1) \psdot(2.5,1) \psline{-}(1.5,0.5)(2.5,0)
\psdot(2.5,0)\psline{-,linestyle=dashed}(2.5,1)(3,1)
\rput(2,0.95){\scriptsize{\emph{2}}}
\rput(2,0.05){\scriptsize{\emph{3}}}\psline{-,linestyle=dashed}(2.5,0)(3,0)
\end{pspicture}
\begin{pspicture}(0,-1)(4.9,1)
\rput(0.2,0.8){\textbf{II5}} \psline{-,linestyle=dashed}(0,0)(0.5,0)
\psdot(0.5,0) \rput(0.45,0.3){\footnotesize{$Q_i$}}
\psline{-}(0.5,0)(1.5,0.5) \psline{-}(0.5,0)(1.5,-0.5)
\psdot(1.5,0.5)
\psdot(1.5,-0.5)\rput(1,0.45){\scriptsize{\emph{1}}}\rput(1,-0.5){\scriptsize{\emph{2}}}
\psline{-}(1.5,0.5)(2.5,0.5)\psdot(2.5,0.5)\psline{-,linestyle=dashed}(2.5,0.5)(3,0.5)
\rput(2,0.7){\scriptsize{\emph{2}}}
\psline{-}(1.5,-0.5)(2.5,-0.5)\psdot(2.5,-0.5)\psline{-,linestyle=dashed}(2.5,-0.5)(3,-0.5)
\rput(2,-0.7){\scriptsize{\emph{1}}}
\end{pspicture}
\\ \\
\begin{pspicture}(0,-1)(4.9,1)
\rput(0.2,0.8){\textbf{II6}} \psline{-,linestyle=dashed}(0,0)(0.5,0)
\psdot(0.5,0) \rput(0.45,0.3){\footnotesize{$Q_i$}}
\psline{-}(0.5,0)(1.5,0.5) \psline{-}(0.5,0)(1.5,-0.5)
\psdot(1.5,0.5)
\psdot(1.5,-0.5)\rput(1,0.45){\scriptsize{\emph{1}}}\rput(1,-0.5){\scriptsize{\emph{2}}}
\psline{-}(1.5,0.5)(2.5,0.5)\psdot(2.5,0.5)\psline{-,linestyle=dashed}(2.5,0.5)(3,0.5)
\rput(2,0.7){\scriptsize{\emph{2}}}
\psline{-}(1.5,-0.5)(2.5,-0.5)\psdot(2.5,-0.5)\psline{-,linestyle=dashed}(2.5,-0.5)(3,-0.5)
\rput(2,-0.7){\scriptsize{\emph{3}}}
\end{pspicture}
\begin{pspicture}(0,-1)(4.9,1)
\rput(0.2,0.8){\textbf{II7}}\psline{-,linestyle=dashed}(0,0)(0.5,0)
\psdot(0.5,0) \rput(0.45,0.3){\footnotesize{$Q_i$}}
\psline{-}(0.5,0)(1.5,0.5) \psline{-}(0.5,0)(1.5,-0.5)
\psdot(1.5,0.5)
\psdot(1.5,-0.5)\rput(1,0.45){\scriptsize{\emph{1}}}\rput(1,-0.5){\scriptsize{\emph{2}}}
\psline{-}(1.5,0.5)(2.5,0.5)\psdot(2.5,0.5)\psline{-,linestyle=dashed}(2.5,0.5)(3,0.5)
\rput(2,0.7){\scriptsize{\emph{3}}}
\psline{-}(1.5,-0.5)(2.5,-0.5)\psdot(2.5,-0.5)\psline{-,linestyle=dashed}(2.5,-0.5)(3,-0.5)
\rput(2,-0.7){\scriptsize{\emph{3}}}
\end{pspicture}
\begin{pspicture}(0,-1)(4.9,1)
\rput(0.2,0.8){\textbf{II8}}\psline{-,linestyle=dashed}(0,0)(0.5,0)
\psdot(0.5,0) \rput(0.45,0.3){\footnotesize{$Q_i$}}
\psline{-}(0.5,0)(1.5,0.5) \psline{-}(0.5,0)(1.5,-0.5)
\psdot(1.5,0.5)
\psdot(1.5,-0.5)\rput(1,0.45){\scriptsize{\emph{1}}}\rput(1,-0.5){\scriptsize{\emph{2}}}
\psline{-}(1.5,0.5)(2.5,1) \psdot(2.5,1) \psline{-}(1.5,0.5)(2.5,0)
\psdot(2.5,0)\psline{-,linestyle=dashed}(2.5,1)(3,1)
\rput(2,0.95){\scriptsize{\emph{2}}}
\rput(2,0.05){\scriptsize{\emph{3}}}\psline{-,linestyle=dashed}(2.5,0)(3,0)
\psline{-}(1.5,-0.5)(2.5,-0.5)\psdot(2.5,-0.5)\psline{-,linestyle=dashed}(2.5,-0.5)(3,-0.5)
\rput(2,-0.7){\scriptsize{\emph{1}}}
\end{pspicture}
\\ \\
\begin{pspicture}(0,-1)(4.9,1)
\rput(0.2,0.8){\textbf{II9}}\psline{-,linestyle=dashed}(0,0)(0.5,0)
\psdot(0.5,0) \rput(0.45,0.3){\footnotesize{$Q_i$}}
\psline{-}(0.5,0)(1.5,0.5) \psline{-}(0.5,0)(1.5,-0.5)
\psdot(1.5,0.5)
\psdot(1.5,-0.5)\rput(1,0.45){\scriptsize{\emph{1}}}\rput(1,-0.5){\scriptsize{\emph{2}}}
\psline{-}(1.5,0.5)(2.5,1) \psdot(2.5,1) \psline{-}(1.5,0.5)(2.5,0)
\psdot(2.5,0)\psline{-,linestyle=dashed}(2.5,1)(3,1)\psline{-,linestyle=dashed}(2.5,0)(3,0)
\rput(2,0.95){\scriptsize{\emph{2}}}\rput(2.8,-0.75){\scriptsize{\emph{\^{3}}}}
\rput(2,0.05){\scriptsize{\emph{3}}}
\psline{-}(1.5,-0.5)(2.5,-0.5)\psdot(2.5,-0.5)\psline{-,linestyle=dashed}(2.5,-0.5)(3,-0.5)
\rput(2,-0.7){\scriptsize{\emph{3}}}
\end{pspicture}
\begin{pspicture}(0,-1)(4.9,1)
\rput(0.2,0.8){\textbf{II10}}\psline{-,linestyle=dashed}(0,0)(0.5,0)
\psdot(0.5,0) \rput(0.45,0.3){\footnotesize{$Q_i$}}
\psline{-}(0.5,0)(1.5,0.5) \psline{-}(0.5,0)(1.5,-0.5)
\psdot(1.5,0.5)
\psdot(1.5,-0.5)\rput(1,0.45){\scriptsize{\emph{1}}}\rput(1,-0.5){\scriptsize{\emph{2}}}
\psline{-}(1.5,-0.5)(2.5,0) \psdot(2.5,0)
\psline{-}(1.5,-0.5)(2.5,-1)
\psdot(2.5,-1)\psline{-,linestyle=dashed}(2.5,-1)(3,-1)
\rput(2,-0.1){\scriptsize{\emph{1}}}
\rput(2,-1){\scriptsize{\emph{3}}}\psline{-,linestyle=dashed}(2.5,0)(3,0)
\psline{-}(1.5,0.5)(2.5,0.5)\psdot(2.5,0.5)
\rput(2,0.7){\scriptsize{\emph{3}}}\psline{-,linestyle=dashed}(2.5,0.5)(3,0.5)
\rput(2.8,-1.25){\scriptsize{\emph{\^{3}}}}
\end{pspicture}
\begin{pspicture}(0,-1)(4.9,1)
\rput(0.2,0.8){\textbf{II11}}\psline{-,linestyle=dashed}(0,0)(0.5,0)
\psdot(0.5,0) \rput(0.45,0.3){\footnotesize{$Q_i$}}
\psline{-}(0.5,0)(1.5,0.5) \psline{-}(0.5,0)(1.5,-0.5)
\psdot(1.5,0.5)
\psdot(1.5,-0.5)\rput(1,0.45){\scriptsize{\emph{1}}}\rput(1,-0.5){\scriptsize{\emph{2}}}
\psline{-}(1.5,0.5)(2.5,0.8) \psdot(2.5,0.8)
\psline{-}(1.5,0.5)(2.5,0.2)
\psdot(2.5,0.2)\psline{-,linestyle=dashed}(2.5,0.8)(3,0.8)
\rput(2,0.85){\scriptsize{\emph{2}}}
\rput(2,0.18){\scriptsize{\emph{3}}}\psline{-,linestyle=dashed}(2.5,0.2)(3,0.2)
\rput(2.8,-1){\scriptsize{\emph{\^{3}}}}
\psline{-}(1.5,-0.5)(2.5,-0.2)\psdot(2.5,-0.2)\psline{-,linestyle=dashed}(2.5,-0.2)(3,-0.2)
\psline{-}(1.5,-0.5)(2.5,-0.8)\psdot(2.5,-0.8)\psline{-,linestyle=dashed}(2.5,-0.8)(3,-0.8)
\rput(2,-0.8){\scriptsize{\emph{3}}}
\rput(2,-0.14){\scriptsize{\emph{1}}}
\end{pspicture}
\\ \\
\begin{pspicture}(0,-1)(4.9,1)
\rput(0.2,0.8){\textbf{III1}}\psline{-,linestyle=dashed}(0,0)(0.5,0)
\psdot(0.5,0) \rput(0.45,0.3){\footnotesize{$Q_i$}}
\psline{-}(0.5,0)(1.5,0.5) \psdot(1.5,0.5)\psline{-}(0.5,0)(1.5,0)
\psdot(1.5,0) \psline{-}(0.5,0)(1.5,-0.5) \psdot(1.5,-0.5)
\rput(1,0.45){\scriptsize{\emph{1}}}\rput(1,-0.5){\scriptsize{\emph{3}}}\rput(1.2,0.15){\scriptsize{\emph{2}}}
\end{pspicture}
\begin{pspicture}(0,-1)(4.9,1)
\rput(0.2,0.8){\textbf{III2}}\psline{-,linestyle=dashed}(0,0)(0.5,0)
\psdot(0.5,0) \rput(0.45,0.3){\footnotesize{$Q_i$}}
\psline{-}(0.5,0)(1.5,0.5) \psdot(1.5,0.5)\psline{-}(0.5,0)(1.5,0)
\psdot(1.5,0) \psline{-}(0.5,0)(1.5,-0.5) \psdot(1.5,-0.5)
\rput(1,0.45){\scriptsize{\emph{1}}}\rput(1,-0.5){\scriptsize{\emph{3}}}\rput(1.2,0.15){\scriptsize{\emph{2}}}
\psline{-}(1.5,0.5)(2.5,0.5)\psdot(2.5,0.5)\psline{-,linestyle=dashed}(2.5,0.5)(3,0.5)
\rput(2,0.7){\scriptsize{\emph{2}}}
\end{pspicture}
\begin{pspicture}(0,-1)(4.9,1)
\rput(0.2,0.8){\textbf{III3}}\psline{-,linestyle=dashed}(0,0)(0.5,0)
\psdot(0.5,0) \rput(0.45,0.3){\footnotesize{$Q_i$}}
\psline{-}(0.5,0)(1.5,0.5) \psdot(1.5,0.5)\psline{-}(0.5,0)(1.5,0)
\psdot(1.5,0) \psline{-}(0.5,0)(1.5,-0.5) \psdot(1.5,-0.5)
\rput(1,0.45){\scriptsize{\emph{1}}}\rput(1,-0.5){\scriptsize{\emph{3}}}\rput(1.2,0.15){\scriptsize{\emph{2}}}
\psline{-}(1.5,0.5)(2.5,1) \psdot(2.5,1)
\psline{-}(1.5,0.5)(2.5,0.5)
\psdot(2.5,0.5)\psline{-,linestyle=dashed}(2.5,1)(3,1)
\rput(2,0.95){\scriptsize{\emph{2}}}
\rput(2,0.3){\scriptsize{\emph{3}}}\psline{-,linestyle=dashed}(2.5,0.5)(3,0.5)
\end{pspicture}
\\ \\
\begin{pspicture}(0,-1)(4.9,1)
\rput(0.2,0.8){\textbf{III4}}\psline{-,linestyle=dashed}(0,0)(0.5,0)
\psdot(0.5,0) \rput(0.45,0.3){\footnotesize{$Q_i$}}
\psline{-}(0.5,0)(1.5,0.5) \psdot(1.5,0.5)\psline{-}(0.5,0)(1.5,0)
\psdot(1.5,0) \psline{-}(0.5,0)(1.5,-0.5) \psdot(1.5,-0.5)
\rput(1,0.45){\scriptsize{\emph{1}}}\rput(1,-0.5){\scriptsize{\emph{3}}}\rput(1.2,0.15){\scriptsize{\emph{2}}}
\psline{-}(1.5,0.5)(2.5,0.5)\psdot(2.5,0.5)\psline{-,linestyle=dashed}(2.5,0.5)(3,0.5)
\rput(2,0.7){\scriptsize{\emph{2}}}\psline{-}(1.5,0)(2.5,0)\psdot(2.5,0)\psline{-,linestyle=dashed}(2.5,0)(3,0)
\rput(2,0.2){\scriptsize{\emph{1}}}
\end{pspicture}
\begin{pspicture}(0,-1)(4.9,1)
\rput(0.2,0.8){\textbf{III5}}\psline{-,linestyle=dashed}(0,0)(0.5,0)
\psdot(0.5,0) \rput(0.45,0.3){\footnotesize{$Q_i$}}
\psline{-}(0.5,0)(1.5,0.5) \psdot(1.5,0.5)\psline{-}(0.5,0)(1.5,0)
\psdot(1.5,0) \psline{-}(0.5,0)(1.5,-0.5) \psdot(1.5,-0.5)
\rput(1,0.45){\scriptsize{\emph{1}}}\rput(1,-0.5){\scriptsize{\emph{3}}}\rput(1.2,0.15){\scriptsize{\emph{2}}}
\psline{-}(1.5,0.5)(2.5,0.5)\psdot(2.5,0.5)\psline{-,linestyle=dashed}(2.5,0.5)(3,0.5)
\rput(2,0.7){\scriptsize{\emph{3}}}\psline{-}(1.5,0)(2.5,0)\psdot(2.5,0)\psline{-,linestyle=dashed}(2.5,0)(3,0)
\rput(2,0.2){\scriptsize{\emph{1}}}
\end{pspicture}
\begin{pspicture}(0,-1)(4.9,1)
\rput(0.2,0.8){\textbf{III6}}\psline{-,linestyle=dashed}(0,0)(0.5,0)
\psdot(0.5,0) \rput(0.45,0.3){\footnotesize{$Q_i$}}
\psline{-}(0.5,0)(1.5,0.5) \psdot(1.5,0.5)\psline{-}(0.5,0)(1.5,0)
\psdot(1.5,0) \psline{-}(0.5,0)(1.5,-0.5) \psdot(1.5,-0.5)
\rput(1,0.45){\scriptsize{\emph{1}}}\rput(1,-0.5){\scriptsize{\emph{3}}}\rput(1.2,0.15){\scriptsize{\emph{2}}}
\psline{-}(1.5,-0.5)(2.5,-0.5)\psdot(2.5,-0.5)\psline{-,linestyle=dashed}(2.5,-0.5)(3,-0.5)
\rput(2,-0.7){\scriptsize{\emph{1}}}\psline{-}(1.5,0)(2.5,0)\psdot(2.5,0)\psline{-,linestyle=dashed}(2.5,0)(3,0)
\rput(2,0.2){\scriptsize{\emph{1}}}
\end{pspicture}
\\ \\
\begin{pspicture}(0,-1)(4.9,1)
\rput(0.2,0.8){\textbf{III7}}\psline{-,linestyle=dashed}(0,0)(0.5,0)
\psdot(0.5,0) \rput(0.45,0.3){\footnotesize{$Q_i$}}
\psline{-}(0.5,0)(1.5,0.5) \psdot(1.5,0.5)\psline{-}(0.5,0)(1.5,0)
\psdot(1.5,0) \psline{-}(0.5,0)(1.5,-0.5) \psdot(1.5,-0.5)
\rput(1,0.45){\scriptsize{\emph{1}}}\rput(1,-0.5){\scriptsize{\emph{3}}}\rput(1.2,0.15){\scriptsize{\emph{2}}}
%\psline{-}(1.5,0.5)(2.5,0.5)\psdot(2.5,0.5)\psline{-,linestyle=dashed}(2.5,0.5)(3,0.5)
%\rput(2,0.7){\scriptsize{\emph{2}}}
\psline{-}(1.5,0)(2.5,0)\psdot(2.5,0)\psline{-,linestyle=dashed}(2.5,0)(3,0)
\rput(2,0.2){\scriptsize{\emph{1}}} \psline{-}(1.5,0.5)(2.5,1)
\psdot(2.5,1) \psline{-}(1.5,0.5)(2.5,0.5)
\psdot(2.5,0.5)\psline{-,linestyle=dashed}(2.5,1)(3,1)
\rput(2,0.95){\scriptsize{\emph{2}}}
\rput(2.3,0.7){\scriptsize{\emph{3}}}\psline{-,linestyle=dashed}(2.5,0.5)(3,0.5)
\end{pspicture}
\begin{pspicture}(0,-1)(4.9,1)
\rput(0.2,0.8){\textbf{III8}}\psline{-,linestyle=dashed}(0,0)(0.5,0)
\psdot(0.5,0) \rput(0.45,0.3){\footnotesize{$Q_i$}}
\psline{-}(0.5,0)(1.5,0.5) \psdot(1.5,0.5)\psline{-}(0.5,0)(1.5,0)
\psdot(1.5,0) \psline{-}(0.5,0)(1.5,-0.5) \psdot(1.5,-0.5)
\rput(1,0.45){\scriptsize{\emph{1}}}\rput(1,-0.5){\scriptsize{\emph{3}}}\rput(1.2,0.15){\scriptsize{\emph{2}}}
\psline{-}(1.5,0.5)(2.5,0.5)\psdot(2.5,0.5)\psline{-,linestyle=dashed}(2.5,0.5)(3,0.5)
\rput(2,0.7){\scriptsize{\emph{3}}}\psline{-}(1.5,0)(2.5,0)\psdot(2.5,0)\psline{-,linestyle=dashed}(2.5,0)(3,0)
\rput(2,0.2){\scriptsize{\emph{1}}}\psline{-}(1.5,-0.5)(2.5,-0.5)\psdot(2.5,-0.5)\psline{-,linestyle=dashed}(2.5,-0.5)(3,-0.5)
\rput(2,-0.7){\scriptsize{\emph{1}}}
\end{pspicture}
\begin{pspicture}(0,-1)(4.9,1)
\rput(0.2,0.8){\textbf{III9}}\psline{-,linestyle=dashed}(0,0)(0.5,0)
\psdot(0.5,0) \rput(0.45,0.3){\footnotesize{$Q_i$}}
\psline{-}(0.5,0)(1.5,0.5) \psdot(1.5,0.5)\psline{-}(0.5,0)(1.5,0)
\psdot(1.5,0) \psline{-}(0.5,0)(1.5,-0.5) \psdot(1.5,-0.5)
\rput(1,0.45){\scriptsize{\emph{1}}}\rput(1,-0.5){\scriptsize{\emph{3}}}\rput(1.2,0.15){\scriptsize{\emph{2}}}
%\psline{-}(1.5,0.5)(2.5,0.5)\psdot(2.5,0.5)\psline{-,linestyle=dashed}(2.5,0.5)(3,0.5)
%\rput(2,0.7){\scriptsize{\emph{3}}}
\psline{-}(1.5,0)(2.5,0)\psdot(2.5,0)\psline{-,linestyle=dashed}(2.5,0)(3,0)
\rput(2,0.2){\scriptsize{\emph{1}}}\psline{-}(1.5,-0.5)(2.5,-0.5)\psdot(2.5,-0.5)\psline{-,linestyle=dashed}(2.5,-0.5)(3,-0.5)
\rput(2,-0.7){\scriptsize{\emph{1}}} \psline{-}(1.5,0.5)(2.5,1)
\psdot(2.5,1) \psline{-}(1.5,0.5)(2.5,0.5)
\psdot(2.5,0.5)\psline{-,linestyle=dashed}(2.5,1)(3,1)
\rput(2,0.95){\scriptsize{\emph{2}}}
\rput(2.3,0.7){\scriptsize{\emph{3}}}\psline{-,linestyle=dashed}(2.5,0.5)(3,0.5)
\end{pspicture}
\begin{center}
\emph{List 1}
\end{center}
In the next step we give an alternative description for
$\chi(E_i^{\circ})$. We first introduce some new notation.
\begin{notation}
\emph{We write $D:=m_i^2- \sum_{j \rightarrow i} m_j^2$ and
$r_{ab}:=m_i - M_{Q_i}(a,b) - M_{Q_i}(b,a)$ for $a,b \in
\{1,2,3\}=\{a,b,c\}$, $a \neq b$. Let $R$ be equal to $\hat{r_{12}}
+ \hat{r_{13}} + \hat{r_{23}}$ where $\hat{r_{ab}}:=$\[
\begin{cases}
\begin{array}{cl}
r_{ab} & \mbox{ if label $c$ does not appear
under $Q_i$;} \\
%\nexists l: [(Q_i \rightarrow Q_l) \mbox{ and } (Q_i(a) \rightarrow Q_l) \mbox{ and } (Q_i(b) \rightarrow Q_l) ]      \\
0 & \mbox{ else}.
\end{array}
\end{cases} \]
We refer to the beginning of Section $5$ for the definition of $T$
and to Table $1$ for the values that $T$ takes.}
\end{notation}
%Although the multiplicities of the points $Q_i(a)$ and $Q_i(b)$ can
%be equal to $0$, we interpret these points as points of the
%constellation. Hence, $\hat{r_{ab}}$ can be $0$ although the
%multiplicity of $Q_i(a)$ is equal to $0$ and $r_{ab} \neq 0$.
\begin{lemma}\label{lemmadrt}
\[\chi(E_i^{\circ}) = D - R + T.\]
\end{lemma}
${}$\\
\emph{Proof.} \quad We will prove that \begin{eqnarray} R=3m_i-2
\sum_{j \rightarrow i}m_j + \sum_{j \rightarrow i}(\sum_{\substack{k
\rightarrow i\\ k \twoheadrightarrow j}} m_k) + \sum_{i \rightarrow
j}(\sum_{\substack{k \rightarrow j\\ k \twoheadrightarrow i}} m_k) -
\sum_{i \rightarrow j}m_i.\end{eqnarray} Let $X$ be the right hand
side in $(3)$, let $X_1:=\sum_{j \rightarrow i}(\sum_{\substack{k
\rightarrow i\\ k \twoheadrightarrow j}} m_k)$ and $X_2:=\sum_{i
\rightarrow j}(\sum_{\substack{k \rightarrow j\\ k
\twoheadrightarrow i}} m_k)$. For $k \rightarrow i$, one has one of
the following situations.
\begin{itemize}
\item There exist exactly two points $Q_{j_1}$ and $Q_{j_2}$ that are proximate to $Q_i$ and
for which $k \twoheadrightarrow j_1$ and $k \twoheadrightarrow j_2$.
Then $m_k$ appears twice as term in $X_1$ and $Q_k$ is not linearly
proximate to $Q_i$, hence $m_k$ does not appear in $X_2$. This
implies that $m_k$ does not show up in $X$.
\item There exists exactly one point $Q_j$ that is proximate to
$Q_i$ and for which $k \twoheadrightarrow j$. We are in the
following situation:\\
\begin{pspicture}(-2,-0.5)(8,0.5)
\psline{-,linestyle=dashed}(0,0)(0.5,0) \psdot(0.5,0)
\rput(0.5,0.3){\footnotesize{$Q_i$}} \psline{-}(0.5,0)(1.5,0)
\psdot(1.5,0) \psline{-, linestyle=dashed}(1.5,0)(6,0)
\rput(1,0.2){\scriptsize{\emph{1}}} \psdot(6,0)
\psline{-}(6,0)(7,0)\psdot(7,0)
\psline{-,linestyle=dashed}(7,0)(7.5,0)
\rput(7,0.3){\footnotesize{$Q_k$}}
\rput(6,0.3){\footnotesize{$Q_j$}}\rput(6.5,0.2){\scriptsize{\emph{2}}}
\rput(4,0.2){\scriptsize{\emph{2}}}
\end{pspicture}
\\
Then $m_k$ appears once in $X_1$ and $k \twoheadrightarrow i$. If
label $3$ appears under $Q_i$, then $m_k$ appears once in $X_2$.
Hence, $m_k$ does not show up in the expression $X$. If there is no
label $3$ under $Q_i$, then $m_k$ does not appear in $X_2$ such that
this $m_k$ appears with coefficient $-1$ in $X$.
\item There exists no point $Q_j$ such that $j \rightarrow i$ and $k
\twoheadrightarrow j$. Then $k \succ i$ and $m_k$ does not appear in
$X_1$. The number of times that $m_k$ appears in $X_2$ depends on
the labels below $Q_i$. It can be once, twice or three times.
\end{itemize}
Notice that the multiplicities $m_k$ of the points $Q_k$ with $k
\rightarrow i$ but not $k \twoheadrightarrow i$ do not appear in
$X$. To analyse further the formula $X$, we now take the labels into
account that appear below $Q_i$. \\ If $Q_i$ is the origin, then the
points $Q_k$ for which $k\succ i$ appear with coefficient $-2$ in
$X$. The other points $Q_j$ for which $j \twoheadrightarrow i$ have
coefficient $-1$. Hence
\[X=3m_i-\sum_{k \twoheadrightarrow i}m_k - \sum_{k \succ i} m_k = r_{12}+r_{13}+r_{23}=R.\]
Also in the other cases one can check that $X=R$: when $1$ is the
only label below $Q_i$, then $X=r_{12}+r_{13}=R$. If the labels
showing up below $Q_i$ are $1$ and $3$, then $X=r_{13}=R$. If the
three labels show up under $Q_i$, then $X=0=R$.\hfill $\blacksquare$
\\ \\
Notice that it follows from the linear proximity relations that $R
\geq 0$. Formula $(1)$ in Subsection $2.3$ shows that $D \geq 0$ and
from Table $1$ it follows that $0 \leq T \leq 3$. In order to find
the cases where $\chi(E_i^{\circ}) \leq 0$, we will investigate when
$R \geq D$. We want to give an estimation for $D$. In particular, we
will determine a lower bound $L$ for $D$ and we will then check when
$R \geq L$. We introduce some terminology.
\begin{definition}
Let $l \in \mathbb{Z}_{> 1}$ and let $n_1,\cdots, n_l, h_1, \cdots,
h_{l-1} \in \mathbb{Z}_{>0}$ such that $n_j = h_j n_{j+1}+n_{j+2}$
where $0 < n_{j+2} < n_{j+1}$, for $1 \leq j \leq l-2$, and such
that $n_{l-1}=h_{l-1}n_l$. If $l$ is even, then set $(a,b)=(3,2)$.
If $l$ is odd, we set $(a,b)=(2,3)$. Let $\mathcal{A}$ be an
idealistic cluster
\\ \begin{pspicture}(0.1,-0.5)(13.3,0.5)
\psline{-,linestyle=dashed}(0.1,0)(0.5,0) \psline{-}(0.5,0)(1.3,0)
\psdot(0.5,0)\psdot(1.3,0) \psline{-}(1.3,0)(2.1,0) \psdot(2.1,0)
\rput(0.9,0.2){\scriptsize{\emph{1}}}\rput(1.7,0.2){\scriptsize{\emph{2}}}
\rput(3,0.2){\scriptsize{\emph{2}}}\rput(3.8,0.2){\scriptsize{\emph{3}}}
\rput(5.4,0.2){\scriptsize{\emph{3}}}\rput(6.2,0.2){\scriptsize{\emph{2}}}
\rput(7.8,0.2){\scriptsize{\emph{2}}}\rput(8.6,0.2){\scriptsize{\emph{3}}}
\rput(10.6,0.2){\scriptsize{a}}\rput(11.4,0.2){\scriptsize{b}}
\rput(12.8,0.2){\scriptsize{b}}
\psline{-,linestyle=dashed}(2.1,0)(2.6,0) \psdot(2.6,0)
\psline{-}(2.6,0)(3.4,0)\psdot(3.4,0)
\psline{-}(3.4,0)(4.2,0)\psdot(4.2,0)
\psline{-,linestyle=dashed}(4.2,0)(5,0) \psdot(5,0)
\psline{-}(5,0)(5.8,0)\psdot(5.8,0) \psline{-}(5.8,0)(6.6,0)
\psdot(6.6,0)\psline{-,linestyle=dashed}(6.6,0)(7.4,0) \psdot(7.4,0)
\psline{-}(7.4,0)(8.2,0) \psdot(8.2,0) \psline{-}(8.2,0)(9,0)
\psdot(9,0)\psline{-,linestyle=dashed}(9,0)(10.2,0) \psdot(10.2,0)
\psline{-}(10.2,0)(11,0) \psdot(11,0) \psline{-}(11,0)(11.8,0)
\psdot(11.8,0) \psline{-,linestyle=dashed}(11.8,0)(12.4,0)
\psdot(12.4,0) \psline{-}(12.4,0)(13.2,0) \psdot(13.2,0)
\rput(0.6,0.3){\footnotesize{$Q_i$}}
\rput(0.6,-0.3){\footnotesize{$n_1$}}
\rput(1.3,-0.3){\footnotesize{$n_2$}}
\rput(2.1,-0.3){\footnotesize{$n_2$}}
\rput(2.6,-0.3){\footnotesize{$n_2$}}
\rput(3.4,-0.3){\footnotesize{$n_3$}}
\rput(4.2,-0.3){\footnotesize{$n_3$}}
\rput(5,-0.3){\footnotesize{$n_3$}}
\rput(5.8,-0.3){\footnotesize{$n_4$}}
\rput(6.6,-0.3){\footnotesize{$n_4$}}
\rput(7.4,-0.3){\footnotesize{$n_4$}}
\rput(8.2,-0.3){\footnotesize{$n_5$}}
\rput(9,-0.3){\footnotesize{$n_5$}}
\rput(10.2,-0.3){\footnotesize{$n_{l-1}$}}
\rput(11,-0.3){\footnotesize{$n_l$}}
\rput(11.8,-0.3){\footnotesize{$n_l$}}
\rput(12.4,-0.3){\footnotesize{$n_l$}}
\rput(13.2,-0.3){\footnotesize{$n_l$}}
\end{pspicture}
\\ where $n_{j}$ appears $h_{j-1}$ consecutive times, $2 \leq j \leq l$.
We call $\mathcal{A}$ a \emph{Euclidean cluster starting in $Q_i$}.
\end{definition}
\begin{definition}
Let $\mathcal{A}$ be a cluster of the form\\ \\
\begin{pspicture}(-3,-1)(7,1)
\psline{-,linestyle=dashed}(-0.5,0)(0.5,0)
\psdot(0.5,0)\psline{-}(0.5,0)(1.5,0) \psdot(1.5,0)
\psline{-}(1.5,0)(2.5,0.5) \psdot(2.5,0.5)
\psline{-}(1.5,0)(2.5,-0.5) \psdot(2.5,-0.5)
\rput(1,0.2){\scriptsize{\emph{1}}}\rput(2,0.45){\scriptsize{\emph{2}}}\rput(2,-0.5){\scriptsize{\emph{3}}}
\rput(3,0.7){\scriptsize{\emph{3}}}
\rput(5,0.7){\scriptsize{\emph{3}}}\rput(6,0.7){\scriptsize{\emph{2}}}
\rput(3,-0.7){\scriptsize{\emph{2}}}
\rput(5,-0.7){\scriptsize{\emph{2}}}\rput(6,-0.7){\scriptsize{\emph{3}}}
\psline{-}(2.5,0.5)(3.5,0.5) \psdot(3.5,0.5)
\psline{-,linestyle=dashed}(3.5,0.5)(4.5,0.5) \psdot(4.5,0.5)
\psline{-}(4.5,0.5)(5.5,0.5)\psdot(5.5,0.5)
\psline{-,linestyle=dashed}(5.5,0.5)(6.5,0.5)
\psline{-}(2.5,-0.5)(3.5,-0.5) \psdot(3.5,-0.5)
\psline{-,linestyle=dashed}(3.5,-0.5)(4.5,-0.5) \psdot(4.5,-0.5)
\psline{-}(4.5,-0.5)(5.5,-0.5)\psdot(5.5,-0.5)
\psline{-,linestyle=dashed}(5.5,-0.5)(6.5,-0.5)
 \rput(0.5,0.3){\footnotesize{$m_1$}}
\rput(0.5,-0.3){\footnotesize{$Q_i$}}
 \rput(1.5,0.3){\footnotesize{$m_2$}}
  \rput(2.5,0.8){\footnotesize{$n_1$}} \rput(3.5,0.8){\footnotesize{$n_2$}} \rput(4.5,0.8){\footnotesize{$n_{l-1}$}}
 \rput(5.5,0.8){\footnotesize{$n_l$}}
 \rput(2.5,-0.8){\footnotesize{$n'_1$}} \rput(3.5,-0.8){\footnotesize{$n'_2$}} \rput(4.5,-0.8){\footnotesize{$n'_{r-1}$}}
 \rput(5.5,-0.8){\footnotesize{$n'_r$}}
\end{pspicture}\\
such that
\\
\begin{pspicture}(0,0)(7,1.5)
\psline{-,linestyle=dashed}(0.1,0.5)(0.5,0.5)
\psdot(0.5,0.5)\psline{-}(0.5,0.5)(1.5,0.5) \psdot(1.5,0.5)
\psline{-}(1.5,0.5)(2.5,0.5) \psdot(2.5,0.5)
\rput(1,0.7){\scriptsize{\emph{1}}}\rput(2,0.7){\scriptsize{\emph{2}}}\rput(3,0.7){\scriptsize{\emph{3}}}
\rput(5,0.7){\scriptsize{\emph{3}}}\rput(5.8,0.7){\scriptsize{\emph{2}}}
\psline{-}(2.5,0.5)(3.5,0.5) \psdot(3.5,0.5)
\psline{-,linestyle=dashed}(3.5,0.5)(4.5,0.5) \psdot(4.5,0.5)
\psline{-}(4.5,0.5)(5.5,0.5)\psdot(5.5,0.5)
\psline{-,linestyle=dashed}(5.5,0.5)(6,0.5)
 \rput(0.5,0.8){\footnotesize{$M_1$}} \rput(1.5,0.8){\footnotesize{$M_2$}}
  \rput(2.5,0.8){\footnotesize{$n_1$}} \rput(3.5,0.8){\footnotesize{$n_2$}} \rput(4.5,0.8){\footnotesize{$n_{l-1}$}}
 \rput(5.5,0.8){\footnotesize{$n_l$}}\rput(6.5,0.5){\mbox{ and}}
\end{pspicture}
\begin{pspicture}(0,0)(7,1.5)
\psline{-,linestyle=dashed}(0.1,0.5)(0.5,0.5)
\psdot(0.5,0.5)\psline{-}(0.5,0.5)(1.5,0.5) \psdot(1.5,0.5)
\psline{-}(1.5,0.5)(2.5,0.5) \psdot(2.5,0.5)
\rput(1,0.7){\scriptsize{\emph{1}}}\rput(2,0.7){\scriptsize{\emph{3}}}\rput(3,0.7){\scriptsize{\emph{2}}}
\rput(5,0.7){\scriptsize{\emph{2}}}\rput(5.8,0.7){\scriptsize{\emph{3}}}
\psline{-}(2.5,0.5)(3.5,0.5) \psdot(3.5,0.5)
\psline{-,linestyle=dashed}(3.5,0.5)(4.5,0.5) \psdot(4.5,0.5)
\psline{-}(4.5,0.5)(5.5,0.5)\psdot(5.5,0.5)
\psline{-,linestyle=dashed}(5.5,0.5)(6,0.5)
 \rput(0.5,0.8){\footnotesize{$M'_1$}} \rput(1.5,0.8){\footnotesize{$M'_2$}}
  \rput(2.5,0.8){\footnotesize{$n'_1$}} \rput(3.5,0.8){\footnotesize{$n'_2$}} \rput(4.5,0.8){\footnotesize{$n'_{r-1}$}}
 \rput(5.5,0.8){\footnotesize{$n'_r$}}
\end{pspicture}\\
are Euclidean clusters, where $M_1:=m_1-\sum_{j=1}^r n'_j$,
$M_2:=m_2-\sum_{j=1}^r n'_j$ and $M'_1:=m_1-\sum_{j=1}^l n_j$,
$M'_2:=m_2-\sum_{j=1}^l n_j$. We call the cluster $\mathcal{A}$ a
\emph{bi-Euclidean cluster starting in  $Q_i$}.
\end{definition}
\begin{example}
\emph{The cluster
\\ \begin{pspicture}(0.1,-0.5)(13.3,0.5)
\psline{-}(3,0)(3.8,0) \psdot(3,0)\psdot(3.8,0)
\psline{-}(3.8,0)(4.6,0) \psdot(4.6,0)
\rput(3.4,0.2){\scriptsize{\emph{1}}}\rput(4.2,0.2){\scriptsize{\emph{2}}}
\rput(5,0.2){\scriptsize{\emph{2}}}\rput(5.8,0.2){\scriptsize{\emph{2}}}
\rput(6.6,0.2){\scriptsize{\emph{3}}}\rput(7.4,0.2){\scriptsize{\emph{2}}}
\rput(8.2,0.2){\scriptsize{\emph{2}}}\rput(9,0.2){\scriptsize{\emph{2}}}
\psline{-}(4.6,0)(5.4,0) \psdot(5.4,0)
\psline{-}(5.4,0)(6.2,0)\psdot(6.2,0)
\psline{-}(6.2,0)(7,0)\psdot(7,0) \psline{-}(7,0)(7.8,0)
\psdot(7.8,0) \psline{-}(7.8,0)(8.6,0)\psdot(8.6,0)
\psline{-}(8.6,0)(9.4,0) \psdot(9.4,0)
\rput(3,0.3){\footnotesize{$Q_i$}}
\rput(3,-0.3){\footnotesize{$19$}}
\rput(3.8,-0.3){\footnotesize{$5$}}
\rput(4.6,-0.3){\footnotesize{$5$}}
\rput(5.4,-0.3){\footnotesize{$5$}}
\rput(6.2,-0.3){\footnotesize{$4$}}
\rput(7,-0.3){\footnotesize{$1$}}
\rput(7.8,-0.3){\footnotesize{$1$}}
\rput(8.6,-0.3){\footnotesize{$1$}}
\rput(9.4,-0.3){\footnotesize{$1$}}
\end{pspicture}
\\ is a Euclidean cluster starting in $Q_i$.
The cluster \\
\begin{pspicture}(-3,-1)(7,1)
\psdot(-0.5,0)\psline{-}(-0.5,0)(0.5,0)
\psdot(0.5,0)\psline{-}(0.5,0)(1.5,0) \psdot(1.5,0)
\psline{-}(1.5,0)(2.5,0.5) \psdot(2.5,0.5)
\psline{-}(1.5,0)(2.5,-0.5) \psdot(2.5,-0.5)
\rput(0,0.2){\scriptsize{\emph{2}}}\rput(1,0.2){\scriptsize{\emph{1}}}\rput(2,0.45){\scriptsize{\emph{2}}}\rput(2,-0.5){\scriptsize{\emph{3}}}
\rput(3,0.7){\scriptsize{\emph{3}}}\rput(4,0.7){\scriptsize{\emph{2}}}
\rput(5,0.7){\scriptsize{\emph{2}}}\rput(6,0.7){\scriptsize{\emph{3}}}
\psline{-}(2.5,0.5)(3.5,0.5) \psdot(3.5,0.5)
\psline{-}(3.5,0.5)(4.5,0.5) \psdot(4.5,0.5)
\psline{-}(4.5,0.5)(5.5,0.5)\psdot(5.5,0.5)
\psline{-}(5.5,0.5)(6.5,0.5)
 \rput(-0.5,0.3){\footnotesize{$88$}}
 \rput(0.5,0.3){\footnotesize{$17$}}
\rput(0.5,-0.3){\footnotesize{$Q_i$}}
 \rput(1.5,0.3){\footnotesize{$12$}}
  \rput(2.5,0.8){\footnotesize{$5$}} \rput(3.5,0.8){\footnotesize{$2$}} \rput(4.5,0.8){\footnotesize{$2$}}
 \rput(5.5,0.8){\footnotesize{$1$}}\rput(6.5,0.8){\footnotesize{$1$}}
 \rput(2.5,-0.8){\footnotesize{$5$}}
\end{pspicture}\\
is a bi-Euclidean cluster starting in $Q_i$.}\hfill $\square$
\end{example}%\vspace*{-1cm}
%\[\]
\begin{definition}
Suppose that $Q$ is a point different from the origin in a
$3$-dimen\-sio\-nal toric constellation $\mathcal{C}$. Let $a \in
\{1,2,3\}$ such that $Q=P(a)$ for a point $P \in \mathcal{C}$ and
suppose that there exists $b \in \{1,2,3\}, a \neq b$, such that
$Q(b) \in \mathcal{C}$. Then we call $Q$ a \emph{switch point}.
\end{definition}
\begin{proposition} \label{lemmacomb}
Let $\mathcal{A}=(\mathcal{C},\underline{m})$ be a $3$-dimensional
toric idealistic cluster. Let $Q_i \in \mathcal{C}$ and suppose that
there exists exactly one point $Q_k \in \mathcal{C}$ for which $k
\succ i$.
%\begin{pspicture}(0,0)(2,0.5)
%\psdot(0.5,0) \rput(0.5,0.3){$Q_i$} \psline{-}(0.5,0)(1.5,0)
%\psdot(1.5,0) \rput(1.5,0.3){$Q_k$} \rput(1,0.2){\emph{1}}}
%\psline{-,linestyle=dashed}(1.5,0)(2,0)
%\end{pspicture}\\
Then the following properties hold:
\begin{enumerate}
\item $m_i m_k  \geq \sum_{j \rightarrow i}m_j^2$;
\item $m_i m_k  = \sum_{j \rightarrow i}m_j^2$ if and only
if $\mathcal{A}$ is a Euclidean cluster or $\mathcal{A}$ is a
bi-Euclidean cluster starting in $Q_i$.
\end{enumerate}
\end{proposition}
${}$ \emph{Proof.} \quad CASE $1$: there exists at most one point
$Q_l$ in $\mathcal{C}$ that is proximate to $Q_i$ and such that $Q_l
\succ Q_k$. Then we can suppose that the cluster is of the form:
\begin{pspicture}(-0.5,0)(2,0.5)
\psline{-,linestyle=dashed}(0,0)(0.5,0) \psdot(0.5,0)
\rput(0.5,0.3){\footnotesize{$m_i$}}
\rput(0.5,-0.3){\footnotesize{$Q_i$}} \psline{-}(0.5,0)(1.5,0)
\psdot(1.5,0) \rput(1.5,0.3){\footnotesize{$m_k$}}
\rput(1.5,-0.3){\footnotesize{$Q_k$}}\rput(1,0.2){\scriptsize{\emph{1}}}
\rput(2,0.2){\scriptsize{\emph{2}}}\psline{-,linestyle=dashed}(1.5,0)(2.5,0)
\psdot(2.5,0) \psline{-,linestyle=dashed}(2.5,0)(3,0)
\end{pspicture}\\ \\
We have $m_i \geq M_{Q_i}(1,2)$ and thus \begin{eqnarray} m_k m_i
\geq \sum_{t \geq 0} m_k m_{Q_i(1,2^t)}.\end{eqnarray} We give lower
bounds for the terms $m_k m_{Q_i(1,2^t)}$ in $(4)$ depending on
whether $Q_i(1,2^t)$ is a switch point or not. If $Q_i(1,2^t)=Q_k$,
then $m_k m_{Q_i(1,2^t)} = m_k^2$. If $P^t:=Q_i(1,2^t)$, $t \neq 0$,
is a switch point, then
\[m_k m_{P^t} \geq
\left(\sum_{s \geq 0}m_{P^{t-1}(2,3^s)}+\sum_{s \geq
0}m_{P^{t-1}(3,2^s)} \right)m_{P^t}.\]
\\If $P^t$ is not a switch point, then we estimate $m_k m_{P^t} \geq
m_{P^t} m_{P^t}$. We fill in these lower bounds in $(4)$ and we get
\begin{eqnarray*}
m_k m_i \geq m_k^2 + \sum_{ \substack{ P^t \text{ not} \\
\text{switch
point}}}m_{P^t}^2 + \sum_{ \substack{ P^t \text{ switch} \\
\text{point, } t \neq 0}}\left(\sum_{s \geq
0}m_{P^{t-1}(2,3^s)}+\sum_{s \geq 0}m_{P^{t-1}(3,2^s)}
\right)m_{P^t}.
\end{eqnarray*}
We iterate this process: whenever we have a product $m_{Q_j}
m_{Q_l}$ with $j < l$, we use the estimations described above for
$m_{Q_j} m_{Q_l}$ according to whether $Q_l$ is a switch point or
not. I.e., if $Q_l$ is a switch point and if $P \in \mathcal{C}$ is
such that $Q_l \succ P$, then set $m_{Q_j} m_{Q_l} \geq (\sum_{t
\geq 0} m_{P(2,3^t)} + \sum_{t \geq 0} m_{P(3,2^t)} )m_{Q_l}$. If
$Q_l$ is not a switch point, then we set $m_{Q_j} m_{Q_l} \geq
m_{Q_l} m_{Q_l}$. This is obviously a finite process and it shows
that
\[m_i m_k  \geq \sum_{j \rightarrow i}m_j^2.\]
We now study when $\sum_{j \rightarrow i}m_j^2=m_i m_k$.
\begin{itemize}
\item If $\mathcal{C}$ is a chain, then it is not difficult to see that $\sum_{j \rightarrow i}m_j^2=m_i
m_k$ if and only if $\mathcal{A}$ is a Euclidean cluster.
\item If $\mathcal{A}$ contains a subcluster of the form
\\
\begin{pspicture}(0,0)(7,1.5)
\psline{-,linestyle=dashed}(0,0.5)(0.5,0.5)
\psdot(0.5,0.5)\psline{-}(0.5,0.5)(1.5,0.5) \psdot(1.5,0.5)
\psline{-,linestyle=dashed}(1.5,0.5)(4.5,0.5) \psdot(4.5,0.5)
\psdot(5.5,0.5)\psline(4.5,0.5)(5.5,0.5)
\rput(1,0.2){\scriptsize{\emph{1/3}}}\rput(3,0.2){\scriptsize{\emph{2}}}
\rput(5,0.2){\scriptsize{\emph{\textbf{2}}}}
\psline{-}(5.5,0.5)(6.5,1)\psdot(6.5,1) \psline{-}(6.5,1)(7.5,1)
\psdot(7.5,1)\psline{-,linestyle=dashed}(7.5,1)(8,1)
\psline{-}(5.5,0.5)(6.5,0)\psdot(6.5,0)\psline{-,linestyle=dashed}(6.5,0)(7,0)
\rput(6,0.95){\scriptsize{\emph{\textbf{2}}}}\rput(7,1.2){\scriptsize{\emph{\textbf{2}}}}
\rput(6,0){\scriptsize{\emph{3}}}
 \rput(1.5,0.8){\footnotesize{$P$}}
  \rput(5.3,0.8){\footnotesize{$P(2^s)$}}
\end{pspicture}\\
where $P(2^{s+1})$ is not a switch point, then at some moment in the
process we get
\begin{eqnarray*}
m_i m_k & \geq & \cdots + m_P\left(\sum_{t \geq 0}m_{P(2^t)}\right) \\
& > & \cdots + \sum_{t=0}^{s-1}m_{P(2^t)}^2 +
m_{P(2^s)}\left(\sum_{t \geq 0}m_{P(2^s,3^t)}\right) +
\sum_{t=s+1}m_{P(2^{t})}^2.
\end{eqnarray*}
Indeed, $m_P > m_{P(2^{s+1})}$.
\item If $\mathcal{A}$ contains a subcluster of the form
\\
\begin{pspicture}(0,0)(7,1.5)
\psline{-,linestyle=dashed}(0,0.5)(0.5,0.5)
\psdot(0.5,0.5)\psline{-}(0.5,0.5)(1.5,0.5) \psdot(1.5,0.5)
\psline{-,linestyle=dashed}(1.5,0.5)(4.5,0.5) \psdot(4.5,0.5)
\psdot(5.5,0.5)\psline(4.5,0.5)(5.5,0.5)
\rput(1,0.2){\scriptsize{\emph{1/3}}}\rput(3,0.2){\scriptsize{\emph{2}}}
\rput(5,0.2){\scriptsize{\emph{\textbf{2}}}}
\psline{-}(5.5,0.5)(6.5,1)\psdot(6.5,1) \psline{-}(6.5,1)(7.5,1)
\psdot(7.5,1)\psline{-,linestyle=dashed}(7.5,1)(8,1)
\psline{-}(5.5,0.5)(6.5,0)\psdot(6.5,0)\psline{-,linestyle=dashed}(6.5,0)(7,0)
\rput(6,0.95){\scriptsize{\emph{\textbf{2}}}}\rput(7,1.2){\scriptsize{\emph{\textbf{3}}}}
\rput(6,0){\scriptsize{\emph{3}}}
 \rput(1.5,0.8){\footnotesize{$P$}}
  \rput(5.3,0.8){\footnotesize{$P(2^s)$}} \rput(6.6,0.3){\footnotesize{$Q$}}
\end{pspicture}\\
then at some moment in the process we get
\begin{eqnarray*}
m_i m_k & \geq & \cdots + m_P\left(\sum_{t \geq 0}m_{P(2^t)}\right) \\
& \geq & \cdots + \sum_{t=0}^{s-1}m_{P(2^t)}^2 +
m_{P(2^s)}\left(\sum_{t \geq 0}m_{P(2^s,3^t)}\right) + \\ & &
m_{P(2^{s+1})} \left(\sum_{t \geq 0}m_{P(2^{s+1},3^t)}+\sum_{t \geq
0}m_{Q(2^t)}\right)
\\
& \geq & \sum_{j \rightarrow i}m_j^2 +
m_{P(2^{s+1})}\sum_{t \geq 0}m_{Q(2^t)}\\
& > &  \sum_{j \rightarrow i}m_j^2.
\end{eqnarray*}
\item If $\mathcal{A}$ contains a subcluster of the form
\\
\begin{pspicture}(0,0)(7,1.5)
\psline{-,linestyle=dashed}(0,0.5)(0.5,0.5)
\psdot(0.5,0.5)\psline{-}(0.5,0.5)(1.5,0.5) \psdot(1.5,0.5)
%\psline{-}(1.5,0.5)(2.5,0.5)
\psline{-,linestyle=dashed}(1.5,0.5)(4.5,0.5) %\psdot(2.5,0.5)
\psdot(4.5,0.5) \psdot(5.5,0.5)\psline(4.5,0.5)(5.5,0.5)
\rput(1,0.2){\scriptsize{\emph{1/3}}}\rput(3,0.2){\scriptsize{\emph{2}}}
\rput(5,0.2){\scriptsize{\emph{\textbf{2}}}}
\psline{-}(5.5,0.5)(6.5,1)\psdot(6.5,1)
\psline{-}(5.5,0.5)(6.5,0)\psdot(6.5,0)\psline{-,linestyle=dashed}(6.5,0)(7,0)
\rput(6,0.95){\scriptsize{\emph{\textbf{2}}}}
\rput(6,0){\scriptsize{\emph{3}}}
 \rput(1.5,0.8){\footnotesize{$P$}}
  \rput(5.3,0.8){\footnotesize{$P(2^s)$}}
\end{pspicture}\\
then at some moment in the process we get
\begin{eqnarray*}
m_i m_k & \geq & \cdots + m_P\left(\sum_{t \geq 0}m_{P(2^t)}\right) \\
& > & \cdots + \sum_{t=0}^{s-1}m_{P(2^t)}^2 + m_{P(2^s)}\sum_{t \geq
0}m_{P(2^s,3^t)} + m_{P(2^{s+1})}^2.
\end{eqnarray*}
Indeed, $m_P > m_{P(2^{s+1})}$.
\end{itemize}
CASE $2$: there exist two points $Q_a$ and $Q_b$ in $\mathcal{C}$
that are proximate to $Q_i$ and such that $Q_a \succ Q_k$ and $Q_b
\succ Q_k$. Then we may suppose that the cluster is of the form:
\\ \begin{pspicture}(-3,-1)(4.9,1)
\psline{-,linestyle=dashed}(0,0)(0.5,0)\psdot(0.5,0)
\rput(0.5,0.3){\footnotesize{$m_i$}}
\rput(0.5,-0.3){\footnotesize{$Q_i$}}\psline{-}(0.5,0)(1.5,0)
\rput(1.5,-0.3){\footnotesize{$Q_k$}} \psdot(1.5,0)
\psline{-}(1.5,0)(2.5,0.5) \psdot(2.5,0.5)
\psline{-}(1.5,0)(2.5,-0.5) \psdot(2.5,-0.5)
\rput(1,0.2){\scriptsize{\emph{1}}}\rput(2,0.45){\scriptsize{\emph{2}}}\rput(2,-0.5){\scriptsize{\emph{3}}}
\psline{-,linestyle=dashed}(2.5,0.5)(3,0.5)
\psline{-,linestyle=dashed}(2.5,-0.5)(3,-0.5)\rput(1.5,0.3){\footnotesize{$m_k$}}
\rput(2.5,0.8){\footnotesize{$n_1$}}
\rput(2.5,0.2){\footnotesize{$Q_a$}}\rput(2.5,-0.8){\footnotesize{$n'_1$}}
\rput(2.5,-0.2){\footnotesize{$Q_b$}}
\end{pspicture}
\\ \\
Define $t:=m_k - M_{Q_k}(2,3) - M_{Q_k}(3,2)$. As the cluster is
idealistic, $t \geq 0$. Then also the clusters
\begin{center}
\begin{pspicture}(0.4,0)(5,1.5)
\psdot(0.5,0.5)\psline{-}(0.5,0.5)(1.5,0.5) \psdot(1.5,0.5)
\psline{-}(1.5,0.5)(2.5,0.5) \psdot(2.5,0.5)
\rput(1,0.7){\scriptsize{\emph{1}}}\rput(2,0.7){\scriptsize{\emph{2}}}
\psline{-,linestyle=dashed}(2.5,0.5)(3,0.5)
 \rput(0.5,0.8){\footnotesize{$M_1$}}
 \rput(1.5,0.8){\footnotesize{$M_2$}}
  \rput(2.5,0.8){\footnotesize{$n_1$}} \rput(4,0.5){\mbox{and}}
\rput(0.5,0.2){\footnotesize{$Q_i$}}
\rput(1.5,0.2){\footnotesize{$Q_k$}}
\rput(2.5,0.2){\footnotesize{$Q_a$}}
\end{pspicture}
\begin{pspicture}(0,0)(5,1.5)
\psdot(0.5,0.5)\psline{-}(0.5,0.5)(1.5,0.5) \psdot(1.5,0.5)
\psline{-}(1.5,0.5)(2.5,0.5) \psdot(2.5,0.5)
\rput(1,0.7){\scriptsize{\emph{1}}}\rput(2,0.7){\scriptsize{\emph{3}}}
\psline{-,linestyle=dashed}(2.5,0.5)(3,0.5)
 \rput(0.5,0.8){\footnotesize{$M'_1$}}
 \rput(1.5,0.8){\footnotesize{$M'_2$}}
  \rput(2.5,0.8){\footnotesize{$n'_1$}}
  \rput(0.5,0.2){\footnotesize{$Q_i$}}
\rput(1.5,0.2){\footnotesize{$Q_k$}}
\rput(2.5,0.2){\footnotesize{$Q_b$}}
\end{pspicture}\end{center}
with $M_1:=m_i-M_{Q_k}(2,3)-t $ and $M_2:=M_{Q_k}(2,3)$,
$M'_1:=m_i-M_{Q_k}(3,2)-t $ and $M'_2:=M_{Q_k}(3,2)$ are idealistic.
They are clusters of the form as in Case $1$, therefore that we can
use the bound that we obtained there:
\begin{eqnarray*}
\sum_{j \rightarrow i}m_j^2 & \leq & M_1 M_2 + M'_1 M'_2 - M_2^2 -
{M'_2}^2 + m_k^2 \\
& = & (M_2 + M'_2)(m_i-m_k) + m_k^2 \\
& = & m_i m_k - t(m_i-m_k) \\
& \leq & m_i m_k.
\end{eqnarray*}
From the previous computations it follows that $\sum_{j \rightarrow
i}m_j^2= m_i m_k$ if and only if the cluster is a bi-Euclidean
cluster starting in $Q_i$.  \hfill $\blacksquare$
\\ \\ This combinatoric result is the key to determine the sign
of $\chi(E_i^{\circ})$.
 ${}$
\newpage  ${}$  \begin{center} \textsc{6. Determination of the sign of
$\chi(E_i^{\circ})$}
\end{center} ${}$\\
In this section we classify the irreducible exceptional components
$E_i, 1 \leq i \leq r$, that arise in the blowing up of some
$3$-dimensional toric idealistic cluster according to the sign of
$\chi(E_i^{\circ})$. As in Lemma \ref{lemmadrt} we write
$\chi(E_i^{\circ})$ as $D-R+T$. For the points $Q_i$ in the clusters
of List $1$, we give a lower bound $L$ for $D$. We will use very
frequently Proposition \ref{lemmacomb}. As upper bound for $R$ we
use that $R \leq r_{12}+r_{13}+r_{23}$. We will study for which
clusters
in List $1$ it then holds that $r_{12}+r_{13}+r_{23} \geq L$.  \\
\indent We mark the name of the constellation by a star if there
exists a cluster with that underlying constellation that yields
$\chi(E_i^{\circ}) \leq 0$. We
refer to Table $1$ for the values of $T$. \\
\\ Let us first make the following observation.
\begin{remark} \label{remarkhulp}
\emph{Suppose $\mathcal{A}=(\mathcal{C},\underline{m})$ is a
$3$-dimensional toric idealistic cluster. Let $Q_i$ be a point of
the constellation $\mathcal{C}$. We define a subconstellation
$S^i\mathcal{C}$ of $\mathcal{C}$ as follows: the origin of
$S^i\mathcal{C}$ is $Q_i$ and $Q_j \in S^i\mathcal{C}$ if and only
if $j \rightarrow i$ in $\mathcal{C}$ or $j=i$. Suppose now that
$Q_k \in S^i\mathcal{C}$, $Q_k \neq Q_i$. We define a cluster
$S^i_k\mathcal{C}=(S^i\mathcal{C},\underline{n})$ with underlying
constellation $S^i\mathcal{C}$: for each point $Q_j \in
S^i\mathcal{C}, j \neq k$, set its multiplicity $n_j:=m_j$ and set
$n_k:=m_k+1$. If $S^i_k\mathcal{C}$ is idealistic, then there always
exists an idealistic cluster
$\tilde{\mathcal{A}}=(\mathcal{C},\tilde{\underline{m}})$ that
contains $S^i_k\mathcal{C}$ as a subcluster. Blowing up the
constellation $\mathcal{C}$ of cluster $\tilde{\mathcal{A}}$ then
yields
\begin{eqnarray*}
\chi(\tilde{E_i^{\circ}}) & = & \tilde{D} - \tilde{R}+ \tilde{T} \\
& = & D - 2m_k-1- (R - x) + T\\
& = & \chi(E_i^{\circ}) - 2m_k-1 + x,
\end{eqnarray*}
where $x$ is equal to $0,1$ or $2$ depending on the constellation
$\mathcal{C}$. \\
It follows that $\chi(\tilde{E_i^{\circ}}) < \chi(E_i^{\circ})$.}
\hfill $\square$\end{remark}
%\vspace*{-1cm}
%\[\]
${}$ This will make it possible to simplify computations. Indeed, as
described above, when we let increase the values of the
multiplicities such that the cluster stays idealistic and if
$\chi(\tilde{E_i^{\circ}}) \geq 0$, then $\chi(E_i^{\circ})>0$. \\
\\
We now proceed to the classification. Firstly we investigate the
constellations of List $1$ where at most one edge is going out of
$Q_i$. Then we consider the ones where exactly two edges leave out
of $Q_i$. We treat constellation II11 and we draw conclusions about
the subconstellations of II11 if possible. We will have to
investigate constellation II7 separately and we then also get the
classification for the constellations II1 and II3. Studying
constellation III9 will be enough to classify the constellations
where three edges are going out of $Q_i$.
\\ \\
\begin{tabular}{p{2cm}p{10.5cm}}
\begin{pspicture}(0,-0.5)(4.9,1)
\rput(0.2,0.8){\textbf{01$^*$}}
\psline{-,linestyle=dashed}(0,0)(0.5,0) \psdot(0.5,0)
\rput(0.5,0.3){\footnotesize{$Q_i$}}
\rput(0.5,-0.3){\footnotesize{$m_i$}}
\end{pspicture}
& \vspace*{-0.9cm} $D = m_i^2$. Does there exist a positive integer
$m_i$ such that $3m_i \geq m_i^2$?
\end{tabular}
\\ \\
We study the exact value of $\chi({E_i^{\circ}})$ if $m_i \in
\{1,2,3\}$. If $Q_i$ is the origin, then $T=3$ and
$\chi({E_i^{\circ}})=m_i^2-3m_i+3>0$. If there exists exactly one
point $Q_j$ such that $i \rightarrow j$, then $T=1$ and
$\chi({E_i^{\circ}})=m_i^2-2m_i+1$. We find that
\underline{$\chi({E_i^{\circ}})=0$} if $m_i=1$. If there exist
exactly two points in the constellation to which $Q_i$ is proximate,
then $T=0$ and $\chi({E_i^{\circ}})=m_i^2-m_i$. We find again that
\underline{$\chi({E_i^{\circ}})=0$} if $m_i=1$. If there are three
points to which $Q_i$ is proximate, then $T=0$ and
$\chi({E_i^{\circ}})=m_i^2>0$.
\\ \\
\begin{tabular}{p{2.5cm}p{10cm}}
\begin{pspicture}(0,-0.5)(4.9,1)
\rput(0.2,0.8){\textbf{I1$^*$}}
\psline{-,linestyle=dashed}(0,0)(0.5,0) \psdot(0.5,0)
\rput(0.5,0.3){\footnotesize{$Q_i$}} \psline{-}(0.5,0)(1.5,0)
\psdot(1.5,0) \rput(1,0.2){\scriptsize{\emph{1}}}
\rput(0.5,-0.3){\footnotesize{$m_i$}}
\rput(1.5,-0.3){\footnotesize{$m_1$}}
\end{pspicture}
& \vspace*{-1.2cm} $D=m_i^2-m_1^2$. Do there exist positive integers
$m_i$ and $m_1$ such that $(m_i-m_1)+(m_i-m_1)+m_i \geq
m_i^2-m_1^2$?
\end{tabular}
\\ \\
If $m_1=m_i$, this inequality holds. Then $R=\hat{r_{23}}$ and thus,
if there is a label $1$ under $Q_i$, one has that $R=0$ and
\underline{$\chi({E_i^{\circ}})=T=0$}. If there is no label $1$
under $Q_i$, then $R=r_{23}=m_i$, so
\underline{$\chi({E_i^{\circ}})=T-m_i$}. If $Q_i$ is the origin,
then we have \underline{$\chi({E_i^{\circ}})=2-m_i$}. If only label
$2$ or only label $3$ appears under $Q_i$, then
\underline{$\chi({E_i^{\circ}})=1-m_i$}. If label $2$ as well as
label $3$ are present under $Q_i$, then
\underline{$\chi({E_i^{\circ}})=-m_i$}.
\\ \\
Suppose now that $m_1 < m_i$ and that the inequality holds. This
implies that
\[(m_i-1)(m_i-3) \geq m_1(m_1-2) \geq m_i(m_i-3).\]
Then $(m_i,m_1)=(3,2)$ or $(m_i,m_1)=(2,1)$. If $(m_i,m_1)=(3,2)$,
then $\chi({E_i^{\circ}})=5-R+T$ and $R \leq 5$. If $R=5$, then
$Q_i$ is the origin. Then $T=2$ and thus $\chi({E_i^{\circ}}) > 0$.
If $(m_i,m_1)=(2,1)$, then $\chi({E_i^{\circ}})=3-R+T$ and $R \leq
4$. If $R \geq 3$, then we should have $\hat{r_{23}}=r_{23}$ and
also say $\hat{r_{12}}=r_{12}$. Thus label $1$ and label $3$ do not
appear under $Q_i$. Then only label $2$ appears below $Q_i$ or $Q_i$
is the origin. However, also under these conditions we have
$\chi({E_i^{\circ}})
> 0$.
\\ \\
\begin{tabular}{p{4cm}p{8.5cm}}
\begin{pspicture}(0,-1)(4.9,1)
\rput(0.2,0.8){\textbf{I2$^*$}}
\psline{-,linestyle=dashed}(0,0)(0.5,0) \psdot(0.5,0)
\rput(0.5,0.3){\footnotesize{$Q_i$}} \psline{-}(0.5,0)(1.5,0)
\psdot(1.5,0) \psline{-}(1.5,0)(2.5,0) \psdot(2.5,0)
\rput(1,0.2){\scriptsize{\emph{1}}}\rput(2,0.2){\scriptsize{\emph{2}}}\psline{-,linestyle=dashed}(2.5,0)(3,0)
\rput(0.5,-0.3){\footnotesize{$m_i$}}
\rput(1.5,-0.3){\footnotesize{$m_1$}}
\end{pspicture}
& \vspace*{-1.6cm} $L = m_i^2-m_i m_1$. Do there exist
multiplicities for which $(m_i-\sum_{t \geq
0}m_{Q_i(1,2^t)})+(m_i-m_1)+m_i \geq m_i^2-m_i m_1$?
\end{tabular}
\\ \\
We rewrite this inequality as $m_1(m_i-2)-\sum_{t \geq
1}m_{Q_i(1,2^t)} \geq m_i(m_i-3)$.
\begin{itemize}
\item If $m_1=m_i-1$, the cluster becomes \\
\begin{pspicture}(-3,-1)(4.9,1)
\psline{-,linestyle=dashed}(0,0)(0.5,0) \psdot(0.5,0)
\rput(0.5,0.3){\footnotesize{$Q_i$}} \psline{-}(0.5,0)(1.5,0)
\psdot(1.5,0) \psline{-}(1.5,0)(2.5,0) \psdot(2.5,0)
\rput(1,0.2){\scriptsize{\emph{1}}}\rput(2,0.2){\scriptsize{\emph{2}}}
\rput(0.5,-0.3){\footnotesize{$m_i$}}
\rput(1.5,-0.3){\footnotesize{$m_i-1$}}
\rput(2.5,-0.3){\footnotesize{$1$}}
\psline{-,linestyle=dashed}(2.5,0)(6.5,0)
\rput(4.5,0.2){\scriptsize{\emph{3}}}
 \psdot(6.5,0)
\rput(6.5,-0.3){\footnotesize{$1$}}
\end{pspicture}\\
where label $3$ appears say $k$ times, with $0 \leq k \leq m_i-2$.
Then $\chi({E_i^{\circ}}) = m_i^2-(m_i-1)^2-(k+1)-R+T$ and $R \leq
m_i+1$, so $\chi({E_i^{\circ}}) \geq 2m_i-2-k-m_i-1+T =m_i-3-k+T$.
If $k=m_i-3$, then $\chi({E_i^{\circ}})$ could only be $0$ if
$R=m_i+1$ and $T=0$ but this is impossible. If $k=m_i-2$ and
$\chi({E_i^{\circ}})\leq 0$, then $R$ should be $m_i$ or $m_i+1$. If
$R=m_i$, then $\chi({E_i^{\circ}}) = T$. We find that
\underline{$\chi({E_i^{\circ}}) = 0$} if label $2$ and $3$ appear
below $Q_i$. When $R=m_i+1$, then $\chi({E_i^{\circ}}) = -1+T$. Then
\underline{$\chi({E_i^{\circ}}) = 0$} if we only have label $3$
under $Q_i$.
\item If $m_i \geq 3$ and if the inequality holds, then certainly $m_1 > m_i-3$.
So suppose now that $m_1=m_i-2$. Then the inequality becomes
$(m_i-2)(m_i-2)-\sum_{t \geq 1}m_{Q_i(1,2^t)} \geq m_i(m_i-3)$ or
$4-\sum_{t \geq 1}m_{Q_i(1,2^t)} \geq m_i$ and so $m_i=3$. The
cluster is then of the form \\
\begin{tabular}{p{5cm}p{1cm}p{5cm}}\begin{pspicture}(-2,-1)(3,1)
\psline{-,linestyle=dashed}(0,0)(0.5,0) \psdot(0.5,0)
\rput(0.5,0.3){\footnotesize{$Q_i$}} \psline{-}(0.5,0)(1.5,0)
\psdot(1.5,0) \psline{-}(1.5,0)(2.5,0) \psdot(2.5,0)
\rput(1,0.2){\scriptsize{\emph{1}}}\rput(2,0.2){\scriptsize{\emph{2}}}
\rput(0.5,-0.3){\footnotesize{$3$}}
\rput(1.5,-0.3){\footnotesize{$1$}}
\rput(2.5,-0.3){\footnotesize{$1$}}
\end{pspicture} &  \vspace*{-1.1cm}or & \begin{pspicture}(0,-1)(4.9,1)
\psline{-,linestyle=dashed}(0,0)(0.5,0) \psdot(0.5,0)
\rput(0.5,0.3){\footnotesize{$Q_i$}} \psline{-}(0.5,0)(1.5,0)
\psdot(1.5,0) \psline{-}(1.5,0)(2.5,0) \psdot(2.5,0)
\rput(1,0.2){\scriptsize{\emph{1}}}\rput(2,0.2){\scriptsize{\emph{2}}}
\rput(0.5,-0.3){\footnotesize{$3$}}
\rput(1.5,-0.3){\footnotesize{$1$}}
\rput(2.5,-0.3){\footnotesize{$1$}} \psline{-}(2.5,0)(3.5,0)
\psdot(3.5,0)\rput(3,0.2){\scriptsize{\emph{2}}}\rput(3.5,-0.3){\footnotesize{$1$}}
\end{pspicture} \end{tabular} \\
In the first picture $\chi({E_i^{\circ}}) = 7-R+T$ and $R \leq 6$,
and thus $\chi({E_i^{\circ}}) > 0$. In the picture at the right,
$\chi({E_i^{\circ}}) = 6-R+T$ and $R \leq 5$, and thus again
$\chi({E_i^{\circ}}) > 0$.
\end{itemize} ${}$ \\
\begin{tabular}{p{3.5cm}p{9cm}}
\begin{pspicture}(0,-1)(4.9,1)
\rput(0.2,0.8){\textbf{I3$^*$}}
\psline{-,linestyle=dashed}(0,0)(0.5,0) \psdot(0.5,0)
\rput(0.5,0.3){\footnotesize{$Q_i$}} \psline{-}(0.5,0)(1.5,0)
\psdot(1.5,0) \psline{-}(1.5,0)(2.5,0.5) \psdot(2.5,0.5)
\psline{-}(1.5,0)(2.5,-0.5) \psdot(2.5,-0.5)
\rput(1,0.2){\scriptsize{\emph{1}}}\rput(2,0.45){\scriptsize{\emph{2}}}\rput(2,-0.5){\scriptsize{\emph{3}}}
\psline{-,linestyle=dashed}(2.5,0.5)(3,0.5)
\psline{-,linestyle=dashed}(2.5,-0.5)(3,-0.5)
\rput(0.5,-0.3){\footnotesize{$m_i$}}
\rput(1.5,-0.3){\footnotesize{$m_1$}}
\end{pspicture}
& \vspace*{-1.6cm} $L = m_i^2-m_im_1$. Do there exist multiplicities
such that $(m_i-\sum_{t \geq 0}m_{Q_i(1,2^t)})+(m_i-\sum_{t \geq
0}m_{Q_i(1,3^t)})+ m_i \geq m_i^2-m_im_1$?
\end{tabular}
\\ \\
We rewrite the inequality as follows:
\[m_1(m_i-2) - \sum_{t \geq 1}m_{Q_i(1,2^t)} - \sum_{t \geq 1}m_{Q_i(1,3^t)} \geq m_i(m_i-3).\]
If this inequality holds, then certainly $m_1=m_i-1$. Let $k \in
\{2,\cdots,m_i-1\}$ be the number of points that are proximate to
$Q_i$ and that are different from $Q_i(1)$. Then we find that
$\chi({E_i^{\circ}}) = m_i^2-(m_i-1)^2-k-R+T$. As $R \leq m_i$, we
have $\chi({E_i^{\circ}}) \geq m_i-1-k+T$. It follows that
\underline{$\chi({E_i^{\circ}}) = 0$} if $k=m_i-1$, $R=m_i$ and when
label $2$ and $3$ appear below $Q_i$.
\\ \\\begin{tabular}{p{4cm}p{8.5cm}}
\begin{pspicture}(0,-1)(3,1)
\rput(0.2,0.8){\textbf{II7$^*$}}\psline{-,linestyle=dashed}(0,0)(0.5,0)
\psdot(0.5,0) \rput(0.45,0.3){\footnotesize{$Q_i$}}
\psline{-}(0.5,0)(1.5,0.5) \psline{-}(0.5,0)(1.5,-0.5)
\psdot(1.5,0.5)
\psdot(1.5,-0.5)\rput(1,0.45){\scriptsize{\emph{1}}}\rput(1,-0.5){\scriptsize{\emph{2}}}
\psline{-}(1.5,0.5)(2.5,0.5)\psdot(2.5,0.5)\psline{-,linestyle=dashed}(2.5,0.5)(3,0.5)
\rput(2,0.7){\scriptsize{\emph{3}}}
\psline{-}(1.5,-0.5)(2.5,-0.5)\psdot(2.5,-0.5)\psline{-,linestyle=dashed}(2.5,-0.5)(3,-0.5)
\rput(2,-0.7){\scriptsize{\emph{3}}}
\rput(0.5,-0.3){\footnotesize{$m_i$}}\rput(1.5,0.7){\footnotesize{$m_1$}}\rput(1.5,-0.8){\footnotesize{$m_2$}}\rput(2.5,0.7)
{\footnotesize{$m_3$}}\rput(2.5,-0.8){\footnotesize{$m_4$}}
\end{pspicture}
& \vspace*{-1.7cm} $L =m_i^2-m_1 \sum_{t \geq 0}m_{Q_i(1,3^t)} - m_2
\sum_{t \geq 0}m_{Q_i(2,3^t)}$. We allow that $m_3$ and $m_4$ are
$0$, thus we include the constellations II1 and II3.
\end{tabular}
\begin{itemize}
\item Suppose $r_{12}=0$. Can the following inequality hold?
\footnotesize{\begin{eqnarray*} \left(m_i-\sum_{t \geq
0}m_{Q_i(1,3^t)}\right) + \left(m_i-\sum_{t \geq
0}m_{Q_i(2,3^t)}\right) & \geq & m_i^2-m_1
\sum_{t \geq 0}m_{Q_i(1,3^t)}\\& &  - m_2 \sum_{t \geq 0}m_{Q_i(2,3^t)} \\
\Updownarrow
\\(m_1-1)\sum_{t \geq 0}m_{Q_i(1,3^t)} + (m_2-1)\sum_{t \geq
0}m_{Q_i(2,3^t)} & \geq & m_i(m_i-2).
\end{eqnarray*}}
\normalsize On the other hand we have \footnotesize{
\begin{eqnarray*}
m_i(m_i-2) & =  & m_i(m_1-1+m_2-1) \\
& \geq &  (m_1-1)\sum_{t \geq 0}m_{Q_i(1,3^t)} + (m_2-1)\sum_{t \geq
0}m_{Q_i(2,3^t)}, \end{eqnarray*}} \normalsize and thus $m_i(m_i-2)=
(m_1-1)\sum_{t \geq 0}m_{Q_i(1,3^t)} + (m_2-1)\sum_{t \geq
0}m_{Q_i(2,3^t)}$. The cluster has then one of the following forms:
\begin{itemize}
\item $m_1=m_2=1$: if the cluster is \\ \begin{pspicture}(-4,-1)(4.9,1)
\psline{-,linestyle=dashed}(0,0)(0.5,0) \psdot(0.5,0)
\rput(0.45,0.3){\footnotesize{$Q_i$}} \psline{-}(0.5,0)(1.5,0.5)
\psline{-}(0.5,0)(1.5,-0.5) \psdot(1.5,0.5)
\psdot(1.5,-0.5)\rput(1,0.45){\scriptsize{\emph{1}}}\rput(1,-0.5){\scriptsize{\emph{2}}}
\rput(0.5,-0.3){\footnotesize{$2$}}\rput(1.5,0.7){\footnotesize{$1$}}\rput(1.5,-0.8){\footnotesize{$1$}}
\end{pspicture}
\\ then $\chi({E_i^{\circ}}) = 2-R+T$ with $R \leq 2$. However, if
$R=2$, then $T > 0$, hence $\chi({E_i^{\circ}}) > 0$. \\ The other
clusters for which $m_1=m_2=1$ will be treated in the next cases.
\item $m_1=1$ and $\sum_{t \geq 0}m_{Q_i(2,3^t)}=m_i$: \\ \begin{pspicture}(-3,-1)(4.9,1)
\psline{-,linestyle=dashed}(0,0)(0.5,0) \psdot(0.5,0)
\rput(0.45,0.3){\footnotesize{$Q_i$}} \psline{-}(0.5,0)(1.5,0.5)
\psline{-}(0.5,0)(1.5,-0.5) \psdot(1.5,0.5)
\psdot(1.5,-0.5)\rput(1,0.45){\scriptsize{\emph{1}}}\rput(1,-0.5){\scriptsize{\emph{2}}}
\psline{-,linestyle=dashed}(1.5,0.5)(2.5,0.5)
\psline{-}(1.5,-0.5)(2.5,-0.5)\psdot(2.5,-0.5)
\rput(2.2,-0.7){\scriptsize{\emph{3}}}\rput(0.5,-0.3){\footnotesize{$m_i$}}\rput(1.5,0.7){\footnotesize{$1$}}\rput(1.5,-0.8){\footnotesize{$m_i-1$}}
\rput(2.5,-0.8){\footnotesize{$1$}}
\psline{-,linestyle=dashed}(2.5,0.5)(4.5,0.5) \psdot(4.5,0.5)
\psline{-,linestyle=dashed}(2.5,-0.5)(5,-0.5) \psdot(5,-0.5)
\rput(4.5,0.7){\footnotesize{$1$}}\rput(5,-0.8){\footnotesize{$1$}}\rput(4,-0.7){\scriptsize{\emph{1}}}\rput(3,0.7){\scriptsize{\emph{3}}}
\end{pspicture}
\\ Suppose that the multiplicity $1$ appears $k \in \{1,\cdots,m_i\}$ times in the upper
chain and that the label $1$ appears $l-1$ times in the lower chain,
$1 \leq l \leq m_i-1$. We have
$\chi({E_i^{\circ}})=m_i^2-(m_i-1)^2-l-k-R+T=2m_i-1-l-k-R+T$. As $R
\leq m_i-k$, we get $\chi({E_i^{\circ}}) \geq m_i-1-l+T$. If
$\chi({E_i^{\circ}})\leq 0$, then we must have that $R=m_i-k$,
$l=m_i-1$ and $T=0$. \\ \indent If $k < m_i$, then label $2$ may not
appear under $Q_i$ (indeed, $R=\check{r_{13}}$) and label $1$ should
certainly appear under $Q_i$ (see Table $1$). We then have that
\underline{$\chi({E_i^{\circ}})=0$}. If $k=m_i$, then also $\sum_{t
\geq 0}m_{Q_i(1,3^t)}=m_i$. This cluster will be treated further on.
\item $m_2=1$ and $\sum_{t \geq 0}m_{Q_i(1,3^t)}=m_i$:
up to permutation this case is the same as the previous case.
\item $\sum_{t \geq 0}m_{Q_i(1,3^t)}=m_i$ and
$\sum_{t \geq 0}m_{Q_i(2,3^t)}=m_i$: in this case $R=0$ and
therefore \underline{$\chi({E_i^{\circ}})=0$} if and only if
$D=T=0$. From Proposition \ref{lemmacomb} it follows that both
chains that leave out of $Q_i$ should be Euclidean clusters. To have
$T=0$, one needs at least two labels under $Q_i$ or exactly one
label under $Q_i$ that then should be $1$ or $2$.
\end{itemize}
\item If $r_{12} \neq 0$, we may suppose that $r_{13}=r_{23}=0$. We
study if the following inequality can hold:
\[m_i-m_1-m_2 \geq m_i^2-m_im_1-m_im_2.\]
We rewrite the inequality as $(m_1+m_2)(m_i-1) \geq m_i(m_i-1)$.
This gives a contradiction to $r_{12} \neq 0$.
\end{itemize}
%\begin{remark}
%\emph{Notice that one can always suppose that $r_{12}=0$ for the
%clusters with underlying constellation II1 and II3. The previous
%computations for the constellation II7 are still valid when $m_3$ or
%$m_4$ are equal to $0$. Hence, we can conclude that
%$\chi(E_i^{\circ})>0$ for the clusters with underlying constellation
%II1 or II3.}
%\end{remark}${}$
${}$ \\ \\ ${}$
\begin{tabular}{p{3.5cm}p{9cm}}
\begin{pspicture}(0,-1)(2.5,1)
\rput(0.2,0.8){\textbf{II11}}\psline{-,linestyle=dashed}(0,0)(0.5,0)
\psdot(0.5,0) \rput(0.45,0.3){\footnotesize{$Q_i$}}
\psline{-}(0.5,0)(1.5,0.5) \psline{-}(0.5,0)(1.5,-0.5)
\psdot(1.5,0.5)
\psdot(1.5,-0.5)\rput(1,0.45){\scriptsize{\emph{1}}}\rput(1,-0.5){\scriptsize{\emph{2}}}
\psline{-}(1.5,0.5)(2.5,0.8) \psdot(2.5,0.8)
\psline{-}(1.5,0.5)(2.5,0.2)
\psdot(2.5,0.2)\psline{-,linestyle=dashed}(2.5,0.8)(3,0.8)
\rput(2,0.85){\scriptsize{\emph{2}}}
\rput(2,0.18){\scriptsize{\emph{3}}}\psline{-,linestyle=dashed}(2.5,0.2)(3,0.2)
\psline{-}(1.5,-0.5)(2.5,-0.2)\psdot(2.5,-0.2)\psline{-,linestyle=dashed}(2.5,-0.2)(3,-0.2)
\psline{-}(1.5,-0.5)(2.5,-0.8)\psdot(2.5,-0.8)\psline{-,linestyle=dashed}(2.5,-0.8)(3,-0.8)
\rput(2,-0.8){\scriptsize{\emph{3}}}\rput(2.8,-1.1){\scriptsize{\emph{\^{3}}}}
\rput(2,-0.14){\scriptsize{\emph{1}}}
\rput(0.5,-0.3){\footnotesize{$m_i$}}\rput(1.5,0.7){\footnotesize{$m_1$}}\rput(1.5,-0.8){\footnotesize{$m_2$}}\rput(2.5,-1){\footnotesize{$m_6$}}
\rput(2.5,-0.4){\footnotesize{$m_5$}}\rput(2.5,0.4){\footnotesize{$m_4$}}\rput(2.5,1){\footnotesize{$m_3$}}
\end{pspicture}
& \vspace*{-1.6cm} We estimate in a rough way and get $L > m_i^2 -
m_im_1-m_im_2$.
\end{tabular}
\\ \\
As described in Remark \ref{remarkhulp}, let the value of $m_2$
increase as long as the cluster stays idealistic.
\begin{itemize}
\item Suppose that $r_{12}=0$. We study if the following inequality
can occur: \footnotesize{
\begin{eqnarray*}
\left(m_i-\sum_{t \geq 0}m_{Q_i(1,3^t)}\right)+(m_i-m_2-m_6) & > &
m_i^2 - m_im_1-m_im_2.
\end{eqnarray*}}
\normalsize We rewrite it as \footnotesize{\begin{eqnarray} -\sum_{t
\geq 0}m_{Q_i(1,3^t)}-m_2-m_6>m_i(m_i-2-m_1-m_2).
\end{eqnarray}}
\normalsize As $2+m_1+m_2 \leq m_i$, this inequality can never hold.
\item Suppose that $r_{12} \neq 0$ and that $r_{23}=0$. Moreover we
can suppose that $r_{13}=0$ (we let increase the value of $m_1$). We
investigate the inequality \footnotesize{
\begin{eqnarray*}
\left(m_i-m_1-\sum_{t \geq 1}m_{Q_i(1,2^t)}-m_2-\sum_{t \geq
1}m_{Q_i(2,1^t)}\right) & > &  m_i^2 - m_im_1 -m_im_2.
\end{eqnarray*}}
\normalsize We rewrite the inequality as
\footnotesize{
\begin{eqnarray} -\sum_{t \geq 0}m_{Q_i(1,2^t)}-\sum_{t \geq 0}m_{Q_i(2,1^t)} >
m_i(m_i-m_1-m_2-1)
\end{eqnarray}}
\normalsize and we see again that this can never happen.
\end{itemize}
\begin{remark}
\emph{If we allow the multiplicities for the constellation II11 to
be $0$, except for $m_i,m_1$ and $m_2$, and if we also suppose that
not both $m_3$ and $m_5$ are $0$, then we also have
$L>m_i^2-m_im_1-m_im_2$. For the clusters with underlying
constellation II2, II4, II5, II6 or II8, we may suppose that
$r_{12}=0$. It follows then from the inequality $(5)$ that
$\chi(E_i^{\circ}) > 0$ for the clusters with underlying
constellation II5 and II8. When $r_{12} \neq 0$, then it follows
from inequality $(6)$ that $\chi(E_i^{\circ}) > 0$ for the clusters
with underlying constellation II9 or II10.} \hfill $\square$
\end{remark} From this remark, it follows that we should study the
case $r_{12}=0$, $m_i < m_1+m_2+2$ for the clusters with underlying
constellation II2, II4, II6, II9 or II10. If $m_i=m_1+m_2$, then we
have a cluster whose underlying constellation is a subconstellation
of II7, and thus already treated. So suppose that $m_i=m_1+m_2+1$.
Then inequality $(5)$ becomes
\begin{eqnarray*}
1-\sum_{t \geq 1}m_{Q_i(1,3^t)}-m_6>0.
\end{eqnarray*}
It follows that $\sum_{t \geq 1}m_{Q_i(1,3^t)}=m_6 =0$ and that the
cluster is like
\\ \\ \begin{tabular}{p{4cm}p{8.5cm}}
\begin{pspicture}(0,-1)(4.9,1)
%\rput(0.2,0.8){\textbf{II6}}
\psline{-,linestyle=dashed}(0,0)(0.5,0) \psdot(0.5,0)
\rput(0.45,0.3){\footnotesize{$Q_i$}} \psline{-}(0.5,0)(1.5,0.5)
\psline{-}(0.5,0)(1.5,-0.5) \psdot(1.5,0.5)
\psdot(1.5,-0.5)\rput(1,0.45){\scriptsize{\emph{1}}}\rput(1,-0.5){\scriptsize{\emph{2}}}
\psline{-}(1.5,0.5)(2.5,0.5)\psdot(2.5,0.5)\psline{-,linestyle=dashed}(2.5,0.5)(3,0.5)
\rput(2,0.7){\scriptsize{\emph{2}}}
\rput(0.5,-0.3){\footnotesize{$m_i$}}\rput(1.5,0.7){\footnotesize{$m_1$}}\rput(1.5,-0.8){\footnotesize{$m_2$}}
\rput(2.5,0.7){\footnotesize{$1$}}
\end{pspicture}
& \vspace*{-1.7cm} We have $L = m_i^2-m_2^2-(m_i-m_2)m_1$.
\end{tabular}
Can the following inequality hold:
\begin{eqnarray*}
m_i-m_1 + m_i-m_2 \geq m_i^2-m_2^2-(m_i-m_2)m_1?
\end{eqnarray*}
Substituting $m_i$ by $m_1+m_2+1$, we get $1 \geq 2m_1m_2+m_2$. This
contradiction allows us to conclude that $\chi(E_i^{\circ})
> 0$.
${}$ \\ \\
\begin{tabular}{p{4cm}p{8.5cm}}
\begin{pspicture}(0,-1)(4.9,1)
\rput(0.2,0.8){\textbf{III9}}\psline{-,linestyle=dashed}(0,0)(0.5,0)
\psdot(0.5,0) \rput(0.45,0.3){\footnotesize{$Q_i$}}
\psline{-}(0.5,0)(1.5,0.5) \psdot(1.5,0.5)\psline{-}(0.5,0)(1.5,0)
\psdot(1.5,0) \psline{-}(0.5,0)(1.5,-0.5) \psdot(1.5,-0.5)
\rput(1,0.45){\scriptsize{\emph{1}}}\rput(1,-0.5){\scriptsize{\emph{3}}}\rput(1.1,0.15){\scriptsize{\emph{2}}}
%\psline{-}(1.5,0.5)(2.5,0.5)\psdot(2.5,0.5)\psline{-,linestyle=dashed}(2.5,0.5)(3,0.5)
%\rput(2,0.7){\emph{3}}}
\psline{-}(1.5,0)(2.5,0)\psdot(2.5,0)\psline{-,linestyle=dashed}(2.5,0)(3,0)
\rput(2,0.2){\scriptsize{\emph{1}}}\psline{-}(1.5,-0.5)(2.5,-0.5)\psdot(2.5,-0.5)\psline{-,linestyle=dashed}(2.5,-0.5)(3,-0.5)
\rput(2,-0.7){\scriptsize{\emph{1}}} \psline{-}(1.5,0.5)(2.5,1)
\psdot(2.5,1) \psline{-}(1.5,0.5)(2.5,0.5)
\psdot(2.5,0.5)\psline{-,linestyle=dashed}(2.5,1)(3,1)
\rput(2,0.95){\scriptsize{\emph{2}}}
\rput(2.2,0.7){\scriptsize{\emph{3}}}\psline{-,linestyle=dashed}(2.5,0.5)(3,0.5)
\rput(0.5,-0.3){\footnotesize{$m_i$}}\rput(1.5,0.7){\footnotesize{$m_1$}}\rput(1.5,0.2){\footnotesize{$m_2$}}
\rput(1.5,-0.8){\footnotesize{$m_3$}}\rput(2.6,1.2){\footnotesize{$m_4$}}\rput(2.6,0.7){\footnotesize{$m_5$}}
\rput(2.6,-0.8){\footnotesize{$m_7$}}\rput(2.6,0.2){\footnotesize{$m_6$}}
\end{pspicture}
& \vspace*{-1.7cm} A rough estimate gives $L > m_i^2-m_im_1-m_2(m_i-
\sum_{t \geq 0}m_{Q_i(1,2^t)}-r_{12})-m_3(m_i-\sum_{t \geq
0}m_{Q_i(1,3^t)}-r_{13})$.
\end{tabular}
\begin{itemize}
\item \normalsize{We let increase the value of $m_1$; suppose that we then get $r_{13}=0$. We now also let $m_2$ increase; suppose that $r_{12}$ becomes $0$. Can
the following inequality then hold:} \footnotesize{
\begin{eqnarray*}
m_i-m_2-m_3 \hspace{-0.35cm} &  > & \hspace{-0.3cm}
m_i^2-m_im_1-m_2\left(m_i- \sum_{t \geq
0}m_{Q_i(1,2^t)}\right)-m_3\left(m_i - \sum_{t \geq
0}m_{Q_i(1,3^t)}\right)?\end{eqnarray*}} \normalsize
 We rewrite the inequality as follows:
\footnotesize{
\begin{eqnarray}
0 > (m_i-m_1-1)(m_i-m_2-m_3)+ m_2 \sum_{t \geq
1}m_{Q_i(1,2^t)}+m_3\sum_{t \geq 1}m_{Q_i(1,3^t)}.
\end{eqnarray}}
\normalsize This inequality can never be true.
\item We let increase $m_1$; suppose that we get $r_{13}=0$. Then we let increase the value of $m_2$ and
$r_{23}$ becomes $0$: can \footnotesize{
\begin{eqnarray*}
r_{12} > m_i^2 -m_im_1-m_2(m_i- \sum_{t \geq
0}m_{Q_i(1,2^t)}-r_{12})-m_3(m_i-\sum_{t \geq 0}m_{Q_i(1,3^t)})?
\end{eqnarray*}}
\normalsize
 As $m_i=m_2+m_3$, we get
\begin{eqnarray}
r_{12} > r_{12}m_2 + m_2\sum_{t \geq 1}m_{Q_i(1,2^t)}+m_3\sum_{t
\geq 1}m_{Q_i(1,3^t)},
\end{eqnarray}
which is never satisfied. %\footnotesize{
%\begin{eqnarray*}
%m_i- \sum_{t=0}^{\infty}m_{Q_i(1,2^t)} -
%\sum_{t=0}^{\infty}m_{Q_i(2,1^t)} & > &  m_i^2-m_im_1 -m_2\left(m_i-
%\sum_{t=0}^{\infty}m_{Q_i(1,2^t)}\right)- \\
%& & m_3\left(m_i-\sum_{t=0}^{\infty}m_{Q_i(1,3^t)}\right)?
%\end{eqnarray*}}
%\normalsize
% As $m_i=m_2+m_3$, we get
%\begin{eqnarray}
%m_3>\sum_{t=1}^{\infty}m_{Q_i(2,1^t)} +
%(m_2+1)\sum_{t=1}^{\infty}m_{Q_i(1,2^t)}+m_1+m_3\sum_{t=1}^{\infty}m_{Q_i(1,3^t)},
%\end{eqnarray}
%which is never satisfied.
\item We let increase $m_1$; suppose that $r_{12}$ becomes $0$. Now we let increase $m_3$
and suppose $r_{23}$ becomes $0$ (we already treated
$r_{12}=r_{13}=0$). Can \footnotesize{
\begin{eqnarray*}
r_{13} > m_i^2-m_im_1-m_2(m_i- \sum_{t \geq
0}m_{Q_i(1,2^t)})-m_3(m_i-\sum_{t \geq 0}m_{Q_i(1,3^t)}-r_{13})?
\end{eqnarray*}}
\normalsize As $m_i=m_2+m_3$, we get
\begin{eqnarray}
r_{13} > r_{13}m_3 +m_2\sum_{t \geq 1}m_{Q_i(1,2^t)} +m_3\sum_{t
\geq 1}m_{Q_i(1,3^t)},
\end{eqnarray}
which can not hold.
%\footnotesize{
%\begin{eqnarray*}
%m_i- \sum_{t=0}^{\infty}m_{Q_i(1,3^t)} -
%\sum_{t=0}^{\infty}m_{Q_i(3,1^t)} & > &  m_i^2-m_im_1 -m_2\left(m_i-
%\sum_{t=0}^{\infty}m_{Q_i(1,2^t)}\right)- \\
%& & m_3\left(m_i-\sum_{t=0}^{\infty}m_{Q_i(1,3^t)}\right)?
%\end{eqnarray*}}
%\normalsize As $m_i=m_2+m_3$, we get
%\begin{eqnarray}
%m_2>\sum_{t=1}^{\infty}m_{Q_i(3,1^t)} +
%(m_3+1)\sum_{t=1}^{\infty}m_{Q_i(1,3^t)}+m_1+m_2\sum_{t=1}^{\infty}m_{Q_i(1,2^t)},
%\end{eqnarray} which can not hold.
\end{itemize}
\begin{remark}
\emph{Notice that one can use the same lower bound for $L$ for the
subconstellations IIIx with $1 \leq x \leq 8$ of constellation III9
and that the inequalities $(7), (8)$ and $(9)$ neither hold for
them.}
%For the subclusters with
%$m_5 \neq 0$ (keep the permutation into account!), we also have that
%equation $(8)$ does not hold. In particular, $(8)$ never holds for
%III3, III5, III7 and III8. Analogously we find that $(9)$ does never
%hold for the constellations III2, III3, III4 and III7.
\hfill $\square$ \end{remark}
This closes the computational part that yields the classification of
the $\chi(E_i^{\circ})$. In particular, we get the following
results.
\begin{theorem}
Let $f$ be a polynomial map that is general with respect to a
$3$-dimensional toric idealistic cluster
$\mathcal{A}=(\mathcal{C},\underline{m})$. If $Q_i \in \mathcal{C}$,
then $\chi(E_i^{\circ})<0$ if and only if the configuration in $E_i
\cong \mathbb{P}^2$ consists of (at least three) lines - possibly
exceptional - that are all going through the same point, i.e. if and
only if $Q_i$ appears in a subcluster of List $2$ in $\mathcal{A}$.
\end{theorem}
\begin{tabular}{p{2cm}p{9cm}}
\begin{pspicture}(0,-0.5)(4.9,1)
\rput(0.2,0.8){\textbf{C1}} \psdot(0.5,0)
\rput(0.5,0.3){\footnotesize{$Q_i$}} \psline{-}(0.5,0)(1.5,0)
\psdot(1.5,0) \rput(1,0.2){\scriptsize{\emph{1}}}
\rput(0.5,-0.3){\footnotesize{$m_i$}}
\rput(1.5,-0.3){\footnotesize{$m_i$}}
\end{pspicture}
& \vspace*{-1.3cm} If $Q_i$ is the origin, then $\chi(E_i^{\circ}) =
2-m_i$. Thus, if $m_i \geq 3$, then $\chi(E_i^{\circ}) < 0$.
\end{tabular}
\\ \\
\begin{tabular}{p{2cm}p{9cm}}
\begin{pspicture}(0,-0.5)(4.9,1)
\rput(0.2,0.8){\textbf{C2}} \psline{-}(0,0)(0.5,0) \psdot(0.5,0)
\rput(0.5,0.3){\footnotesize{$Q_i$}} \psline{-}(0.5,0)(1.5,0)
\psdot(1.5,0) \rput(1,0.2){\scriptsize{\emph{1}}}
\rput(0.5,-0.3){\footnotesize{$m_i$}}
\rput(1.5,-0.3){\footnotesize{$m_i$}}
\end{pspicture}
& \vspace*{-1.3cm} If only label $2$ or only label $3$ appears under
$Q_i$, then $\chi(E_i^{\circ})  =  1-m_i$. So, if $m_i \geq 2$, then
$\chi(E_i^{\circ}) < 0$.
\end{tabular}
\\ \\
\begin{tabular}{p{2cm}p{9cm}}
\begin{pspicture}(0,-0.5)(4.9,1)
\rput(0.2,0.8){\textbf{C3}} \psline{-}(0,0)(0.5,0) \psdot(0.5,0)
\rput(0.5,0.3){\footnotesize{$Q_i$}} \psline{-}(0.5,0)(1.5,0)
\psdot(1.5,0) \rput(1,0.2){\scriptsize{\emph{1}}}
\rput(0.5,-0.3){\footnotesize{$m_i$}}
\rput(1.5,-0.3){\footnotesize{$m_i$}}
\end{pspicture}
& \vspace*{-1.3cm} If only label $2$ and label $3$ appear under
$Q_i$, then $\chi(E_i^{\circ})  =  -m_i$ and thus $\chi(E_i^{\circ})
< 0$.
\end{tabular}
\begin{center}
\emph{List $2$}
\end{center}
\begin{example}\emph{
The surface with equation $x^{2m_i}+y^{m_i}+z^{m_i}=0$ is an example
of a surface that is general with respect to the cluster
\textbf{C1}}. \hfill $\square$
\end{example}
In the general case of surfaces, there exist much more
configurations that yield a negative $\chi(E_i^{\circ})$. In
\cite{Veysconfigurations} are given such examples.
\begin{theorem}
Let $f$ be a polynomial map that is general with respect to a
$3$-dimensional toric idealistic cluster
$\mathcal{A}=(\mathcal{C},\underline{m})$. If $Q_i \in \mathcal{C}$,
then $\chi(E_i^{\circ})=0$ if and only if $Q_i$ appears in a
subcluster of List $3$ in $\mathcal{A}$.
\end{theorem}
\begin{tabular}{p{2cm}p{9cm}}
\begin{pspicture}(0,-0.5)(4.9,1)
\rput(0.2,0.8){\textbf{C4}} \psline{-}(0,0)(0.5,0) \psdot(0.5,0)
\rput(0.5,0.3){\footnotesize{$Q_i$}}
\rput(0.5,-0.3){\footnotesize{$m_i$}}
\end{pspicture}
& \vspace*{-1.3cm} If there exists exactly one or exactly two points
to which $Q_i$ is proximate and if $m_i=1$, then $\chi(E_i^{\circ})
= 0$.
\end{tabular}
\\ \\
\begin{tabular}{p{2cm}p{9cm}}
\begin{pspicture}(0,-0.5)(4.9,1)
\rput(0.2,0.8){\textbf{C5}} \psdot(0.5,0)
\rput(0.5,0.3){\footnotesize{$Q_i$}} \psline{-}(0.5,0)(1.5,0)
\psdot(1.5,0) \rput(1,0.2){\scriptsize{\emph{1}}}
\rput(0.5,-0.3){\footnotesize{$m_i$}}
\rput(1.5,-0.3){\footnotesize{$m_i$}}
\end{pspicture}
& \vspace*{-1.3cm} If $Q_i$ is the origin and if $m_i=2$, then
$\chi(E_i^{\circ}) = 0$.
\end{tabular}
\\ \\
\begin{tabular}{p{2cm}p{9cm}}
\begin{pspicture}(0,-0.5)(4.9,1)
\rput(0.2,0.8){\textbf{C6}} \psline{-}(0,0)(0.5,0) \psdot(0.5,0)
\rput(0.5,0.3){\footnotesize{$Q_i$}} \psline{-}(0.5,0)(1.5,0)
\psdot(1.5,0) \rput(1,0.2){\scriptsize{\emph{1}}}
\rput(0.5,-0.3){\footnotesize{$m_i$}}
\rput(1.5,-0.3){\footnotesize{$m_i$}}
\end{pspicture}
& \vspace*{-1.3cm} If only label $2$ or only label $3$ appears under
$Q_i$ and if $m_i=1$, then $\chi(E_i^{\circ})  =  0$.
\end{tabular}
\\ \\
\begin{tabular}{p{2cm}p{9cm}}
\begin{pspicture}(0,-0.5)(4.9,1)
\rput(0.2,0.8){\textbf{C7}} \psline{-}(0,0)(0.5,0) \psdot(0.5,0)
\rput(0.5,0.3){\footnotesize{$Q_i$}} \psline{-}(0.5,0)(1.5,0)
\psdot(1.5,0) \rput(1,0.2){\scriptsize{\emph{1}}}
\rput(0.5,-0.3){\footnotesize{$m_i$}}
\rput(1.5,-0.3){\footnotesize{$m_i$}}
\end{pspicture}
& \vspace*{-1cm} If at least label $1$ appears under $Q_i$, then
$\chi(E_i^{\circ})  =  0$.
\end{tabular}
\\ \\
\begin{tabular}{p{4cm}p{8cm}}
\begin{pspicture}(0,-1)(3,1.5)
\rput(0.2,0.8){\textbf{C8}}\psline{-}(0,0)(0.5,0) \psdot(0.5,0)
\rput(0.45,0.3){\footnotesize{$Q_i$}} \psline{-}(0.5,0)(1.5,0.5)
\psline{-}(0.5,0)(1.5,-0.5) \psdot(1.5,0.5)
\psdot(1.5,-0.5)\rput(1,0.45){\scriptsize{\emph{1}}}\rput(1,-0.5){\scriptsize{\emph{2}}}
\psline{-}(1.5,0.5)(2.5,0.5)\psdot(2.5,0.5)\psline{-,linestyle=dashed}(2.5,0.5)(3,0.5)
\rput(2,0.7){\scriptsize{\emph{3}}}
\psline{-}(1.5,-0.5)(2.5,-0.5)\psdot(2.5,-0.5)\psline{-,linestyle=dashed}(2.5,-0.5)(3,-0.5)
\rput(2,-0.7){\scriptsize{\emph{3}}}
\rput(0.5,-0.3){\footnotesize{$m_i$}}\rput(1.5,0.7){\footnotesize{$m_1$}}\rput(1.5,-0.8){\footnotesize{$m_2$}}
\end{pspicture}
& \vspace*{-1.8cm}  If $m_1+m_2=m_i$, if the upper chain and the
lower chain are Euclidean clusters and \newline \textbf{A}. if only
label $1$ or only label $2$ appears under $Q_i$, then
$\chi(E_i^{\circ}) = 0$; or
\newline \textbf{B}. if at least two different labels appear
under $Q_i$, then $\chi(E_i^{\circ}) = 0$.
\end{tabular}\\
\begin{tabular}{p{6.2cm}p{5.8cm}}
\begin{pspicture}(0,-1)(7,1)
\rput(0.2,0.8){\textbf{C9}} \psline{-}(0,0)(0.5,0) \psdot(0.5,0)
\rput(0.5,0.3){\footnotesize{$Q_i$}} \psline{-}(0.5,0)(1.5,0)
\psdot(1.5,0) \psline{-}(1.5,0)(2.5,0) \psdot(2.5,0)
\rput(1,0.2){\scriptsize{\emph{1}}}\rput(2,0.2){\scriptsize{\emph{2}}}
\rput(0.5,-0.3){\footnotesize{$m_i$}}
\rput(1.5,-0.3){\footnotesize{$m_i-1$}}
\rput(2.5,-0.3){\footnotesize{$1$}}
\rput(3.7,0.2){\scriptsize{\emph{3}}}
\psline{-,linestyle=dashed}(2.5,0)(5,0) \psdot(5,0)
\rput(5,-0.3){\scriptsize{1}}
\end{pspicture}
& \vspace*{-1.6cm}  If label $3$ appears $m_i-2$ times and \newline
\textbf{A}. if only label $3$ appears under $Q_i$, then
$\chi(E_i^{\circ}) = 0$; or
\newline \textbf{B}. if only label $2$
and $3$ appear under $Q_i$, then $\chi(E_i^{\circ}) = 0$.
\end{tabular}
\\
\begin{tabular}{p{3.2cm}p{8.8cm}}
\begin{pspicture}(0,-1)(3,1)
\rput(0.2,0.8){\textbf{C10}} \psline{-}(0,0)(0.5,0) \psdot(0.5,0)
\rput(0.5,0.3){\footnotesize{$Q_i$}}
\rput(1.5,0.3){\footnotesize{$P$}}\psline{-}(0.5,0)(1.5,0)
\psdot(1.5,0) \psline{-}(1.5,0)(2.5,0.5) \psdot(2.5,0.5)
\psline{-}(1.5,0)(2.5,-0.5) \psdot(2.5,-0.5)
\rput(1,0.2){\scriptsize{\emph{1}}}\rput(2,0.45){\scriptsize{\emph{2}}}\rput(2,-0.5){\scriptsize{\emph{3}}}
\psline{-,linestyle=dashed}(2.5,0.5)(3,0.5)
\psline{-,linestyle=dashed}(2.5,-0.5)(3,-0.5)
\rput(0.5,-0.3){\footnotesize{$m_i$}}
\rput(1.4,-0.3){\footnotesize{$m_i-1$}}
\end{pspicture}
& \vspace*{-1.6cm} If $\# \{s \mbox{ $|$ } s \in \mathbb{Z}_{\geq
0}, P(2,3^s) \in \mathcal{C}\} + \# \{s \mbox{ $|$ } s \in
\mathbb{Z}_{\geq 0}, P(3,2^s) \in \mathcal{C}\}=m_i-1$ and if only
label $2$ and label $3$ appear under $Q_i$, then $\chi(E_i^{\circ})
= 0$.
\end{tabular}
\\ \\ \\ \\\begin{tabular}{p{5cm}p{7cm}}
\begin{pspicture}(0,0)(4.9,1)
\rput(0.2,0.8){\textbf{C11}} \psline{-}(0,0)(0.5,0) \psdot(0.5,0)
\rput(0.45,0.3){\footnotesize{$Q_i$}}
\rput(1.45,-0.3){\footnotesize{$P$}} \psline{-}(0.5,0)(1.5,0.5)
\psline{-}(0.5,0)(1.5,-0.5) \psdot(1.5,0.5)
\psdot(1.5,-0.5)\rput(1,0.45){\scriptsize{\emph{1}}}\rput(1,-0.5){\scriptsize{\emph{2}}}
\psline{-}(1.5,0.5)(2.5,0.5)\psdot(2.5,0.5)
\rput(2,0.7){\scriptsize{\emph{3}}}
\psline{-}(1.5,-0.5)(2.5,-0.5)\psdot(2.5,-0.5)
\rput(2.2,-0.7){\scriptsize{\emph{3}}}\rput(0.5,-0.3){\footnotesize{$m_i$}}\rput(1.5,0.7){\footnotesize{$1$}}\rput(1.5,-0.8){\footnotesize{$m_i-1$}}
\rput(2.5,0.7){\footnotesize{$1$}}\rput(2.5,-0.8){\footnotesize{$1$}}
\psline{-,linestyle=dashed}(2.5,0.5)(4.5,0.5) \psdot(4.5,0.5)
\psline{-}(2.5,-0.5)(3.5,-0.5)
\psline{-,linestyle=dashed}(3.5,-0.5)(5,-0.5) \psdot(5,-0.5)
\rput(4.5,0.7){\footnotesize{$1$}}\rput(5,-0.8){\footnotesize{$1$}}\rput(4.3,-0.7){\scriptsize{\emph{1}}}\rput(3.5,0.7){\scriptsize{\emph{3}}}
\psdot(3.5,-0.5)\rput(3,-0.7){\scriptsize{\emph{1}}}\rput(3.5,-0.8){\footnotesize{$1$}}
\end{pspicture}
& \vspace*{-1.5cm} If $\# \{s \mbox{ $|$ } s \in \mathbb{Z}_{\geq
0}, P(3,1^s) \in \mathcal{C}\} =m_i-1$, if $\# \{s \mbox{ $|$ } s
\in \mathbb{Z}_{\geq 0}, Q_i(1,3^s) \in \mathcal{C}\} =k$ for $1
\leq k \leq m_i$ and \newline \textbf{A}. if only label $1$ appears
under $Q_i$, then $\chi(E_i^{\circ}) = 0$; or
\newline \textbf{B}. if only label $1$ and label $3$ appear under $Q_i$, then
$\chi(E_i^{\circ}) = 0$.
\end{tabular}
\begin{center}
\emph{List $3$}
\end{center}
\begin{example}\emph{
The ideal
\[I=(x^9,y^5,z^5,x^6y,x^5y^2,x^3y^3,x^2y^4,y^4z,y^3z^2,y^2z^3,yz^4,xz^4,\]\[x^2z^3,x^5z^2,x^7z,xyz^3,xy^2z^2,xy^3z,x^3yz^2,x^3y^2z,x^5yz)\]
is the complete finitely supported ideal that corresponds to the
cluster
\\
\begin{pspicture}(-5,-1)(3,1)
%\psline{-,linestyle=dashed}(-1,0)(-0.5,0)
\psline{-}(-0.5,0)(0.5,0)\psdot(-0.5,0) \psdot(0.5,0)
%\rput(0.45,0.3){\footnotesize{$Q_i$}}
\psline{-}(0.5,0)(1.5,0.5) \psline{-}(0.5,0)(1.5,-0.5)
\psdot(1.5,0.5)
\psdot(1.5,-0.5)\rput(1,0.45){\scriptsize{\emph{1}}}\rput(1,-0.5){\scriptsize{\emph{2}}}
\psline{-}(1.5,0.5)(2.5,0.5)\psdot(2.5,0.5)\psline{-}(2.5,0.5)(3.5,0.5)\psdot(3.5,0.5)
\rput(2,0.7){\scriptsize{\emph{3}}}\rput(3,0.7){\scriptsize{\emph{3}}}
\psline{-}(1.5,-0.5)(2.5,-0.5)\psdot(2.5,-0.5)\psline{-}(2.5,-0.5)(3.5,-0.5)\psdot(3.5,-0.5)
\rput(2,-0.7){\scriptsize{\emph{3}}}\rput(3,-0.7){\scriptsize{\emph{1}}}
\rput(-0.5,-0.3){\footnotesize{$5$}}\rput(0.5,-0.3){\footnotesize{$3$}}\rput(1.5,0.7){\footnotesize{$1$}}\rput(2.5,0.7){\footnotesize{$1$}}\rput(3.5,0.7){\footnotesize{$1$}}
\rput(1.5,-0.8){\footnotesize{$2$}}\rput(2.5,-0.8){\footnotesize{$1$}}\rput(3.5,-0.8){\footnotesize{$1$}}\rput(0,0.2){\scriptsize{\emph{1}}}
%\rput(1.5,-0.15){\footnotesize{$Q_j$}}
\end{pspicture} \\ A general element of $I$ illustrates a surface with a
singularity as in cluster \textbf{C8}. \\ \\ Let $J$ be the ideal
\[(x^6,y^6,z^9,x^5y,x^4y^2,x^3y^3,x^2y^4,xy^5,y^5z,y^4z^2,y^3z^3,y^2z^5,yz^7,x^5z,x^4z^3,x^3z^4,\]
\[x^2z^6,xz^7,xyz^6,xy^2z^4,xy^3z^2,xy^4z,x^2yz^4,x^2y^2z^2,x^2y^3z,x^3yz^2,x^3y^2z,x^4yz).\]
This is the complete finitely supported ideal corresponding with the
cluster
\\
\begin{pspicture}(-5,-1)(4,1)
\psdot(-0.5,0) \psline{-}(-0.5,0)(0.5,0) \psdot(0.5,0)
% \rput(0.5,0.3){\footnotesize{$Q_i$}}
\psline{-}(0.5,0)(1.5,0) \psdot(1.5,0) \psline{-}(1.5,0)(2.5,0)
\psdot(2.5,0)
\rput(1,0.2){\scriptsize{\emph{1}}}\rput(2,0.2){\scriptsize{\emph{2}}}
\rput(-0.5,-0.3){\footnotesize{$6$}}
\rput(0.5,-0.3){\footnotesize{$3$}}
\rput(1.5,-0.3){\footnotesize{$2$}}
\rput(2.5,-0.3){\footnotesize{$1$}}
\rput(3,0.2){\scriptsize{\emph{3}}} \psline{-}(2.5,0)(3.5,0)
\rput(0,0.2){\scriptsize{\emph{3}}} \psdot(3.5,0)
\rput(3.5,-0.3){\footnotesize{$1$}}
\end{pspicture}
\\ A general element of $J$ illustrates a surface with a singularity as in cluster \textbf{C9}.
}\end{example}\vspace*{-0.3cm} \hfill $\square$
\begin{remark}
\emph{Let $Q_l$ be a point with multiplicity $1$ in a
$3$-dimensional toric idealistic constellation $\mathcal{C}$ and let
$Q_k \in \mathcal{C}$ be such that $Q_l \succ Q_k$. Suppose that
$Q_l$ is lying only on the irreducible exceptional component $E_{k}
\cong \mathbb{P}^2$. Then obviously $C_k$ has normal crossings in
$Q_l$. Suppose that $Q_l$ is lying on exactly two exceptional
components $E_k \cong \mathbb{P}^2$ and $E_j$. If $C_k$ does not
have normal crossings in $Q_l$ then $E_k \cap E_j$ should be the
tangent line to $C_k$ in $Q_l$. After blowing up in the point $Q_l$,
one needs at least one more blowing up to obtain an embedded
resolution. By iterating this argument, we can conclude that
studying the cluster C9 is enough to know the poles of the
topological zeta function associated to the blowing up of the
clusters C11. Neither we have to consider the cluster C4 and the
cluster C6.} \hfill $\square$ \end{remark}
 ${}$  ${}$ \begin{center} \textsc{7. The monodromy conjecture for candidate poles of order $1$}
\end{center} ${}$ \\
For the sake of completeness, we recall the short proof of the next
lemma (see also \cite{CRAS}). Recall that, given a candidate pole
$-\nu_j/N_j=a/b$ with $a$ and $b$ coprime, $J_b$ then denotes the
subset of indices $\{ 1 \leq i \leq r \mbox{ $|$ } b \mbox{ divides
} N_i \}$.
\begin{lemma}\label{lemmahulp}
Let $\chi(E_t^{\circ}) < 0$ such that we are in the situation
\begin{center}
\begin{pspicture}(-4,-0.5)(4,0.7)
\psline{-,linestyle=dashed}(-4,0)(-3,0) \psline{-}(-3,0)(-2,0)
\psline{-}(-2,0)(-1,0)
\psline{-,linestyle=dashed}(-0.6,0)(0.2,0)\psline{-}(-1,0)(-0.6,0)
\psline{-}(0.2,0)(1.2,0)\psline{-,linestyle=dashed}(1.3,0)(1.9,0)
\psline{-}(1.9,0)(3.3,0) \psline{-,linestyle=dashed}(3.3,0)(4,0)
 \psdot(-3,0)\psdot(-2,0)\psdot(-1,0)\psdot(3.5,0)\psdot(2.5,0)\psdot(0.7,0)\rput(-3,0.3){\footnotesize{$m_i$}}\rput(-2,0.3){\footnotesize{$m_i$}}\rput(-1,0.3){\footnotesize{$m_i$}}
 \rput(0.7,0.3){\footnotesize{$m_i$}}
\rput(2.5,0.30){\footnotesize{$m_i$}}\rput(3.6,0.34){m'}
 \rput(-2.5,0.2){\footnotesize{$3$}}\rput(-1.5,0.2){\footnotesize{$3$}}\rput(-0.2,0.2){\footnotesize{$3$}}\rput(1.5,0.2){\footnotesize{$3$}}
\rput(3,0.2){\footnotesize{$3$}}
\rput(-3,-0.3){\footnotesize{$Q_t$}}\rput(-2,-0.3){\footnotesize{$Q_{t+1}$}}\rput(-1,-0.3){\footnotesize{$Q_{t+2}$}}\rput(0.7,-0.3){\footnotesize{$Q_{j}$}}
\rput(2.5,-0.3){\footnotesize{$Q_l$}}
\rput(3.5,-0.3){\footnotesize{$Q_{l+1}$}} \rput(6.3,0){$(11)$}
\end{pspicture}
\end{center}
where $Q_t$ is the point in the chain with the lowest level for
which an edge with label $3$ is leaving and where $Q_l$ is the point
in this chain with the highest level for which its multiplicity is
equal to $m_i$.
\begin{enumerate}
\item If a set $J_b$ contains the index $t$, then it also
contains the indices in $\{t+1,\cdots,l\}$.
\item If $\frac{\nu_l}{N_l}=\frac{c}{d}$ with $c$ and $d$ coprime, then $t
\notin J_d$.
\end{enumerate}
\end{lemma}
\emph{Proof.} \quad If we denote the numerical data of $E_t$ by
$(\nu,N)$, then, independently of the number of points $Q_s$ for
which $t \rightarrow s$, one easily computes that the numerical data
for $i \in \{t+1,\cdots,l\}$ are
\[E_i((i-t+1)\nu-(i-t),(i-t+1)N).\] Now the first assertion follows
immediately. \\ To see the second claim, suppose that $t \in J_d$.
Then $d \mbox{ $|$ } N$ which implies that
\[l-t+1 | (l-t+1)\nu - (l-t).\]
This contradiction closes the proof. \hfill $\blacksquare$
\\ \\ We can now prove one of the most important properties concerning the surfaces we study.
(In \cite{CRAS} we proved this result for a more restricted class of
surfaces.)
\begin{theorem} \label{nietslecht}
 If $\chi(E_j^{\circ}) > 0$, then $e^{-2\pi i
\frac{\nu_j}{N_j}}$ is an eigenvalue of monodromy of $f$.
 \end{theorem}
\emph{Proof.} \quad Suppose that $E_j$ is an exceptional component
for which $\chi(E_j^{\circ}) > 0$. To prove that $e^{-2\pi i
\nu_j/N_j}$ is an eigenvalue of monodromy of $f$, we show that
$e^{-2\pi i \nu_j/N_j}$ is a pole of $\zeta_f$. We write $\nu_j/N_j$
as $a/b$ with $a$ and $b$ coprime. If $J_b$ does not contain an
index $t$ for which $\chi(E_t^{\circ}) < 0$, then there is nothing
to verify. So suppose now that $\chi(E_t^{\circ}) < 0$ and that $t
\in J_b$. From Lemma \ref{lemmahulp} it follows that $E_j \neq E_l$
and that $l \in J_b$. We will show that
$\chi(E_t^{\circ})+\chi(E_l^{\circ}) \geq 0$. The configuration in
$E_t \cong \mathbb{P}^2$ is as follows:
\begin{center}
\includegraphics{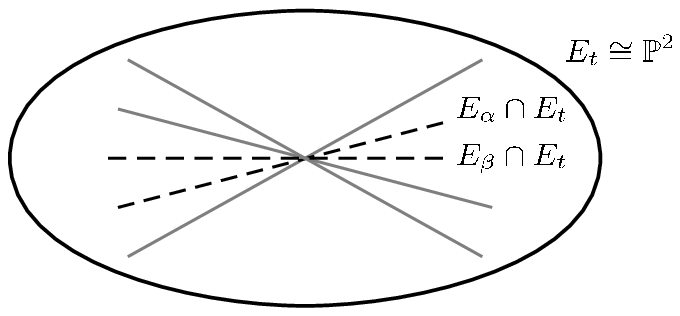}
\end{center}
\begin{enumerate}
\item If $Q_t$ is the origin of the constellation, then
$\chi(E_t^{\circ})=2-m_i$. For $\chi({E_l}^{\circ})$ we find that
$\chi(E_l^{\circ})=D-R+T$ with $D \geq m_i^2-m_im'$ and $R \leq
2m_i-2m'$. We get
\begin{eqnarray*}
\chi(E_t^{\circ})+\chi(E_l^{\circ}) & \geq &
2-m_i+m_i^2-m_im'-2m_i+2m'\\
& = & 2+2m'+m_i(m_i-m'-2).
\end{eqnarray*}
If $m_i-m' \geq 2$, then $\chi(E_t^{\circ})+\chi(E_l^{\circ}) > 0$.
If $m_i-m'=1$, then $\chi(E_t^{\circ})+\chi(E_l^{\circ}) \geq
2+2(m_i-1)-2m_i=0$. One can even check that also here
$\chi(E_t^{\circ})+\chi(E_l^{\circ}) > 0$. Hence we always have
$\chi(E_t^{\circ})+\chi(E_l^{\circ}) > 0$.
\item If there is exactly one point, say $Q_{\alpha}$, for which $t
\rightarrow \alpha$, then $\chi(E_t^{\circ}) = 1-m_i$. For
$\chi(E_l^{\circ})$ we find that $\chi(E_l^{\circ})=D-R+T$ with $D
\geq m_i^2-m_im'$ and $R \leq m_i-m'$ and we obtain
\begin{eqnarray*}
\chi(E_t^{\circ})+\chi(E_l^{\circ}) & \geq &
1-m_i+m_i^2-m_im'-m_i+m'\\
& = & 1+m'+m_i(m_i-m'-2).
\end{eqnarray*}
If $m_i-m' \geq 2$, then $\chi(E_t^{\circ})+\chi(E_l^{\circ}) > 0$.
If $m_i-m'=1$, we get $\chi(E_t^{\circ})+\chi(E_l^{\circ}) \geq 0$.
\item Finally, if there exist two points, say $Q_{\alpha}$ and
$Q_{\beta}$, for which $t \rightarrow \alpha$ and $t \rightarrow
\beta$, then $\chi(E_t^{\circ}) = -m_i$. In this case $R=0$ and we
get
\begin{eqnarray*}
\chi(E_t^{\circ})+\chi(E_l^{\circ}) & \geq &
-m_i+m_i^2-m_im'\\
& = & m_i(m_i-m'-1)\\
& \geq & 0.
\end{eqnarray*}
\end{enumerate}
This study permits us to conclude that $\sum_{i \in J_b}
\chi(E_i^{\circ})> 0$. Hence, $e^{-2\pi i \frac{\nu_j}{N_j}}$ is an
eigenvalue of monodromy of $f$.
 \vspace*{-1cm}
\[\]\hfill $\blacksquare$ ${}$
\\ \\ ${}$In the general case of surfaces it can happen that positive $\chi(E_j^{\circ})$ does not imply
that $e^{-2 \pi i \nu_j/N_j}$ is an eigenvalue of monodromy of $f$.
\begin{corollary}
If $-\nu_j/N_j$ is a candidate pole of $Z_{top,f}$ of order $1$ that
is a pole, then $e^{-2 \pi i \nu_j/N_j}$ is an eigenvalue of
monodromy of $f$.
\end{corollary}
\emph{Proof.} \quad In \cite{Veys3} it is shown that then there
exists an exceptional component $E_k$ for which $\nu_k / N_k = \nu_j
/ N_j$ and $\chi(E_k^{\circ})> 0$. The result follows now
immediately from Theorem \ref{nietslecht}. \hfill $\blacksquare$
\\ \\Actually the second author shows in \cite{Veys3} in particular that if $E_j$ is created by blowing up a point and if $\chi(E_j^{\circ})<
0$, then the contribution of $E_j$ to the residue of $-\nu_j /N_j$
for $Z_{top,f}$ is equal to $0$. In this particular setting, this is
a consequence of Proposition \ref{propnondeg}. \\ \indent We first
recall the notion for a polynomial to be \emph{nondegenerate with
respect to its Newton polyhedron}. Let $f \in
\mathbb{C}[x_1,\cdots,x_d]$ be a non-constant polynomial vanishing
in the origin. Write $\underline{x}^{\underline{k}}
:=x_1^{k_1}\cdots x_d^{k_d}$ and $f:=\sum_{\underline{k} \in
\mathbb{N}^d}c_{\underline{k}} \underline{x}^{\underline{k}}$. The
\emph{support of $f$} is supp$(f):=\{\underline{k} \in \mathbb{N}^d
\mbox{ $|$ } c_{\underline{k}} \neq 0\}$. %The \emph{global Newton
%polyhedron} $\Gamma_{gl}$ of $f$ is the convex hull of supp$(f)$ and
The \emph{Newton polyhedron} $\Gamma$ of $f$ is the convex hull of
supp$(f) + \mathbb{R}^d_{\geq 0}$. For a face $\tau$ of $\Gamma$ we
write $f_{\tau}:=\sum_{\underline{k} \in \tau}c_{\underline{k}}
\underline{x}^{\underline{k}}$. A polynomial $f$ is called
\emph{nondegenerate with respect to $\Gamma$} if for every compact
face $\tau$ of $\Gamma$, the polynomials $f_{\tau}$ and $\partial
f_{\tau}/\partial x_i$ have no common zeroes in $(\mathbb{C}^*)^d$,
$1 \leq i \leq d$.\\${}$ \vspace*{-0.2cm}
\begin{proposition} \label{propnondeg}
Every hypersurface that is general with respect to some
$3$-dimensional toric idealistic cluster is nondegenerate with
respect to its Newton polyhedron.
\end{proposition}
%${}$\\
\emph{Proof.} \quad Let $\mathcal{A}=(\mathcal{C},\underline{m})$ be
a toric idealistic cluster such that $f$ is general with respect to
$\mathcal{A}$. Suppose that $f$ is degenerate with respect to
$\mathcal{N}(f)$.\\ Let $\tau$ be a compact face of $\mathcal{N}(f)$
for which there exists a point $p \in (\mathbb{C}^*)^3$ such that
$f_{\tau}(p)=\partial f_{\tau}/\partial x(p)=
\partial f_{\tau}/\partial y(p)=\partial f_{\tau}/\partial z(p)=0$.
\\ If $\tau$ is a facet, then $\tau$ corresponds to some exceptional
irreducible component created by the blowing up of the
constellation, say to $E_i$. More specifically, the strict transform
of $f_{\tau}$ is equal to $E_0 \cap E_i$. As $p$ is not an orbit, it
follows that there exists a point in which $E_0 \cap E_i$ does not
have normal crossings and that is not an orbit. If the dimension of
$\tau$ is one and if $\tau$ is the intersection of two compact
facets, then analogously we have that there exist two irreducible
exceptional components $E_i$ and $E_j$ such that $E_0 \cap E_i \cap
E_j$ does not have normal crossings in a point that is not an orbit.
Remains the case that $\tau$ is the intersection of a compact facet
and a coordinate plane. Suppose that that compact facet corresponds
to $E_i$ and that the coordinate plane is given by $\{x=0\}$. Again
we get that then $E_0 \cap E_i$ does not have normal crossings in a
point that is not an orbit. Indeed, if $E_i$ has equation $y=0$ in
some affine chart, then there is a point $(0,0,p_z)$ with $p_z \neq
0$ in which there are no normal crossings. \hfill $\blacksquare$
\\ ${}$ \vspace*{-0.2cm}
\begin{corollary}\label{corolneg}
If $\chi(E_j^{\circ})< 0$, then the contribution of $E_j$ to the
residue of $-\nu_j /N_j$ for $Z_{top,f}$ is equal to $0$.
\end{corollary}
\emph{Proof.} \quad Denef and Loeser show in \cite{DenefLoeser1}
that the poles of $Z_{top,f}$ are of the form $-\nu(a)/N(a)$ where
$a$ is orthogonal to a facet of $\mathcal{N}(f)$. The compact facets
of $\mathcal{N}(f)$ correspond to the Rees valuations of the
complete ideal of hypersurfaces that pass through the points of the
constellation with at least the given multiplicity. The result now
follows from Proposition \ref{propnondeg} and Equation $(2)$ in
Section 2.3. Indeed, if $\chi(E_i^{\circ}) < 0$, then $m_i^2 =
\sum_{j \rightarrow i} m_j^2$. \hfill $\blacksquare$
%\begin{example}
%The surface with equation $f:=x^7 + y^5 + z^5 + x^3 y^2 + x^4 z^2 =
%0$ has an embedded resolution given by the constellation
%\begin{center}
%\begin{pspicture}(-4,-0.5)(4,0.7)
%\psline{-}(-3,0)(-2,0) \psline{-}(-2,0)(-1,0) \psline{-}(-1,0)(-0,0)
%\psline{-}(0,0)(1,0)
% \psdot(-3,0)\psdot(-2,0)\psdot(-1,0)\psdot(0,0)\psdot(1,0)\rput(-3,0.3){5}\rput(-2,0.3){2}\rput(-1,0.3){2}
% \rput(0,0.3){1}
% \rput(1,0.3){1} \rput(-2.5,0.2){\footnotesize{$1$}}\rput(-1.5,0.2){\footnotesize{$3$}}\rput(-0.5,0.2){\footnotesize{$3$}}\rput(0.5,0.2){\footnotesize{$2$}}
%\rput(-3,-0.3){\footnotesize{$Q_1$}}\rput(-2,-0.3){\footnotesize{$Q_{2}$}}\rput(-1,-0.3){\footnotesize{$Q_{3}$}}\rput(0,-0.3){\footnotesize{$Q_{4}$}}
%\rput(1,-0.3){\footnotesize{$Q_{5}$}}
% \end{pspicture}
%\end{center}
%For the $\chi({E_j}^{\circ})$ we find
%\[\chi(\overset{\circ}{E_1})=9, \quad \chi(\overset{\circ}{E_2})=-1,
%\quad \chi(\overset{\circ}{E_3})=1, \quad
%\chi(\overset{\circ}{E_4})= -1, \quad
%\chi(\overset{\circ}{E_5})=1.\] The numerical data are
%\[E_1(5,3), \quad E_2(7,5), \quad E_3(14,9), \quad E_4(20,13), \quad E_5(40,25).\]
%Our results show that $E_1, E_3$ and $E_5$ give rise to an
%eigenvalue of monodromy of $f$. \\ For the topological and the
%monodromy zeta function we get
%\[Z_{top,f}(s)=\frac{9(106s^2+107s+25)}{5(s+1)(14s+9)(8s+5)} \] and
%\[\zeta_f(t)=\frac{(1-t^7)(1-t^{20})}{(1-t^5)^9 (1-t^{14})(1-t^{40})}. \]
% \hfill $\square$
%\end{example}
\\ \\
Although the surfaces that we work with are all nondegenerate with
respect to their Newton polyhedron, our proof covers many new cases.
We recall the numerical conditions that the nondegenerate
polynomials should satisfy in the proof of the monodromy conjecture
that Loeser gave for them. Suppose that the blowing ups of $Q_i$ and
$Q_j$ give rise to Rees valuations and thus to facets $F_i$ and
$F_j$ of the Newton polyhedron. Suppose that their equations are
\begin{eqnarray*}
a_1(F_i)x_1+a_2(F_i)x_2+a_3(F_i)x_3 & =  & N_i \\
a_1(F_j)x_1+a_2(F_j)x_2+a_3(F_j)x_3 & =  & N_j
\end{eqnarray*}
and that these faces have a non-empty intersection. Let $a_{ij}$ be
the greatest common divisor of the determinants of the $2 \times
2$-matrices in the matrix
\[\left(\begin{array}{ccc}
a_1(F_i) & a_2(F_i) & a_3(F_i) \\
a_1(F_j) & a_2(F_j) & a_3(F_j) \end{array}\right).\] Then to be
covered by the proof of Loeser, it should hold that
\[\frac{\nu_i - \frac{\nu_j}{N_j}N_i}{a_{ij}} \notin \mathbb{Z} \qquad and \qquad \nu_i/N_i \notin \mathbb{Z}.\]
Already very simple toric clusters, such as for example the blowing
up of two points $Q_1$ and $Q_2$ with multiplicity $m_1=6$ and
$m_2=2$, do not satisfy these conditions. Also candidate poles of
order at least $2$ are not included. \newpage  ${}$ \begin{center}
\textsc{8. The monodromy conjecture for candidate poles of order $2$
or $3$}
\end{center} ${}$\\
Let us now study when the topological zeta function can have a
candidate pole of order at least $2$. Suppose a $3$-dimensional
toric idealistic cluster is given and suppose that the blowing up of
the cluster provides an embedded resolution for the hypersurface
$\{f=0\}$. %Its zeta function of monodromy is equal to
%\[\zeta_f(t)=\prod_{j=1}^r
%(1-t^{N_j})^{-\chi(E_j^{\circ})}.\]
Let $s$ be a candidate pole of
order at least $2$ of the topological zeta function associated to
$f$, say $s=-\nu_i/N_i=-\nu_j/N_j$, $1 \leq i,j \leq r$. We write
$s$ as $a/b$ such that $a$ and $b$ are coprime. If $J_b$ is the set
$\{j \in \{1,\cdots,r\} \mbox{ $|$
    } \mbox{ }  b \mbox{ divides } N_j \}$, then we study when $\sum_{j \in J_b} \chi(E_j^{\circ}) =
    0$. Recall that $e^{2 i \pi s}$ is not an eigenvalue of monodromy if this sum is $0$.\\
\indent As we are looking for candidate poles of order at least $2$
that are poles, it follows that $m_i^2$ should be different from
$\sum_{j \rightarrow i}m_j^2$ for one of the exceptional components
$E_i$ that yield that candidate pole. It follows now from Theorem
\ref{nietslecht} that we should study two cases. Firstly there are
the clusters with candidate poles of order at least two provided by
intersecting exceptional components $E_i$ and $E_j$ for which
$\chi(E_i^{\circ})=\chi(E_j^{\circ})=0$. Secondly we study the
clusters with candidate poles of order at least two provided by
intersecting exceptional components $E_i$ and $E_j$ for which
$\chi(E_i^{\circ})=0$ and $\chi(E_j^{\circ}) < 0$. In the following
subsections we proceed with the study of these cases. %\newpage  ${}$
\\ ${}$ \begin{center} \textsc{8.1. $\chi(E_i^{\circ})=\chi(E_j^{\circ})=0$}
\end{center}
\begin{proposition}
If $s_0=-\nu_i/N_i=-\nu_j/N_j$ is a candidate pole of $Z_{top,f}$ of
order at least $2$ that is a pole, and if
$\chi(E_i^{\circ})=\chi(E_j^{\circ})=0$, then $e^{2 \pi i s_0}$ is
an eigenvalue of monodromy of $f$.
\end{proposition}
\emph{Proof.} \quad Suppose that $j \rightarrow i$. We study the
possible combinations from List $3$.
\begin{itemize}
\item C8A and C9A: \\
we can only combine the cluster \quad
\begin{picture}(0,0)(3,2) \psdot(0.5,0) \psline{-}(0.5,0)(2.5,0)
\psdot(1.5,0) \psdot(2.5,0)
\rput(0.5,0.3){\footnotesize{$2$}}\rput(1.5,0.3){\footnotesize{$1$}}
\rput(2.5,0.3){\footnotesize{$1$}}
\rput(1,0.2){\scriptsize{\emph{3}}}\rput(2,0.2){\scriptsize{\emph{2}}}
\end{picture}\\ \\
of the form C9A with a cluster of the form C8A and then we get:\\
\begin{tabular}{p{5cm}p{8cm}}
\begin{pspicture}(-1,-1)(3,1)
\psline{-,linestyle=dashed}(-1,0)(-0.5,0)
\psline{-}(-0.5,0)(0.5,0)\psdot(-0.5,0) \psdot(0.5,0)
\rput(0.45,0.3){\footnotesize{$Q_i$}} \psline{-}(0.5,0)(1.5,0.5)
\psline{-}(0.5,0)(1.5,-0.5) \psdot(1.5,0.5)
\psdot(1.5,-0.5)\rput(1,0.45){\scriptsize{\emph{1}}}\rput(1,-0.5){\scriptsize{\emph{2}}}
\psline{-}(1.5,0.5)(2.5,0.5)\psdot(2.5,0.5)\psline{-}(2.5,0.5)(3.5,0.5)\psdot(3.5,0.5)
\rput(2,0.7){\scriptsize{\emph{3}}}\rput(3,0.7){\scriptsize{\emph{2}}}
\psline{-}(1.5,-0.5)(2.5,-0.5)\psdot(2.5,-0.5)\psline{-}(2.5,-0.5)(3.5,-0.5)\psdot(3.5,-0.5)
\rput(2,-0.7){\scriptsize{\emph{3}}}\rput(3,-0.7){\scriptsize{\emph{3}}}
\rput(0.5,-0.3){\footnotesize{$3$}}\rput(1.5,0.7){\footnotesize{$2$}}\rput(2.5,0.7){\footnotesize{$1$}}\rput(3.5,0.7){\footnotesize{$1$}}
\rput(1.5,-0.8){\footnotesize{$1$}}
\rput(2.5,-0.8){\footnotesize{$1$}}\rput(3.5,-0.8){\footnotesize{$1$}}\rput(0,0.2){\scriptsize{\emph{1}}}\rput(1.5,0.15){\footnotesize{$Q_j$}}
\end{pspicture}
& \vspace*{-1.6cm} Suppose only label $1$ appears under $Q_i$.
\end{tabular}
\\
If not, the upper chain in C8A would not be a Euclidean cluster. We
can write that the numerical data of $E_i$ are equal to
$(2i+1,\sum_{l=1}^i m_l)$ and that the ones of $E_j$ are equal to
$(2i+3,\sum_{l=1}^i m_l + 2)$. Hence, if $E_i$ and $E_j$ give rise
to the same candidate pole, we should have
\[\frac{2i+1}{\sum_{l=1}^i m_l}=\frac{2i+3}{\sum_{l=1}^i m_l + 2}.\]
If this equality holds, then $2i+1=\sum_{l=1}^i m_l \geq 4(i-1)+3$
and then $i$ should be equal to $1$. This is a contradiction because
$Q_i$ is not the origin.
\item C8A and C9B: there are two possibilities.
\begin{enumerate}
\item \begin{tabular}{p{5cm}p{5.5cm}}
\begin{pspicture}(-1,-1)(3,1)
\psline{-,linestyle=dashed}(-1,0)(-0.5,0)
\psline{-}(-0.5,0)(0.5,0)\psdot(-0.5,0) \psdot(0.5,0)
\rput(0.45,0.3){\footnotesize{$Q_i$}} \psline{-}(0.5,0)(1.5,0.5)
\psline{-}(0.5,0)(1.5,-0.5) \psdot(1.5,0.5)
\psdot(1.5,-0.5)\rput(1,0.45){\scriptsize{\emph{1}}}\rput(1,-0.5){\scriptsize{\emph{2}}}
\psline{-}(1.5,0.5)(2.5,0.5)\psdot(2.5,0.5)\psline{-}(2.5,0.5)(3.5,0.5)\psdot(3.5,0.5)
\rput(2,0.7){\scriptsize{\emph{3}}}\rput(3,0.7){\scriptsize{\emph{3}}}
\psline{-}(1.5,-0.5)(2.5,-0.5)\psdot(2.5,-0.5)\psline{-}(2.5,-0.5)(3.5,-0.5)\psdot(3.5,-0.5)
\rput(2,-0.7){\scriptsize{\emph{3}}}\rput(3,-0.7){\scriptsize{\emph{1}}}
\rput(0.5,-0.3){\footnotesize{$3$}}\rput(1.5,0.7){\footnotesize{$1$}}\rput(2.5,0.7){\footnotesize{$1$}}\rput(3.5,0.7){\footnotesize{$1$}}
\rput(1.5,-0.8){\footnotesize{$2$}}\rput(2.5,-0.8){\footnotesize{$1$}}\rput(3.5,-0.8){\footnotesize{$1$}}\rput(0,0.2){\scriptsize{\emph{1}}}
\rput(1.5,-0.15){\footnotesize{$Q_j$}}
\end{pspicture}
& \vspace*{-1.6cm} Suppose only label $1$ appears under $Q_i$.
\end{tabular}
\\We can write that the numerical data of $E_i$ are equal to
$(2i+1,\sum_{l=1}^{i-1} m_l+3)$. The numerical data of $E_j$ are
then equal to $(4i+1,2\sum_{l=1}^{i-1} m_l+3+2)$. If these
exceptional components give rise to the same candidate pole, then we
find that $2i-2=\sum_{l=1}^{i-1} m_l \geq 5(i-1)$. This can only
hold when $i=1$ but $Q_i$ is not the origin.
\item
\begin{tabular}{p{7.5cm}p{3cm}}
\begin{pspicture}(-1,-1)(6.5,1)
\psline{-,linestyle=dashed}(-1,0)(-0.5,0)
\psline{-}(-0.5,0)(0.5,0)\psdot(-0.5,0) \psdot(0.5,0)
\rput(0.45,0.3){\footnotesize{$Q_i$}} \psline{-}(0.5,0)(1.5,0.5)
\psline{-}(0.5,0)(1.5,-0.5) \psdot(1.5,0.5)
\psdot(1.5,-0.5)\rput(1,0.45){\scriptsize{\emph{1}}}\rput(1,-0){\scriptsize{\emph{2}}}
\psline{-,linestyle=dashed}(1.5,0.5)(2.5,0.5)\psdot(2.5,0.5)\psline{-,linestyle=dashed}(2.5,0.5)(4.5,0.5)%\psdot(3.5,0.5)%\psline{-}(3.5,0.5)(4.5,0.5)
\psdot(4.5,0.5)
\rput(2,0.7){\scriptsize{\emph{3}}}\rput(3.5,0.7){\scriptsize{\emph{3}}}
\psline{-}(4.5,0.5)(5.5,0.5)\psdot(5.5,0.5)
\rput(5,0.7){\scriptsize{\emph{2}}}\psline{-}(5.5,0.5)(6.5,0.5)\psdot(6.5,0.5)
\rput(6,0.7){\scriptsize{\emph{3}}}
\psline{-}(1.5,-0.5)(2.5,-0.5)\psdot(2.5,-0.5)\psline{-,linestyle=dashed}(2.5,-0.5)(3,-0.5)
\rput(2,-0.7){\scriptsize{\emph{3}}}
\rput(1.3,-0.8){\footnotesize{$3n-1$}}\rput(1.5,0.7){\footnotesize{$3$}}\rput(2.5,0.7){\footnotesize{$3$}}%\rput(3.5,0.7){\footnotesize{$3$}}
\rput(4.5,0.7){\footnotesize{$2$}}
\rput(5.5,0.7){\footnotesize{$1$}}\rput(6.5,0.7){\footnotesize{$1$}}\rput(0,0.2){\scriptsize{\emph{1}}}\rput(4.5,0.15){\footnotesize{$Q_j$}}
\rput(0.3,-0.3){\footnotesize{$3n+2$}}
\end{pspicture}
& \vspace*{-1.6cm} Suppose only label $1$ appears under $Q_i$.
\end{tabular}\\
If the numerical data of $E_i$ are $(2i+1,\sum_{l=1}^{i} m_l)$ and
if there are $n \geq 1$ points with multiplicity $3$ between $Q_i$
and $Q_j$, then the numerical data of $E_j$ are
$(2(n+1)i+(2n+3),(n+1)\sum_{l=1}^{i} m_l+(3n+2))$. If $E_i$ and
$E_j$ give the same candidate pole, then one should have
\[6in+4i+3n+2 = (n+2)\sum_{l=1}^{i} m_l \geq (n+2)((i-1)(6n+1)+3n+2)\]
or
\[8n+2i+3n^2 \geq 7in + 6in^2.\]
As $i \geq 2$, this inequality can never hold and thus $E_i$ and
$E_j$ can not give rise to the same candidate pole.\end{enumerate}
\item C8B and C9B: again there are two possibilities.
\begin{enumerate}
\item
\begin{tabular}{p{5cm}p{5.5cm}}
\begin{pspicture}(0,-1)(3,1)
\psline{-,linestyle=dashed}(0,0)(0.5,0) \psdot(0.5,0)
\rput(0.45,0.3){\footnotesize{$Q_i$}} \psline{-}(0.5,0)(1.5,0.5)
\psline{-}(0.5,0)(1.5,-0.5) \psdot(1.5,0.5)
\psdot(1.5,-0.5)\rput(1,0.45){\scriptsize{\emph{1}}}\rput(1,-0.5){\scriptsize{\emph{2}}}
\psline{-}(1.5,0.5)(2.5,0.5)\psdot(2.5,0.5)\psline{-}(2.5,0.5)(3.5,0.5)\psdot(3.5,0.5)
\rput(2,0.7){\scriptsize{\emph{3}}}\rput(3,0.7){\scriptsize{\emph{2}}}
\psline{-}(1.5,-0.5)(2.5,-0.5)\psdot(2.5,-0.5)\psline{-}(2.5,-0.5)(3.5,-0.5)\psdot(3.5,-0.5)
\rput(2,-0.7){\scriptsize{\emph{3}}}\rput(3,-0.7){\scriptsize{\emph{3}}}
\rput(0.5,-0.3){\footnotesize{$3$}}\rput(1.5,0.7){\footnotesize{$2$}}\rput(2.5,0.7){\footnotesize{$1$}}\rput(3.5,0.7){\footnotesize{$1$}}
\rput(1.5,-0.8){\footnotesize{$1$}}
\rput(2.5,-0.8){\footnotesize{$1$}}\rput(3.5,-0.8){\footnotesize{$1$}}\rput(1.5,0.15){\footnotesize{$Q_j$}}
%\rput(1.5,0.15){\footnotesize{$Q_j$}}\rput(1.5,0.15){\footnotesize{$Q_j$}}
\end{pspicture}
& \vspace*{-1.6cm} Suppose exactly label $1$ and label $2$ appear
under $Q_i$.
\end{tabular}
\\
In this situation $E_i$ and $E_j$ can give rise to the same
candidate pole, as shown in the following example:
\\
\begin{pspicture}(-4.5,-1)(3,1)
\psdot(-1.5,0)\psline{-}(-1.5,0)(-0.5,0)\rput(-1,0.2){\scriptsize{\emph{2}}}
\rput(-1.5,0.25){\footnotesize{$14$}}\rput(-0.5,0.25){\footnotesize{$5$}}
\psline{-}(-0.5,0)(0.5,0)\psdot(-0.5,0) \psdot(0.5,0)
\rput(0.45,0.3){\footnotesize{$Q_i$}} \psline{-}(0.5,0)(1.5,0.5)
\psline{-}(0.5,0)(1.5,-0.5) \psdot(1.5,0.5)
\psdot(1.5,-0.5)\rput(1,0.45){\scriptsize{\emph{1}}}\rput(1,-0.5){\scriptsize{\emph{2}}}
\psline{-}(1.5,0.5)(2.5,0.5)\psdot(2.5,0.5)\psline{-}(2.5,0.5)(3.5,0.5)\psdot(3.5,0.5)
\rput(2,0.7){\scriptsize{\emph{3}}}\rput(3,0.7){\scriptsize{\emph{2}}}
\psline{-}(1.5,-0.5)(2.5,-0.5)\psdot(2.5,-0.5)\psline{-}(2.5,-0.5)(3.5,-0.5)\psdot(3.5,-0.5)
\rput(2,-0.7){\scriptsize{\emph{3}}}\rput(3,-0.7){\scriptsize{\emph{3}}}
\rput(0.5,-0.3){\footnotesize{$3$}}\rput(1.5,0.7){\footnotesize{$2$}}\rput(2.5,0.7){\footnotesize{$1$}}\rput(3.5,0.7){\footnotesize{$1$}}
\rput(1.5,-0.8){\footnotesize{$1$}}
\rput(2.5,-0.8){\footnotesize{$1$}}\rput(3.5,-0.8){\footnotesize{$1$}}\rput(0,0.2){\scriptsize{\emph{1}}}\rput(1.5,0.15){\footnotesize{$Q_j$}}
\rput(3.5,0.15){\footnotesize{$Q_k$}}
\rput(2.5,0.15){\footnotesize{$Q_l$}}
\end{pspicture}
\\
%The exceptional components created by the blowing up of this
%constellation give rise to the following numerical data:
%\[(3,14), \quad (5,19), \quad (9,36), \quad (13,52), \quad (25,103), \quad (47,192)\]
%\[(15,56), \quad (29,112), \quad (43,168).\]
We find $\nu_i/N_i=\nu_j/N_j=-1/4$ and
%5427980183808s^7+8272603186848 s^6+5281511149120
%s^5+1814909840788s^4+357047147776 s^3+39023442131 s^2+2038682248
%s+29420385
\footnotesize{
\[Z_{top}(s)=\frac{A(s)}
{9(14s+3)(192s+47)(168s+43)(19s+5)(s+1)(103s+25)(4s+1)},\]}
\normalsize with $A$ a polynomial in $s$. However, we have $N_k=192$
and thus also $k \in J_b$. As $\chi(E_k^{\circ})=1>0$, we can
conclude that $e^{-2 i \pi / 4 }$ is an eigenvalue of monodromy.
This phenomenon is
true in general as we will see now.\\
We call $Q_l:=Q_j(3)$ and $Q_k:=Q_l(2)$. We show that if
$\nu_i/N_i=\nu_j/N_j=a/b$ with $a$ and $b$ coprime, then $b \mid
N_k$. Let $Q_2$ be the point with the highest level under $Q_i$ for
which $Q_2(2)$ is a point of the constellation. Let $(\nu_2,N_2)$ be
the numerical data of the point $Q_2$. Then we have that
\begin{eqnarray*}
N_k & = & N_i+N_j+N_l+1
\\
& = & N_i+(N_i+N_2+2)+(N_i+(N_i+N_2+2)+N_2+1)+1
\\ & = & 4N_i + 3N_2+ 6 \\
& = & N_i + 3N_j.
\end{eqnarray*}
Since $b \mid N_i$ and $b\mid N_j$, also $b\mid N_k$. As
$\chi(E_k^{\circ})=1>0$ and $E_k$ does not play the role of $E_l$ in
cluster $(11)$, it follows by the proof of Theorem \ref{nietslecht}
that $e^{2 \pi i s_0}$ is always an eigenvalue of monodromy.
\item
\begin{tabular}{p{7.5cm}p{3cm}}
\begin{pspicture}(0,-1)(6.5,1)
\psline{-,linestyle=dashed}(0,0)(0.5,0) \psdot(0.5,0)
\rput(0.45,0.3){\footnotesize{$Q_i$}} \psline{-}(0.5,0)(1.5,0.5)
\psline{-}(0.5,0)(1.5,-0.5) \psdot(1.5,0.5)
\psdot(1.5,-0.5)\rput(1,0.45){\scriptsize{\emph{1}}}\rput(1,-0.5){\scriptsize{\emph{2}}}
\psline{-}(1.5,0.5)(2.5,0.5)\psdot(2.5,0.5)\psline{-,linestyle=dashed}(2.5,0.5)(3.5,0.5)\psdot(3.5,0.5)\psline{-}(3.5,0.5)(4.5,0.5)
\psdot(4.5,0.5)
\rput(2,0.7){\scriptsize{\emph{3}}}\rput(3,0.7){\scriptsize{\emph{3}}}
\rput(4,0.7){\scriptsize{\emph{3}}}\psline{-}(4.5,0.5)(5.5,0.5)\psdot(5.5,0.5)
\rput(5,0.7){\scriptsize{\emph{2}}}\psline{-}(5.5,0.5)(6.5,0.5)\psdot(6.5,0.5)
\rput(6,0.7){\scriptsize{\emph{3}}}
\psline{-}(1.5,-0.5)(2.5,-0.5)\psdot(2.5,-0.5)\psline{-,linestyle=dashed}(2.5,-0.5)(3,-0.5)
\rput(2,-0.7){\scriptsize{\emph{3}}}
%\rput(0.5,-0.3){\footnotesize{$3$}}
\rput(1.5,0.7){\footnotesize{$3$}}\rput(2.5,0.7){\footnotesize{$3$}}\rput(3.5,0.7){\footnotesize{$3$}}
\rput(4.5,0.7){\footnotesize{$2$}}
\rput(5.5,0.7){\footnotesize{$1$}}\rput(6.5,0.7){\footnotesize{$1$}}\rput(4.5,0.15){\footnotesize{$Q_j$}}
\rput(3.5,0.15){\footnotesize{$Q_3$}}\rput(5.5,0.15){\footnotesize{$Q_l$}}\rput(6.5,0.15){\footnotesize{$Q_k$}}

\end{pspicture}
& \vspace*{-1.6cm} Suppose exactly label $1$ and label $3$ appear
under $Q_i$.
\end{tabular}
\\ We call $Q_l:=Q_j(2)$ and $Q_k:=Q_l(3)$. Let $Q_3$ be the point
such that $Q_j=Q_3(3)$ and let its associated numerical data be
$(\nu_3,N_3)$. Then we get
\begin{eqnarray*}
N_k & = & N_i+N_j+N_l+1
\\
& = & N_i+(N_i+N_3+2)+(N_i+N_3+(N_i+N_3+2)+1)+1
\\ & = & 4N_i + 3N_3+ 6 \\
& = & N_i + 3N_j.
\end{eqnarray*}
Again we can conclude that $e^{2 \pi i s_0}$ is always an eigenvalue
of monodromy.
\end{enumerate}
\item C9A and C7:\\
\begin{tabular}{p{8cm}p{3.5cm}}
\begin{pspicture}(-1,-1)(7,1)
\psline{-}(-0.5,0)(0.5,0) \psdot(-0.5,0)
\psline{-,linestyle=dashed}(-1,0)(-0.5,0) \psdot(0.5,0)
\rput(0.5,0.3){\footnotesize{$Q_i$}}\rput(1.5,0.3){\footnotesize{$P$}}
\psline{-}(0.5,0)(1.5,0) \psdot(1.5,0) \psline{-}(1.5,0)(2.5,0)
\psdot(2.5,0)
\rput(1,0.2){\scriptsize{\emph{1}}}\rput(2,0.2){\scriptsize{\emph{2}}}
\rput(0.5,-0.3){\footnotesize{$m_i$}}
\rput(1.5,-0.3){\footnotesize{$m_i-1$}}
\rput(2.5,-0.3){\footnotesize{$1$}} \psline{-}(2.5,0)(3.5,0)
\psdot(3.5,0)\rput(3,0.2){\scriptsize{\emph{3}}}\rput(3.5,-0.3){\footnotesize{$1$}}
\psline{-,linestyle=dashed}(3.5,0)(5.5,0) \psdot(5.5,0)
\psline{-}(5.5,0)(6.5,0)
\psdot(6.5,0)\rput(6,0.2){\scriptsize{\emph{3}}}
\rput(6.5,-0.3){\footnotesize{$1$}}\rput(0,0.2){\scriptsize{\emph{3}}}
\end{pspicture}
& \vspace*{-1.6cm} Suppose that only label $3$ appears under $Q_i$.
\end{tabular}
\\ Let $P:=Q_i(1)$ and $Q_j$ be the point $P(2,3^l)$ with $l \in \{0,1,\cdots,m_i-3\}$.
If the numerical data of $Q_i$ are equal to
$(2i+1,\sum_{s=1}^i{m_s})$, then we have the following numerical
data corresponding to the points
\begin{eqnarray*}
P & : & (4i+1,2\sum_{s=1}^{i-1}{m_s}+2m_i-1)
\\  P(2) & : & (8i+1,4\sum_{s=1}^{i-1}{m_s}+3m_i)     \\
P(2,3^l) & : & (i(8+6l)+2l+1,(4+3l)\sum_{s=1}^{i-1}{m_s}+(3+3l)m_i).
\end{eqnarray*}
We check if there exists an $l \in \{0,1,\cdots,m_i-2\}$ such that
\[\frac{2i+1}{\sum_{s=1}^{i}m_s}=\frac{i(8+6l)+2l+1}{(4+3l)\sum_{s=1}^{i-1}{m_s}+(3+3l)m_i}.\]
If this equality holds, then
\begin{eqnarray*}
2im_i-lm_i-2m_i & =  & (l+3)\sum_{s=1}^{i-1}{m_s} \\
& \geq & (l+3)(i-1)(2m_i-1) \\
& = & 2ilm_i + 6im_i - 2lm_i - 6m_i -il -3i +l + 3.
\end{eqnarray*}
We rewrite this and we get
\[3(i-1)+il \geq (m_i(2i-1)+1)l+4m_i(i-1).\]
As $3 < 4m_i$ and $i < m_i(2i-1)+1$, we get a contradiction. We
conclude that $Q_i$ and $Q_j$ can not give rise to the same
candidate pole.
\item C9A and C9B:\\
\begin{tabular}{p{5cm}p{6cm}}
\begin{pspicture}(-1,-1)(7,1)
\psline{-}(-0.5,0)(0.5,0) \psdot(-0.5,0)
\psline{-,linestyle=dashed}(-1,0)(-0.5,0) \psdot(0.5,0)
\rput(0.5,0.3){\footnotesize{$Q_i$}} \psline{-}(0.5,0)(1.5,0)
\psdot(1.5,0)\rput(1.5,0.3){\footnotesize{$Q_j$}}
\psline{-}(1.5,0)(2.5,0) \psline{-}(2.5,0)(3.5,0)\psdot(2.5,0)
\psdot(3.5,0)\rput(1,0.2){\scriptsize{\emph{1}}}\rput(2,0.2){\scriptsize{\emph{2}}}
\rput(0.5,-0.3){\footnotesize{$3$}}
\rput(1.5,-0.3){\footnotesize{$2$}}
\rput(2.5,-0.3){\footnotesize{$1$}}
\rput(0,0.2){\scriptsize{\emph{3}}}\rput(3,0.2){\scriptsize{\emph{3}}}\rput(3.5,-0.3){\footnotesize{$1$}}
\end{pspicture}
& \vspace*{-1.6cm} Suppose that only label $3$ appears under $Q_i$.
\end{tabular}\\
If $E_i$ has numerical data $(2i+1,\sum_{s=1}^{i}{m_s})$, then $E_j$
has numerical data $(4i+1,2\sum_{s=1}^{i-1}{m_s}+5)$. If they give
rise to the same candidate pole, then one should have
\[2i-2 = \sum_{s=1}^{i-1}{m_s} \geq 5(i-1).\]
As $Q_i$ is not the origin, this inequality can never be fulfilled.
\item C9B and C7: there are two possibilities.
\begin{enumerate}
\item
\begin{tabular}{p{7cm}p{3.5cm}}
\begin{pspicture}(0,-1)(7,1)
\psline{-,linestyle=dashed}(0,0)(0.5,0) \psdot(0.5,0)
\rput(0.5,0.3){\footnotesize{$Q_i$}}\rput(1.5,0.3){\footnotesize{$P$}}
\psline{-}(0.5,0)(1.5,0) \psdot(1.5,0) \psline{-}(1.5,0)(2.5,0)
\psdot(2.5,0)
\rput(1,0.2){\scriptsize{\emph{1}}}\rput(2,0.2){\scriptsize{\emph{2}}}
\rput(0.5,-0.3){\footnotesize{$m_i$}}
\rput(1.5,-0.3){\footnotesize{$m_i-1$}}
\rput(2.5,-0.3){\footnotesize{$1$}} \psline{-}(2.5,0)(3.5,0)
\psdot(3.5,0)\rput(3,0.2){\scriptsize{\emph{3}}}\rput(3.5,-0.3){\footnotesize{$1$}}
\psline{-,linestyle=dashed}(3.5,0)(5.5,0) \psdot(5.5,0)
\psline{-}(5.5,0)(6.5,0)
\psdot(6.5,0)\rput(6,0.2){\scriptsize{\emph{3}}}
\rput(6.5,-0.3){\footnotesize{$1$}}
\end{pspicture}
& \vspace*{-1.6cm} Suppose that exactly label $2$ and label $3$
appear under $Q_i$.
\end{tabular}
\\ Let $P:=Q_i(1)$ with numerical data $(\nu_1,N_1)$ and $Q_j$ be the point $P(2,3^l)$ with $l \in \{0,1,\cdots,m_i-3\}$.
Let $Q_3$, resp. $Q_2$, be the point with the highest level such
that $i
> 3$, resp. $i > 2$, and such that $Q_3(3)$, resp. $Q_2(2)$, is a point of the constellation. We
denote its numerical data by $(\nu_3,N_3)$, resp. $(\nu_2,N_2)$.
Suppose now that $\nu_i/N_i=\nu_j/N_j=a/b$ with $a$ and $b$ coprime.
Let $Q_k:=Q_1(2,3^{k-1})$. We show that $b \mid N_k$ when $k
> j$.\\
We have that
\begin{eqnarray*}
N_i & = & N_3+N_2+m_i  \quad \mbox{and}\\
N_1 & = & 2N_2+2N_3+2m_i-1 = 2N_i-1.
\end{eqnarray*}
If $Q_j=Q_i(1,2)$, then $N_j=N_3+N_i+N_1+1$ and so
\begin{eqnarray*}
N_k & = & (k-1)N_i+(k-1)N_1+N_j+(k-1) \\
 & = & (k-1)N_i+(k-1)(2N_i-1) +N_j+(k-1) \\
 & = & 3(k-1)N_i+N_j
\end{eqnarray*}
and we can conclude that $b \mid N_k$.
\\ If $Q_j=Q_i(1,2,3^l)$ for $l \neq 0$, then $N_j=lN_i +lN_1  + l + N_3 + N_i + N_1 +
1=(l+1)N_i+(l+1)N_1+N_3+(l+1)=(l+1)N_i+(l+1)(2N_i-1)+N_3+(l+1)=3(l+1)N_i+N_3$
and so
\begin{eqnarray*}
N_k & = & (k-1)N_i+(k-1)N_1+(N_3+N_i+N_1+1)+(k-1) \\
 & = & (k-1)N_i+(k-1)(2N_i-1) +(N_3+N_i+2N_i)+(k-1) \\
 & = & 3kN_i+N_3.
\end{eqnarray*}
As $b\mid N_i$ and $b\mid N_j$, we have that also $b\mid N_3$ and so
$b\mid N_k$. As $\chi(E_k^{\circ})=1>0$ for $k=m_i-1$ and
$Q_{m_i-1}$ can not play the role of $Q_l$ in cluster $(11)$, it
follows that $e^{2 \pi i s_0}$ is an eigenvalue of monodromy.
\item
In the previous cluster $Q_j$ can also be $Q_1$, but then $m_i$
should be equal to $2$.\\
\begin{tabular}{p{5cm}p{5.5cm}}
\begin{pspicture}(0,-1)(7,1)
\psline{-,linestyle=dashed}(0,0)(0.5,0) \psdot(0.5,0)
\rput(0.5,0.3){\footnotesize{$Q_i$}} \psline{-}(0.5,0)(1.5,0)
\psdot(1.5,0) \psline{-}(1.5,0)(2.5,0) \psdot(2.5,0)
\rput(1,0.2){\scriptsize{\emph{1}}}\rput(2,0.2){\scriptsize{\emph{2}}}
\rput(0.5,-0.3){\footnotesize{$2$}}
\rput(1.5,-0.3){\footnotesize{$1$}}
\rput(2.5,-0.3){\footnotesize{$1$}}
\end{pspicture}
& \vspace*{-1.6cm} Suppose that exactly label $2$ and label $3$
appear under $Q_i$.
\end{tabular}\\
We then have that $\nu_i/N_i=\nu_1/N_1=(2\nu_i-1)/(2N_i-1)$ if and
only if $\nu_i/N_i=1$. As $1$ is always an eigenvalue of monodromy,
this cluster does not give any problem.
\end{enumerate}
\item C10 and C7: \\
this case is completely analogous to the combination C9B and C7.
\item C10 and C8B: \\
\\
\begin{tabular}{p{4cm}p{8.8cm}}
\begin{pspicture}(0,-1)(3,1)
\psline{-,linestyle=dashed}(0,0)(0.5,0) \psdot(0.5,0)
\rput(0.5,0.3){\footnotesize{$Q_i$}}
\rput(1.5,0.3){\footnotesize{$P$}}\psline{-}(0.5,0)(1.5,0)
\psdot(1.5,0) \psline{-}(1.5,0)(2.5,0.5) \psdot(2.5,0.5)
\psline{-}(1.5,0)(2.5,-0.5) \psdot(2.5,-0.5)
\rput(1,0.2){\scriptsize{\emph{1}}}\rput(2,0.45){\scriptsize{\emph{2}}}\rput(2,-0.5){\scriptsize{\emph{3}}}
\psline{-}(2.5,0.5)(3.5,0.5) \psline{-}(2.5,-0.5)(3.5,-0.5)
\psdot(3.5,0.5)\psdot(3.5,-0.5) \rput(3.5,-0.8){\footnotesize{$1$}}
\rput(3.5,0.8){\footnotesize{$1$}}
\rput(2.5,-0.8){\footnotesize{$1$}}\rput(2.5,0.8){\footnotesize{$1$}}
\rput(0.5,-0.3){\footnotesize{$3$}}
\rput(1.5,-0.3){\footnotesize{$2$}}
\rput(3,-0.7){\scriptsize{\emph{1}}}
\rput(3,0.7){\scriptsize{\emph{1}}}
\end{pspicture}
& \vspace*{-1.6cm} Only label $2$ and label $3$ appear under $Q_i$.
\end{tabular}
\\
Let $Q_3$, resp. $Q_2$, be the point with the highest level such
that $i
> 3$, resp. $i > 2$, and such that $Q_3(3)$, resp. $Q_2(2)$, is a point of the constellation. We
denote its numerical data by $(\nu_3,N_3)$, resp. $(\nu_2,N_2)$.
Then we have that
\begin{eqnarray*}
N_i & = & N_2 + N_2 + 3 \\
N_j & = & N_i + N_2 + N_3 + 2 = 2 N_i - 1 \\
\nu_i & = & \nu_2 + \nu_3 + 1 \\
\nu_j & = & \nu_i+\nu_2 + \nu_3=2 \nu_i - 1.
\end{eqnarray*}
Hence if $\nu_i/N_i = \nu_j / N_j$, then $-\nu_i/N_i=-1$ and $1$ is
always an eigenvalue of monodromy.
\end{itemize} ${}$
\hfill $\blacksquare$ \\ ${}$ \begin{center} \textsc{8.2.
$\chi(E_i^{\circ})=0$ and $\chi(E_j^{\circ})<0$}
\end{center}
\begin{proposition}
If $s_0=-\nu_i/N_i=-\nu_j/N_j$ is a candidate pole of $Z_{top,f}$ of
order at least $2$ that is a pole, and if $\chi(E_i^{\circ})=0$ and
$\chi(E_j^{\circ})<0$, then $e^{2 \pi i s_0}$ is an eigenvalue of
monodromy of $f$.
\end{proposition}
\emph{Proof.} \quad We take List $2$ and List $3$ and we look for
the combinations that are possible to obtain $\sum_{k \in
J_b}\chi(E_k^{\circ})=0$. Recall that we proved in Theorem
\ref{nietslecht} that $\sum_{k \in J_b}\chi(E_k^{\circ})=0$ implies
that the value of $m'$ in cluster $(11)$ should be equal to $m_i-1$.
The only possible combination where at least $\nu_i$ or at least
$\nu_j$ is Rees, is the following one.
\begin{itemize}
\item C9A and C3:\\
\begin{tabular}{p{5cm}p{6cm}}
\begin{pspicture}(-1,-1)(4,1)
\psline{-}(-0.5,0)(0.5,0) \psdot(-0.5,0)
\psline{-,linestyle=dashed}(-1,0)(-0.5,0) \psdot(0.5,0)
\rput(0.5,0.3){\footnotesize{$Q_i$}} \psline{-}(0.5,0)(1.5,0)
\psdot(1.5,0)\rput(1.5,0.3){\footnotesize{$Q_j$}}
\psline{-}(1.5,0)(2.5,0) \psdot(2.5,0)
\rput(1,0.2){\scriptsize{\emph{1}}}\rput(2,0.2){\scriptsize{\emph{2}}}
\rput(0.5,-0.3){\footnotesize{$2$}}
\rput(1.5,-0.3){\footnotesize{$1$}}
\rput(2.5,-0.3){\footnotesize{$1$}}
\rput(0,0.2){\scriptsize{\emph{3}}}
\end{pspicture}
& \vspace*{-1.6cm} Suppose that only label $3$ appears under $Q_i$.
\end{tabular}\\
If the numerical data of $E_i$ are equal to $(2i+1,\sum_{s=1}^i
m_s)$, then the ones of $Q_j$ are equal to $(4i+1,2\sum_{s=1}^{i-1}
m_s+2+1)$. If $E_i$ and $E_j$ give rise to the same candidate pole,
then one should have
\[2i-1 = \sum_{s=1}^{i-1}m_s \geq 3(i-1)\]
which can only be true if $i=2$ and if the multiplicity of the
origin is $3$. Then we have the cluster
\\
\begin{pspicture}(-5,-1)(4,1)
\psline{-}(-0.5,0)(0.5,0) \psdot(-0.5,0) \psdot(0.5,0)
\rput(0.5,0.3){\footnotesize{$Q_i$}} \psline{-}(0.5,0)(1.5,0)
\psdot(1.5,0)\rput(1.5,0.3){\footnotesize{$Q_j$}}
\psline{-}(1.5,0)(2.5,0) \psdot(2.5,0)
\rput(1,0.2){\scriptsize{\emph{1}}}\rput(2,0.2){\scriptsize{\emph{2}}}
\rput(0.5,-0.3){\footnotesize{$2$}}
\rput(1.5,-0.3){\footnotesize{$1$}}
\rput(2.5,-0.3){\footnotesize{$1$}}
\rput(0,0.2){\scriptsize{\emph{3}}}\rput(-0.5,-0.3){\footnotesize{$3$}}
\end{pspicture}
\\
The candidate pole provided by $E_i$ and $E_j$ is then equal to
$-1$. Remember that $1$ is an eigenvalue of monodromy. \hfill
$\blacksquare$ \end{itemize}
 ${}$ \\ \\
Hence we can conclude with the following result.
\begin{theorem}
If $s_0$ is a candidate pole of $Z_{top,f}$ of order at least $2$
that is a pole, then $e^{2 \pi i s_0}$ is an eigenvalue of monodromy
of $f$.
\end{theorem}${}$ \\ ${}$
\begin{center} \textsc{9. The holomorphy conjecture}
\end{center} ${}$
\\To prove the holomorphy conjecture, we first prove the following
lemma. It gives us a set of orders of eigenvalues of monodromy.
\begin{lemma} \label{lemmahol}
If $\chi(E_j^{\circ}) > 0$, then $e^{2 \pi i /N_j}$ is an eigenvalue
of monodromy of $f$ at some point of the hypersurface $f=0$.
\end{lemma}
\emph{Proof.} \quad To prove that $e^{2 \pi i / N_j}$ is an
eigenvalue of monodromy, we will show that $\sum_{N_j\mid
N_i}\chi(E_i^{\circ}) \neq 0$. So suppose that $N_j \mid N_t$ and
$\chi(E_t^{\circ}) < 0$. Then we are in the situation
\begin{center}
\begin{pspicture}(-4,-0.5)(4,0.7)
\psline{-,linestyle=dashed}(-4,0)(-3,0) \psline{-}(-3,0)(-2,0)
\psline{-}(-2,0)(-1,0)
\psline{-,linestyle=dashed}(-0.6,0)(0.2,0)\psline{-}(-1,0)(-0.6,0)
\psline{-}(0.2,0)(1.2,0)\psline{-,linestyle=dashed}(1.3,0)(1.9,0)
\psline{-}(1.9,0)(3.3,0) \psline{-,linestyle=dashed}(3.3,0)(4,0)
 \psdot(-3,0)\psdot(-2,0)\psdot(-1,0)\psdot(3.5,0)\psdot(2.5,0)\psdot(0.7,0)\rput(-3,0.3){\footnotesize{$m_i$}}\rput(-2,0.3){\footnotesize{$m_i$}}\rput(-1,0.3){\footnotesize{$m_i$}}
 \rput(0.7,0.3){\footnotesize{$m_i$}}
\rput(2.5,0.30){\footnotesize{$m_i$}}\rput(3.6,0.34){m'}
 \rput(-2.5,0.2){\footnotesize{$3$}}\rput(-1.5,0.2){\footnotesize{$3$}}\rput(-0.2,0.2){\footnotesize{$3$}}\rput(1.5,0.2){\footnotesize{$3$}}
\rput(3,0.2){\footnotesize{$3$}}
\rput(-3,-0.3){\footnotesize{$Q_t$}}\rput(-2,-0.3){\footnotesize{$Q_{t+1}$}}\rput(-1,-0.3){\footnotesize{$Q_{t+2}$}}\rput(0.7,-0.3){\footnotesize{$Q_{j}$}}
\rput(2.5,-0.3){\footnotesize{$Q_l$}}
\rput(3.5,-0.3){\footnotesize{$Q_{l+1}$}} \rput(6.3,0){$(11)$}
\end{pspicture}
\end{center}
where $Q_t$ is the point in the chain with the lowest level for
which an edge with label $3$ is leaving and where $Q_l$ is the point
in this chain with the highest level for which its multiplicity is
equal to $m_i$.
\\In Lemma \ref{lemmahulp} we proved that then also $N_j \mid N_i$, for $i \in
\{j+1,\cdots,l\}$. As $N_l>N_t$, it follows that $N_l \nmid N_t$,
and hence $E_j \neq E_l$. In Theorem \ref{nietslecht} we then proved
that $\chi(E_t^{\circ})+\chi(E_l^{\circ}) \geq 0$, and thus we
obtain $\sum_{N_j\mid N_i}\chi(E_i^{\circ}) > 0$.
\vspace*{-1.5cm}\[\]\hfill $\blacksquare$
\begin{theorem}
If $r \in \mathbb{Z}_{>0}$ does not divide the order of any
eigenvalue of monodromy of $f$ at some point of the hypersurface
$f=0$, then $Z_{top,f}^{(r)}$ is holomorphic on $\mathbb{C}$.
\end{theorem}
\emph{Proof.} \quad Suppose that $Z_{top,f}^{(r)}$ is not
holomorphic, hence has a pole, say $s_0$. Let $E_i$ be an
exceptional component that gives rise to this pole of
$Z_{top,f}^{(r)}$ and let $(\nu_i,N_i)$ be its numerical data. If
$\chi(E_i^{\circ}) > 0$, then it follows from Lemma \ref{lemmahol}
that there is an eigenvalue of monodromy of order $N_i$. This
contradicts the given condition on $r$.
\\ If $\chi(E_i^{\circ}) < 0$, then we can set $E_i=E_t$ as in the
cluster above. Thus we also have $r \mid N_l$. However, as
$\chi(E_l^{\circ}) > 0$, it follows that $N_l$ is the order of an
eigenvalue of monodromy. \\
This implies that if $r \mid N_i$, then $\chi(E_i^{\circ})=0$. If
all these components are disjoint, then we get $Z_{top,f}^{(r)}=0$.
We may now suppose that at least two such components intersect each
other, and that at least one of them is Rees (it is shown in
\cite{DenefLoeser1} that only facets in the Newton polyhedron can
give rise to poles of $Z_{top,f}^{(r)}$). Then our cluster must
contain one of the following combinations of subclusters (see also
Section $8.1$.).
\begin{itemize}
\item C8A and C9A: we computed $N_j=N_i+2$, hence if $r\mid N_i$ and
$r\mid N_j$, then $r\mid2$. Set $Q_k:=Q_j(3,2)$, then $N_k=4N_i+6$
and $\chi(E_k^{\circ}) >0$. Lemma \ref{lemmahol} tells us that $N_k$
is the order of an eigenvalue of monodromy, which contradicts the
choice of $r$.
\item C8A and C9B: \begin{enumerate} \item we obtained $N_j=2N_i-1$. If $r\mid N_i$ and $r\mid N_j$,
then $r=1$, which divides the order of any eigenvalue of monodromy.
\item We had
$N_j=(n+1)N_i + (3n+2)$. Set $Q_k:=Q_j(2,3)$, then
$N_k=(3n+4)N_i+9n+6$. If $r$ divides $N_i$ and $N_j$, then it
follows that $r$ also divides $N_k$. As $\chi(E_k^{\circ})>0$, we
can conclude by Lemma \ref{lemmahol} that there is an eigenvalue of
order $N_k$. Again we get a contradiction.
\end{enumerate}
\item C8B and C9B: let $Q_k:=Q_j(3,2)$ as in that cluster in Section
$8.1$. We found $N_k=N_i+3N_j$. Analogously, we find that $E_i$ and
$E_j$ do not give rise to poles of $Z_{top,f}^{(r)}$, if $r\mid N_i$
and $r\mid N_j$.
\\ Also the other combination of C8B and C9B in Section $8.1$ gives
this contradiction.
\item C9A and C7: for $Q_j=P(2,3^l)$, we computed
$N_j=(4+3l)N_i-m_i$. So if $r\mid N_i$ and $r\mid N_j$, then $r\mid
m_i$. Let $Q_k:=P(2,3^{m_i-2})$ be the maximal point. Then
$N_k=(4+3k)N_i-m_i$, hence $r\mid N_k$, but as
$\chi(E_k^{\circ})>0$, we get a contradiction.
\item C9A and C9B: in this cluster we had $N_j=2N_i-1$, but then $r$
should be equal to $1$.
\item C9B and C7: again we can use the maximal point $Q_k:=P(2,3^{m_i-2})$. In Section
$8.1$ we saw already that $\chi(E_k^{\circ})>0$ and if $r$ divides
$N_i$ and $N_j$, that $r$ then also divides $N_k$.
\item C10 and C7: this case is exactly the same as the previous one.
\item C10 and C8B: we found that $N_j=2N_i-1$, thus it follows that
when $r$ divides $N_i$ and $N_j$, then $r=1$.
\end{itemize}
Hence, we find that $Z_{top,f}^{(r)}$ can neither have a pole coming
from an exceptional component for which $\chi(E_i^{\circ})=0$. This
ends the proof. \hfill $\blacksquare$
\\ ${}$
%\begin{center} \textsc{10. Remarks}
%\end{center} ${}$\\
%Notice first of all that in this context, the clusters where some
%$\chi(E_i^{\circ})$ is negative or zero are very rare. In
%particular, we find here that $\chi(E_i^{\circ}) < 0$ if and only if
%the configuration in $E_i \cong \mathbb{P}^2$ consists of (at least
%three) lines - possibly exceptional - that are all going through the
%same point. In the general case of surfaces, there exist much more
%configurations that yield a negative $\chi(E_i^{\circ})$. In
%\cite{Veysconfigurations} are given such examples. It also often
%happens that positive $\chi(E_j^{\circ})$ does not imply that
%$e^{-2\pi i \nu_j/N_j}$ is an eigenvalue of monodromy of $f$.
%\\
%Notice also that our surfaces never give rise to a pole of order $2$
%or $3$ if the concerning $\chi(E_{\cdot}^{\circ})$ are zero or
%negative.
\\ \\ Notice that if $r \mid N_i$ and $r \mid N_j$ with
$\chi(E_i^{\circ})=\chi(E_j^{\circ})=0$ and $E_i \cap E_j \neq
\emptyset$, then we found that $r=1$ or that there exists another
component $E_k$ with $r \mid N_k$ and $\chi(E_k^{\circ}) > 0$. For
general surfaces such a component $E_k$ does not necessarily exist.
\\ \\ \\
%\emph{Acknowledgement:} The first author would like to express her
%gratitude to Antonio Campillo for the valuable discussions.

\footnotesize{

\end{document}